\documentclass[a4paper,12pt]{article}
\usepackage{amsmath}
\usepackage{amscd}
\usepackage[xdvi]{graphicx}
\usepackage[dvips]{color}
\usepackage{times}
\usepackage{txfonts}

\newcommand{\bm}[1]{\mbox{\boldmath $#1$}}

\makeatletter

\renewcommand{\appendix}{\par
  \setcounter{section}{0}%
  \setcounter{subsection}{0}%
  \renewcommand{\thesection}{\@Alph\c@section}}

\@addtoreset{equation}{section}

\begin{document}
\begin{center}
\textbf{\LARGE{Periodic orbits and chaos in fast-slow systems with Bogdanov-Takens type fold points}}
\end{center}

\begin{center}
Faculty of Mathematics

Kyushu University, Fukuoka, 819-0395, Japan

Hayato CHIBA \footnote{E mail address : chiba@math.kyushu-u.ac.jp}
\end{center}
\begin{center}

Revised Mar 16 2010
\end{center}

\begin{center}
\textbf{Abstract}
\end{center}

The existence of stable periodic orbits and chaotic invariant sets of singularly perturbed problems of fast-slow type
having Bogdanov-Takens bifurcation points in its fast subsystem is proved by means of the geometric 
singular perturbation method and the blow-up method.
In particular, the blow-up method is effectively used for analyzing the flow near
the Bogdanov-Takens type fold point in order to show that a slow manifold
near the fold point is extended along the Boutroux's tritronqu\'{e}e
solution of the first Painlev\'{e} equation in the blow-up space.
\\[0.2cm]


\section{Introduction}

Let $(x_1, \cdots  ,x_n, y_1, \cdots ,y_m) \in \mathbf{R}^{n+m}$ be the Cartesian coordinates.
A system of singularly perturbed ordinary differential equations of the form
\begin{equation}
\left\{ \begin{array}{l}
\dot{x}_1 = f_1(x_1, \cdots  ,x_n, y_1, \cdots ,y_m, \varepsilon ), \\
\quad \vdots \\
\dot{x}_n = f_n(x_1, \cdots  ,x_n, y_1, \cdots ,y_m, \varepsilon ), \\ 
\dot{y}_1 = \varepsilon g_1(x_1, \cdots  ,x_n, y_1, \cdots ,y_m, \varepsilon ), \\
\quad \vdots \\
\dot{y}_m = \varepsilon g_m(x_1, \cdots  ,x_n, y_1, \cdots ,y_m, \varepsilon ),
\end{array} \right.
\label{1-1}
\end{equation}
is called a \textit{fast-slow system}, where the dot $(\,\dot{\,\,}\,)$
denotes the derivative with respect to time $t$, and where $\varepsilon >0$ is a small parameter.
Fast-slow systems are characterized by two different time scales, fast and slow time.
In other words, the dynamics consists of fast motions ($(x_1, \cdots  ,x_n)$ direction in the above system)
and slow motions ($(y_1 ,\cdots  ,y_m)$ direction).
This structure yields nonlinear phenomena such as a relaxation oscillation,
which is observed in many physical, chemical and biological problems.
See Grasman~\cite{Gra}, Hoppensteadt and Izhikevich~\cite{Hop} and references therein for applications of fast-slow systems.
To analyze the fast-slow system,
the unperturbed system (\textit{fast system}) of Eq.(\ref{1-1}) is defined to be
\begin{equation}
\left\{ \begin{array}{l}
\dot{x}_1 = f_1(x_1, \cdots  ,x_n, y_1, \cdots ,y_m, 0), \\
\quad \vdots \\
\dot{x}_n = f_n(x_1, \cdots  ,x_n, y_1, \cdots ,y_m, 0), \\ 
\dot{y}_1 = 0, \\
\quad \vdots \\
\dot{y}_m = 0.
\end{array} \right.
\label{1-2}
\end{equation}
The set of fixed points of the unperturbed system is called a \textit{critical manifold}, which is defined by
\begin{equation}
\mathcal{M} = \{ (x_1, \cdots  ,x_n, y_1, \cdots ,y_m) \in \mathbf{R}^{n+m} \, | \,
f_i(x_1, \cdots  ,x_n, y_1, \cdots ,y_m,0) = 0,\,\, i= 1, \cdots  ,n\}.
\label{1-3}
\end{equation}
Typically $\mathcal{M}$ is an $m$-dimensional manifold.
Fenichel~\cite{Fen} proved that if $\mathcal{M}$ is normally hyperbolic, then the original system (\ref{1-1})
with sufficiently small $\varepsilon >0$ has a locally invariant
manifold $\mathcal{M}_\varepsilon $ near $\mathcal{M}$, 
and that dynamics on $\mathcal{M}_\varepsilon $ is approximately given by the $m$-dimensional system
\begin{equation}
\left\{ \begin{array}{l}
\dot{y}_1 = \varepsilon g_1(x_1, \cdots  ,x_n, y_1, \cdots ,y_m, 0), \\
\quad \vdots \\
\dot{y}_m = \varepsilon g_m(x_1, \cdots  ,x_n, y_1, \cdots ,y_m, 0),
\end{array} \right.
\label{1-4}
\end{equation}
where $(x_1, \cdots  ,x_n, y_1, \cdots ,y_m) \in \mathbf{R}^{n+m}$ is restricted to the critical manifold $\mathcal{M}$.
The $\mathcal{M}_\varepsilon $ is diffeomorphic to $\mathcal{M}$ and called the \textit{slow manifold}.
The dynamics of (\ref{1-1}) approximately consists of the fast motion governed by (\ref{1-2})
and the slow motion governed by (\ref{1-4}).
His method for constructing an approximate flow is called the \textit{geometric singular perturbation method}.

However, if the critical manifold $\mathcal{M}$ has degenerate points $\bm{x}_0\in \mathcal{M}$ in the sense that
the Jacobian matrix $\partial \bm{f}/\partial \bm{x},\,\, \bm{f} = (f_1 , \cdots , f_n),\, \bm{x} = (x_1, \cdots  ,x_n)$
at $\bm{x}_0$ has eigenvalues on the imaginary axis, then $\mathcal{M}$ is not normally hyperbolic near the $\bm{x}_0$
and Fenichel's theory is no longer applicable.
The most common case is that $\partial \bm{f}/\partial \bm{x}$ has one zero-eigenvalue at $\bm{x}_0$ and the critical
manifold $\mathcal{M}$ is folded at the point (\textit{fold point}).
In this case, orbits on the slow manifold $\mathcal{M}_\varepsilon $ may jump and get away from $\mathcal{M}_\varepsilon $
in the vicinity of $\bm{x}_0$. As a result, the orbit repeatedly switches between fast motions and slow motions, and
complex dynamics such as a relaxation oscillation can occur.
See Mishchenko and Rozov~\cite{Mis} and Jones~\cite{Jon} for treatments of jump points and the existence of relaxation oscillations
based on the boundary layer technique and the geometric singular perturbation method.

The blow-up method was developed by Dumortier~\cite{Dum1} to investigate local flows near non-hyperbolic fixed points and it was
applied to singular perturbed problems by Dumortier and Roussarie~\cite{Dum2}.
The most typical example is the system of the form
\begin{equation}
\left\{ \begin{array}{l}
\dot{x} = -y + x^2, \\
\dot{y} = \varepsilon g(x,y),
\end{array} \right.
\label{2dim}
\end{equation}
where $(x,y) \in \mathbf{R}^2$.
The critical manifold is a graph of $y = x^2$ and the origin is the fold point,
at which the Jacobian matrix of the fast system has a zero-eigenvalue.
Indeed, the fast system $\dot{x} = -y + x^2$ undergoes a saddle-node bifurcation as $y$ varies.
To analyze this family of vector fields, the trivial equation $\dot{\varepsilon } = 0$ is attached as
\begin{equation}
\left\{ \begin{array}{l}
\dot{x} = -y + x^2, \\
\dot{y} = \varepsilon g(x,y),  \\
\dot{\varepsilon } = 0.
\end{array} \right.
\end{equation}
Then, the Jacobian matrix at the origin $(0,0,0)$ degenerates as
\begin{equation}
\left(
\begin{array}{@{\,}ccc@{\,}}
0 & -1 & 0 \\
0 & 0 & g(0,0) \\
0 & 0 & 0
\end{array}
\right)
\end{equation}
with the Jordan block.
The blow-up method is used to desingularize such singularities based on certain coordinate transformations.
The most simple case $g(0,0) \neq 0$ is deeply investigated by Krupa and Szmolyan et al.~\cite{Kru1, Gil}
with the aid of a geometric view point. Straightforward extensions to higher dimensional cases
are done by Szmolyan and Wechselberger~\cite{Szm2} for $n=1, m=2$ and by Mishchenko and Rozov~\cite{Mis} for any $n$ and $m$.
Under the assumptions that $\partial \bm{f}/\partial \bm{x}$ has only one zero-eigenvalue at a fold point
and that the slow dynamics (\ref{1-4}) has no fixed points near the fold point, they show that
in the blow-up space, the system is reduced to the Riccati equation $dx/dy = y-x^2$  
for any $n \geq 1$ and $m\geq 1$, and a certain special solution of the 
Riccati equation plays an important role to extend a slow manifold $\mathcal{M}_\varepsilon $ to a neighborhood of the fold point, which
guides jumping orbits.
It is to be noted that the classical work of Mishchenko and Rozov~\cite{Mis} is essentially equivalent to the blow-up method.

On the other hand, if the dynamics (\ref{1-4}) has fixed points on (a set of) fold points,
for example, if $g(0,0) = 0$ in Eq.(\ref{2dim}), then more complex phenomena such as canard explosion
can occur. Such situations are investigated by~\cite{Dum2, Kru1, Szm1, Kru3, Mae} by using the blow-up method.
For example, for Eq.(\ref{2dim}) with $g(0,0) = 0$, the original system is reduced to the system
$\dot{x} = -y + x^2,\, \dot{y} = x$ in the blow-up space.
If the dimension $m$ of slow direction is larger than $1$, there are many types of fixed points of (\ref{1-4})
and thus we need more hard analysis as is done in~\cite{Kru3}.

The fast system for Eq.(\ref{2dim}) undergoes a saddle-node bifurcation at the fold point.
Thus we call the fold point the \textit{saddle-node type fold point}.
The cases that fast systems undergo a transcritical bifurcation and a pitchfork bifurcation
are studied in~\cite{Kru2}.
It is shown that in the blow-up space, systems are reduced to the equations 
$dx/dy = x^2 - y^2 + \lambda $ and $dx/dy = xy - x^3$, respectively,
whose special solutions are used to construct slow manifolds near fold points.

Despite many works, behavior of flows near fold points at which the Jacobian matrix $\partial \bm{f}/\partial \bm{x}$ 
of the fast system has more than one zero-eigenvalues is not understood well.
The purpose of this article is to investigate a three dimensional fast-slow system of the form
\begin{equation}
\left\{ \begin{array}{l}
\dot{x} = f_1(x,y,z,\varepsilon ,\delta ), \\ 
\dot{y} = f_2(x,y,z,\varepsilon ,\delta ), \\
\dot{z} = \varepsilon g(x,y,z,\varepsilon ,\delta),
\end{array} \right.
\label{1-5}
\end{equation}
whose fast system has fold points with two zero-eigenvalues,
where $f_1, f_2, g$ are $C^\infty$ functions, $\varepsilon >0$ is a small parameter, and where $\delta >0$
is a small parameter which controls the strength of the stability of the critical manifold (see the assumption (C5) in Sec.2).
Note that the critical manifold
\begin{equation}
\mathcal{M} (\delta )= \{ (x,y,z) \in \mathbf{R}^3 \, | \, f_1(x,y,z,0,\delta ) = f_2(x,y,z,0,\delta ) = 0\}
\label{1-6}
\end{equation}
gives curves on $\mathbf{R}^3$ in general.
We consider the situation that at a fold point $(x_0, y_0, z_0) \in \mathbf{R}^3$ on $\mathcal{M}$,
the Jacobian matrix $\partial (f_1, f_2)/\partial (x, y)$ has two zero-eigenvalues with the Jordan block,
and the two dimensional unperturbed system (fast system) undergoes a Bogdanov-Takens bifurcation.
We call such a fold point the \textit{Bogdanov-Takens type fold point}.
For this system, we will show that the first Painlev\'{e} equation
\begin{eqnarray*}
\frac{d^2y}{dz^2} = y^2 -z 
\end{eqnarray*}
appears in the blow-up space and
plays an important role in the analysis of a local flow near the Bogdanov-Takens type fold points. 
This is in contrast with the fact that the Riccati equation appears in the case of saddle-node type fold points.
It is shown that in the blow-up space, the slow manifold is extended along one of the special solutions, the 
Boutroux's tritronqu\'{e}e solution~\cite{Bou, Jos}, of the first Painlev\'{e} equation.
One of the main results in this article is that a transition map of Eq.(\ref{1-5}) near the Bogdanov-Takens type fold point
is constructed, in which an asymptotic expansion and a pole of the Boutroux's tritronqu\'{e}e solution
are essentially used.
This result shows that the distance between a solution of (\ref{1-5}) near the Bogdanov-Takens type fold point
and a solution of its unperturbed system is of order $O(\varepsilon ^{4/5})$ as $\varepsilon \to 0$ (see Theorem 1 and Theorem 3.2),
while it is of $O(\varepsilon ^{2/3})$ for a saddle-node type fold point (see Mishchenko and Rozov~\cite{Mis}).

It is remarkable that all equations appeared in the blow-up space are related to the Painlev\'{e} theory.
For example, the equation $dx/dy = y - x^2$ obtained from the saddle-node type fold point
is transformed into the Airy equation $du/dy = uy$ by putting $x = (du/dy)/u$, which gives 
classical solutions of the second Painlev\'{e} equation.
The equation $dx/dy = x^2 - y^2 + \lambda $ obtained from the transcritical type fold point
is transformed into the Hermite equation 
\begin{eqnarray*}
\frac{d^2u}{dy^2} + 2y \frac{du}{dy} + (\lambda  + 1) u = 0
\end{eqnarray*}
by putting $x+y = -(du/dy)/u$, which gives classical solutions of the fourth Painlev\'{e} equation.
For other cases listed above, we also see that equations appeared in the blow-up space have the
Painlev\'{e} property~\cite{Con, Inc}; that is, all movable singularities (in the sense of the theory of ODEs on the complex plane)
are poles, not branch points and essential singularities.
This seems to be common for a wide class of fast-slow systems.
Painlev\'{e} equations have many good properties \cite{Con}.
For example, poles of solutions of Painlev\'{e} equations can be transformed into zeros of solutions of 
certain analytic systems by analytic transformations, which allow us to prove that 
the dominant part of the transition map near the Bogdanov-Takens type fold point 
is given by an analytic function describing a position of poles of the first Painlev\'{e} equation.

We also investigate global behavior of the system.
Under some assumptions, we will prove that there exists a stable periodic orbit (relaxation oscillation)
if $\varepsilon >0$ is sufficiently small for fixed $\delta $, 
and further that there exists a chaotic invariant set if $\delta >0$ is also small in comparison with small $\varepsilon $.
Roughly speaking, $\delta $ controls the strength of the stability of stable branches of the critical manifolds.
While chaotic attractors on $3$-dimensional fast-slow systems are reported by Guckenheimer, Wechselberger and Young~\cite{Guc}
in the case of $n=1, m=2$, our system is of $n=2, m=1$.
In the situation of~\cite{Guc}, the chaotic attractor arises according to the theory of H\'{e}non-like maps.
On the other hand, in our system, the mechanism of the onset of a chaotic invariant set is 
similar to that in Silnikov's works~\cite{Ovs, Sil1, Sil2},
in which the existence of a hyperbolic horseshoe is shown for a $3$-dimensional system 
which have a saddle-focus fixed point with a homoclinic orbit. See also Wiggins~\cite{Wig}.
Indeed, in our situation, the critical manifold $\mathcal{M}(\delta )$ plays a similar role to a saddle-focus fixed point
in the Silnikov's system.
Thus the proof of the existence of a relaxation oscillation in our system will be done in usual way:
the Poincar\'{e} return map proves to be contractive, while the proof of the existence of chaos is done in a similar way
to that of the Silnikov's system: as $\delta $ decreases, the Poincar\'{e} return map becomes non-contractive,
undergoes a cascade of bifurcations, and horseshoes are created.
When one want to prove the existence of a stable periodic orbit, it is sufficient to show that the image of the 
return map is exponentially small.
However, to prove the existence of a horseshoe, one has to show that the image of a rectangle under the return map
becomes a horseshoe-shaped (ring-shaped).
Thus our analysis for constructing the return map involves hard calculations, which can be avoided when proving
only a periodic orbit.

Our chaotic invariant set seems to be attracting as that in~\cite{Guc}, however, it remains unsolved.
See Homburg~\cite{Hom} for the proof of the existence of chaotic attractors in the Silnikov's system.

The results in the present article are used in \cite{Chi2} to investigate chaotic invariant sets on the Kuramoto model,
which is one of the most famous models to explain synchronization phenomena.
In \cite{Chi2}, it is shown that the Kuramoto model with appropriate assumptions can be reduced to a three dimensional fast-slow system
by using the renormalization group method \cite{Chi1}.

This paper is organized as follows.
In section 2, we give statements of our theorems on the existence of a periodic orbit and a chaotic invariant set.
An intuitive explanation of the theorems is also shown with an example.
In section 3, local analysis near the Bogdanov-Takens type fold point is given by means of the blow-up method.
Section 4 is devoted to global analysis, and proofs of main theorems are given.
Concluding remarks are included in section 5.


\section{Main results}

To obtain a local result and the existence of relaxation oscillations, the parameter $\delta $ in Eq.(\ref{1-5})
does not play a role. Thus we consider the system of the form
\begin{equation}
\left\{ \begin{array}{l}
\dot{x} = f_1(x,y,z,\varepsilon ), \\ 
\dot{y} = f_2(x,y,z,\varepsilon ), \\
\dot{z} = \varepsilon g(x,y,z,\varepsilon ),
\end{array} \right.
\label{local1}
\end{equation}
with $C^\infty$ functions $f_1, f_2, g : U \times I \to \mathbf{R}$, where $U \subset \mathbf{R}^3$ is an open domain in $\mathbf{R}^3$
and $I \subset \mathbf{R}$ is a small interval containing zero. The unperturbed system is given as
\begin{equation}
\left\{ \begin{array}{l}
\dot{x} = f_1(x,y,z,0 ), \\ 
\dot{y} = f_2(x,y,z,0 ), \\
\dot{z} = 0.
\end{array} \right.
\label{local2}
\end{equation}
Since $z$ is a constant, this system is regarded as a family of $2$-dimensional systems.
The critical manifold is the set of fixed point of (\ref{local2}) defined to be
\begin{equation}
\mathcal{M} = \{ (x,y,z) \in U \, | \, f_1(x,y,z,0 ) = f_2(x,y,z,0 ) = 0\}.
\label{local3}
\end{equation}
The reduced flow on the critical manifold is defined as
\begin{equation}
\dot{z} = \varepsilon g(x,y,z,0)|_{(x,y,z) \in \mathcal{M}}.
\label{local4}
\end{equation}
To investigate a Bogdanov-Takens type fold point, we make the following assumptions.
\\[0.2cm]
\textbf{(A1)}\, The critical manifold $\mathcal{M}$ has a smooth component $S^+ = S_a^+ \cup \{L^+\} \cup S_r^+$,
where $S_a^+$ consists of stable focus fixed points, $S_r^{+}$ consists of saddle fixed points,
and where $L^{+}$ is a fold point.
\\[0.2cm]
\textbf{(A2)} \, The $L^{+}$ is a Bogdanov-Takens type fold point; that is,
$L^{+}$ is a Bogdanov-Takens bifurcation point of the vector field $(f_1(x,y,z,0), f_2(x,y,z,0))$.
In particular, Eq.(\ref{local2}) has a cusp at $L^{+}$.
\\[0.2cm]
\textbf{(A3)} \, The reduced flow (\ref{local4}) on $S^+_a$ is directed toward the fold point $L^+$ and $g(L^+,0) \neq 0$.
\\[-0.2cm]

A few remarks are in order.
It is easy to see from (A1) that the Jacobian matrix $\partial (f_1, f_2)/\partial (x,y)$ has two zero
eigenvalues at $L^{+}$ since $S^{+}_r$ and $S^{+}_a$ are saddles and focuses, respectively.
Thus there exists a coordinate transformation $(x,y,z) \mapsto (X,Y,Z)$ defined near $L^+$ such that
$L^+$ is placed at the origin and Eq.(\ref{local2}) takes the following normal form
\begin{equation}
\left\{ \begin{array}{l}
\dot{X} = a_1(Z) + a_2(Z) Y^2 + a_3(Z)XY + O(X^3,X^2Y, XY^2,Y^3),  \\
\dot{Y} = b_1(Z) + b_2(Z)X + O(X^3,X^2Y, XY^2,Y^3), \\
\dot{Z} = 0, \\
\end{array} \right.
\label{2-2}
\end{equation}
where $a_1(0) = b_1(0) = 0$ so that the origin is a fixed point
(for the normal form theory, see Chow, Li and Wang~\cite{Cho}).
Then the assumption (A2) means that $a_2(0) \neq 0,\,a_3(0) \neq 0,\, b_2(0) \neq 0$.
In this case, it is well known that the flow of Eq.(\ref{2-2}) has a cusp at the origin
(see also Lemma 3.1).
Since Eq.(\ref{2-2}) has a cusp at $L^+$, there exists exactly one orbit $\alpha ^+$ emerging from $L^+$.
The assumption (A3) means that if the critical manifold is locally convex downward (resp. convex upward ),
then $g(x,y,z,0) <0$ (resp. $g(x,y,z,0) >0$) on $S^+_a \cup \{L^+\}$.
Thus an orbit of the reduced flow on $S^+_a$ reaches $L^+$ in finite time.
As a result, an orbit of (\ref{local1}) may jump in the vicinity of $L^+$.
The next theorem describes an asymptotic behavior of such a jumping orbit.
\\[0.2cm]
\textbf{Theorem 1.} \, Suppose that the system (\ref{local1}) satisfies assumptions (A1) to (A3).
Consider a solution $\mathbf{x}(t)$ whose initial point is in the vicinity of $S^+_a$.
Then, there exist $t_0, t_1 > 0$ such that the distance between $\mathbf{x}(t),\, t_0 < t < t_1$ and 
the orbit $\alpha ^+$ of the unperturbed system emerging from $L^+$ is of $O(\varepsilon ^{4/5})$ as $\varepsilon \to 0$.
\\[-0.2cm]

\begin{figure}[h]
\begin{center}
\includegraphics[scale=1.0]{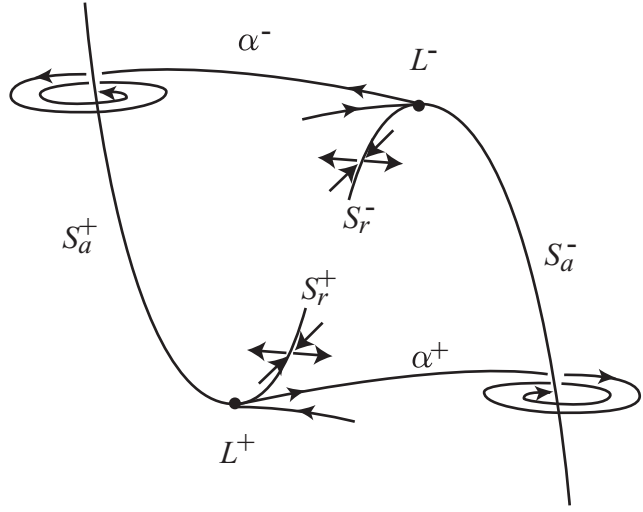}
\end{center}
\caption{Critical manifold and the flow of Eq.(\ref{local1}) with the assumptions (B1) to (B4).}
\label{fig1}
\end{figure}

Note that for a saddle-node type fold points, the distance between $\mathbf{x}(t)$ and an orbit emerging from
a fold point is of $O(\varepsilon ^{2/3})$.
To prove the existence of relaxation oscillations, we need global assumptions for the system (\ref{local1}).
\\[0.2cm]
\textbf{(B1)}\, The critical manifold $\mathcal{M}$ has two smooth components $S^+ = S_a^+ \cup \{L^+\} \cup S_r^+$ 
and $S^- = S_a^- \cup \{L^-\} \cup S_r^-$, where $S_a^{\pm}$ consist of stable focus fixed points, $S_r^{\pm}$ consist
of saddle fixed points, and where $L^{\pm}$ are fold points (see Fig.\ref{fig1}).
\\[0.2cm]
\textbf{(B2)} \, The $L^{\pm}$ are Bogdanov-Takens type fold points; that is,
$L^{\pm}$ are Bogdanov-Takens bifurcation points of the vector field $(f_1(x,y,z,0), f_2(x,y,z,0))$.
In particular, Eq.(\ref{local2}) has cusps at $L^{\pm}$.
\\[0.2cm]
\textbf{(B3)}\, Eq.(\ref{local2}) has two heteroclinic orbits $\alpha ^+$ and $\alpha ^-$ 
which connect $L^+, L^-$ with points on $S_a^-, S_a^+$, respectively.
\\[0.2cm]
\textbf{(B4)} \, The reduced flow (\ref{local4}) on $S_a^{\pm}$ is directed toward the fold points 
$L^{\pm}$ and $g(L^{\pm},0) \neq 0$, respectively.
\\[-0.2cm]

Assumptions (B1) and (B2) assure that $S^{\pm}$ are locally expressed as parabolas, and thus they are of ``J-shaped".
Components $S^+$ and $S^-$ are allowed to be connected.
In this case, $S^+\cup S^-$ is of ``S-shaped".
As was mentioned above, since (\ref{local2}) has cusps at $L^{\pm}$, there exist two orbits $\alpha ^+$ and $\alpha ^-$ of Eq.(\ref{local2})
emerging from $L^+$ and $L^-$.
The assumption (B3) means that these orbits are connected to $S^-_a$ and $S^+_a$, respectively.
If $S^+\cup S^-$ is of ``S-shaped", the assumption (B3) is  typically satisfied
because at least near the fold points, the unperturbed system (\ref{local2}) has heteroclinic orbits
connecting each point on $S^{\pm}_r$ to $S^{\pm}_a$, respectively, due to the basic bifurcation theory.
Note that the assumption (B3) also determines a positional relationship between $S^+$ and $S^-$.
For example, if $S^+$ is convex downward, $S^-$ should be convex upward.
By applying Thm.1 combined with the geometric singular perturbation (boundary layer technique), we can obtain the following result.
\\[0.2cm]
\textbf{Theorem 2.} \, Suppose that the system (\ref{local1}) satisfies assumptions (B1) to (B4).
Then there exists a positive number $\varepsilon _0 $ such that Eq.(\ref{local1}) has a hyperbolically stable 
periodic orbit near $S_a^+ \cup \alpha ^+ \cup S^-_a \cup \alpha ^-$ if $0 < \varepsilon < \varepsilon _0$.
\\[-0.2cm]

To prove the existence of a periodic orbit, the local assumptions are not so important, though a positional relationship between
components of the critical manifold and the existence of heteroclinic orbits 
$\alpha ^{\pm}$ are essential. Indeed, similar results for fast-slow systems having saddle-node type
fold points are obtained by many authors.

To prove the existence of chaos, we have to control the strength of the stability of $S^{\pm}_a$.
Let us consider the system (\ref{1-5}) with $C^\infty$ functions $f_1, f_2, g : U \times I \times I' \to \mathbf{R}$,
where $U$ and $I$ as above and $I' \subset \mathbf{R}$ is a small interval containing zero.
The unperturbed system of Eq.(\ref{1-5}) is given by
\begin{equation}
\left\{ \begin{array}{l}
\dot{x} = f_1(x,y,z,0 ,\delta ), \\ 
\dot{y} = f_2(x,y,z,0 ,\delta ), \\
\dot{z} = 0.
\end{array} \right.
\label{2-1}
\end{equation}
The critical manifold $\mathcal{M}(\delta )$ defined by (\ref{1-6}) is parameterized by $\delta $.
At first, we suppose that the assumptions (B1) to (B4) are satisfied uniformly in $\delta $.
More exactly, we assume following.
\\[0.2cm]
\textbf{(C1)}\, There exists $\delta _0$ such that for every $\delta \in [0, \delta _0)$,
the critical manifold $\mathcal{M}(\delta )$ has two smooth components 
$S^+(\delta ) = S_a^+(\delta ) \cup \{L^+(\delta )\} \cup S_r^+(\delta )$ 
and $S^-(\delta ) = S_a^-(\delta ) \cup \{L^-(\delta )\} \cup S_r^-(\delta )$.
When $\delta >0$, $S_a^{\pm}(\delta )$ consist of stable focus fixed points, $S_r^{\pm}(\delta )$ consist
of saddle fixed points, and $L^{\pm}(\delta )$ are fold points (see Fig.\ref{fig1}).
Further, the $\delta $ family $\mathcal{M}(\delta )$ is smooth with respect to $\delta \in [0, \delta _0)$.
\\[0.2cm]
\textbf{(C2)} \, For every $\delta \in [0, \delta _0)$, $L^{\pm}(\delta )$ are Bogdanov-Takens type fold points; that is,
$L^{\pm}(\delta )$ are Bogdanov-Takens bifurcation points of the vector field $(f_1(x,y,z,0,\delta ), f_2(x,y,z,0,\delta ))$.
In particular, Eq.(\ref{2-1}) has cusps at $L^{\pm}(\delta )$.
\\[0.2cm]
\textbf{(C3)}\, For every $\delta \in (0, \delta _0)$, Eq.(\ref{2-1}) has two heteroclinic orbits 
$\alpha ^+(\delta )$ and $\alpha ^-(\delta )$ 
which connect $L^+(\delta ), L^-(\delta )$ with points on $S_a^-(\delta ), S_a^+(\delta )$, respectively.
\\[0.2cm]
\textbf{(C4)} \, For every $\delta \in [0, \delta _0)$, 
the reduced flow on $S_a^{\pm} (\delta )$ is directed toward the fold points 
$L^{\pm} (\delta )$ and $g(L^{\pm},0,\delta ) \neq 0$, respectively.
\\[-0.2cm]

In addition to the assumptions above, we make the assumptions for the strength of the stability of $S^{\pm}_a$ as follows:
\\[0.2cm]
\textbf{(C5)} \, For every $\delta \in [0, \delta _0)$, 
eigenvalues of the Jacobian matrix $\partial (f_1, f_2)/\partial (x,y)$ of Eq.(\ref{2-1}) at $(x,y,z) \in S^{+}_a(\delta )$
and at $(x,y,z) \in S^{-}_a(\delta )$ are expressed by $-\delta \cdot \mu^{+}(z, \delta ) \pm \sqrt{-1} \omega^+ (z, \delta ) $ and 
$-\delta \cdot \mu^{-}(z, \delta ) \pm \sqrt{-1} \omega^- (z, \delta )$, respectively, where $\mu^{\pm}$ and $\omega ^{\pm}$ are 
real-valued functions satisfying
\begin{equation}
\mu^{\pm} (z, 0) > 0,\,\, \omega ^{\pm} (z, 0) \neq 0.
\label{2-3}
\end{equation}

The assumption (C5) means that the parameter $\delta $ controls the strength of the stability 
of stable focus fixed points on $S^{\pm}_a(\delta )$.

Finally, we suppose that the basin of $S^{\pm}_a (\delta )$ of the unperturbed system can be taken uniformly in $\delta \in (0, \delta_0 )$:
By the assumption (C5), there exist open sets $V^{\pm} \supset S^{\pm}_a (\delta )$ such that 
real parts of eigenvalues of the Jacobian matrix $\partial (f_1, f_2)/\partial (x,y)$ on $V^\pm$ is of order $O(\delta )$.
In general, the ``size" of $V^{\pm}$ depend on $\delta $ and they may tend to zero as $\delta \to 0$.
To prove Theorem 3 below, we assume following.
\\[0.2cm]
\textbf{(C6)}\, There exist open sets $V^{\pm} \supset S^{\pm}_a (\delta )$, which is independent of $\delta $, such that 
real parts of eigenvalues of the Jacobian matrix $\partial (f_1, f_2)/\partial (x,y)$ on $V^\pm$ are negative and of order $O(\delta )$
as $\delta \to 0$.
\\[-0.2cm]

This assumption also assures that the attraction basin of $S^{\pm}_a (\delta )$ of the unperturbed system
can be taken uniformly in $\delta \in (0, \delta_0 )$, see an example below.
It seems that (C6) is satisfied in many situations because of the existence of heteroclinic orbits $\alpha ^\pm (\delta )$.
\\[0.2cm]
\textbf{Theorem 3.} \, Suppose that the system (\ref{1-5}) satisfies assumptions (C1) to (C6).
Then, there exist a positive number $\varepsilon _0$ and positive valued functions 
$\delta _1(\varepsilon ),\, \delta _2(\varepsilon )$ such that 
if $0<\varepsilon <\varepsilon _0$ and $\delta _1(\varepsilon ) < \delta < \delta _2(\varepsilon )$,
then Eq.(\ref{1-5}) has a chaotic invariant set near $S_a^+(\delta ) \cup \alpha ^+(\delta ) \cup S^-_a(\delta ) \cup \alpha ^-(\delta )$,
where $\delta _{1,2} (\varepsilon ) \to 0$ as $\varepsilon \to 0$.
More exactly, the Poincar\'{e} return map $\Pi$ along the flow of (\ref{1-5}) near 
$S_a^+(\delta ) \cup \alpha ^+(\delta ) \cup S^-_a(\delta ) \cup \alpha ^-(\delta )$ is well-defined,
and $\Pi$ has a hyperbolic horseshoe (an invariant Cantor set, on which $\Pi$ is topologically
conjugate to the full shift on two symbols).
\\[-0.2cm]

\begin{figure}[h]
\begin{center}
\includegraphics[scale=1.0]{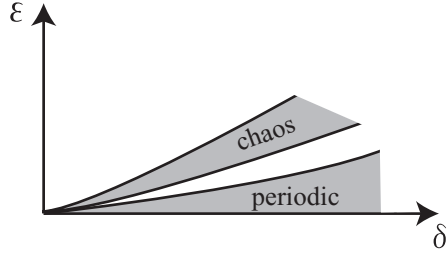}
\end{center}
\caption{Typical bifurcation diagram of (\ref{1-5}) with assumptions (C1) to (C6).}
\label{fig1b}
\end{figure}

Theorems 2 and 3 mean that if $\varepsilon >0$ is sufficiently small for a fixed $\delta $,
then there exists a stable periodic orbit. However, as $\delta $ decreases,
the periodic orbit undergoes a succession of bifurcations
and if $\delta $ gets sufficiently small in comparison with $\varepsilon $, then a chaotic invariant set appears.
In our proof in Sec.4, $\delta $ will be assumed to be of $O(\varepsilon (-\log \varepsilon )^{1/2})$.
We conjecture that this chaotic invariant set is attracting, although the proof is not given in this paper.
In general, given fast-slow systems do not have the parameter $\delta $ explicitly.
However, Theorem 3 suggests that as $\varepsilon $ increases for fixed $\delta $, a periodic orbit undergoes bifurcations
and a chaotic invariant set may appears, see Fig.\ref{fig1b}. 
Obviously the assumptions (C1) to (C4) include assumptions (A1) to (A3) and (B1) to (B4).
In what follows, we consider the system (\ref{1-5}) with the parameter $\delta $.
When proving Theorems 1 and 2, $\delta $ is assumed to be constant, and when proving Theorem 3, 
$\delta $ is assumed to be of $\delta \sim O(\varepsilon (-\log \varepsilon )^{1/2})$ as $\varepsilon \to 0$.
Note that $\varepsilon <\!< \varepsilon (-\log \varepsilon )^{1/2} <\!<1$ as $\varepsilon \to 0$.
Although $\delta >0$ is also small, uniformity assumptions on $\delta $ and the fact $\varepsilon <\!< \delta $
allow us to use the perturbation techniques with respect to only on $\varepsilon $.

\begin{figure}[h]
\begin{center}
\includegraphics[scale=1.0]{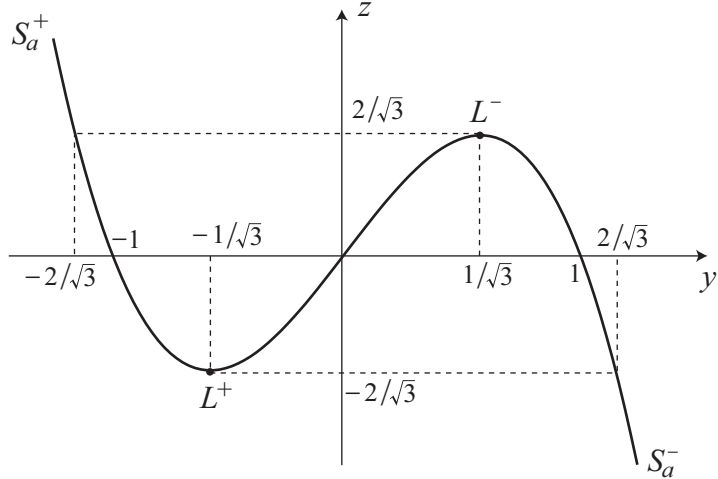}
\end{center}
\caption{Critical manifold of the system (\ref{2-5}).}
\label{fig2}
\end{figure}

In the rest of this section, we give an intuitive explanation of the theorems with an example.
Consider the system
\begin{equation}
\left\{ \begin{array}{l}
\displaystyle \dot{x} = z + 3(y^3 -y) + \delta x(\frac{1}{3} - y^2), \\[0.1cm]
\dot{y} = -x,  \\
\displaystyle \dot{z} = \varepsilon \sin \left( \frac{5}{2} y\right).
\end{array} \right.
\label{2-5}
\end{equation}
The critical manifold $\mathcal{M} = \mathcal{M}(\delta )$ is given by the curve $z = 3(y - y^3), x=0$, and the fold points are given by
$\displaystyle L^{\pm} = (0, \mp \frac{1}{\sqrt{3}}, \mp \frac{2}{\sqrt{3}})$, see Fig.\ref{fig2}.

It is easy to verify that the assumptions (C1), (C2), (C4) and (C5) are satisfied for (\ref{2-5}).
The assumption (C3) of existence of heteroclinic orbits are verified numerically (we do not give a proof here).

The assumption (C6) is also verified by a straightforward calculation.
Now we show that (C6) implies that the attraction basin of $S^{\pm}_a (\delta )$ of the unperturbed system
can be taken uniformly in $\delta \in (0, \delta_0 )$.
We change the coordinates by an affine transformation $(x,y,z) \mapsto (X,Y,Z)$
so that the point $(0, -2/\sqrt{3}, 2/\sqrt{3})$ is placed at the origin
and the linear part of Eq.(\ref{2-5}) is diagonalized.
Then the unperturbed system of Eq.(\ref{2-5}) is rewritten as
\begin{equation}
\frac{d}{dt} \left(
\begin{array}{@{\,}c@{\,}}
X \\
Y
\end{array}
\right) = \frac{\sqrt{-1}}{2}\sqrt{36 - \delta ^2} \left(
\begin{array}{@{\,}cc@{\,}}
-1 & 0 \\
0 & 1
\end{array}
\right) \left(
\begin{array}{@{\,}c@{\,}}
X \\
Y
\end{array}
\right) - \frac{\delta }{2} \left(
\begin{array}{@{\,}cc@{\,}}
1 & 0 \\
0 & 1
\end{array}
\right) \left(
\begin{array}{@{\,}c@{\,}}
X \\
Y
\end{array}
\right) + h(X,Y, \delta ),
\label{2-6}
\end{equation}
where the explicit form of the polynomial $h$, whose degree is greater than one,
is too complicated to be written here.
However, one can verify that $h$ is of the form
\begin{equation}
h(X,Y, \delta ) = \sqrt{-1}h_1(X,Y,\delta ) + \delta h_2(X,Y,\delta ),
\label{2-7}
\end{equation}
where $h_1$ and $h_2$ are polynomials with respect to $X$ and $Y$ such that all coefficients of $h_1$ are real.
Note that $\sqrt{36 - \delta ^2}/2$ and $\delta /2$
correspond to $\omega ^+(z, \delta )$ and $\delta  \mu^+ (z, \delta )$, respectively, in the assumption (C5).
 
Now we bring Eq.(\ref{2-6}) into the normal form with respect to the first term of the right hand side.
There exist a neighborhood $W$ of the origin, which is independent of $\delta $, and a coordinate transformation
$(X,Y) \mapsto (r, \theta )$ defined on $W$ such that Eq.(\ref{2-6}) is put in the form
\begin{equation}
\left\{ \begin{array}{l}
\displaystyle \dot{r} = -\frac{\delta }{2}r + a_3 r^3 + a_5 r^5 + \cdots , \\[0.1cm]
\displaystyle \dot{\theta } = \sqrt{36- \delta ^2}/2 + O(r^2).  \\
\end{array} \right.
\label{2-8}
\end{equation}
Note that the equation of the radius $r$ is independent of $\theta $ (see Chow, Li and Wang~\cite{Cho}).
In our case, $a_3$ is given by
\begin{equation}
a_3 = \delta \frac{-180 + 29 \delta ^2}{6(36 - \delta ^2)^2}.
\label{2-9}
\end{equation}
Further, we can prove that $a_i \sim O(\delta ), \, i=3,5,\cdots $ as $\delta \to 0$
by using the induction together with the property that $h(X,Y,0)$
takes purely imaginary values if $(X,Y) \in \mathbf{R}^2$ (see Eq.(\ref{2-7})).
See Chiba~\cite{Chi1} for explicit formulas of normal forms which are convenient for induction.
Thus the derivative of the right hand side of Eq.(\ref{2-8}) is calculated as
\begin{equation}
\frac{d}{dr}\left(-\frac{\delta }{2}r + a_3r^3 + \cdots  \right)
 = -\frac{\delta }{2}(1+b_3 r^2 + b_5r^4 + \cdots ) + O(\delta ^2),
\label{2-10}
\end{equation}
where $b_3, b_5,\cdots $ are $\delta $-independent constants.
It proves that there exists a $\delta $-independent positive constant $C$ such that if $|r(0)| < C$,
then $r(t)$ decays as $|r(t)| \sim O(e^{-\delta t/2})$ for small $\delta >0$.
The same property can be verified for any system with the assumption (C6) by means of the normal form.

\begin{figure}[h]
\begin{center}
\includegraphics[scale=1.0]{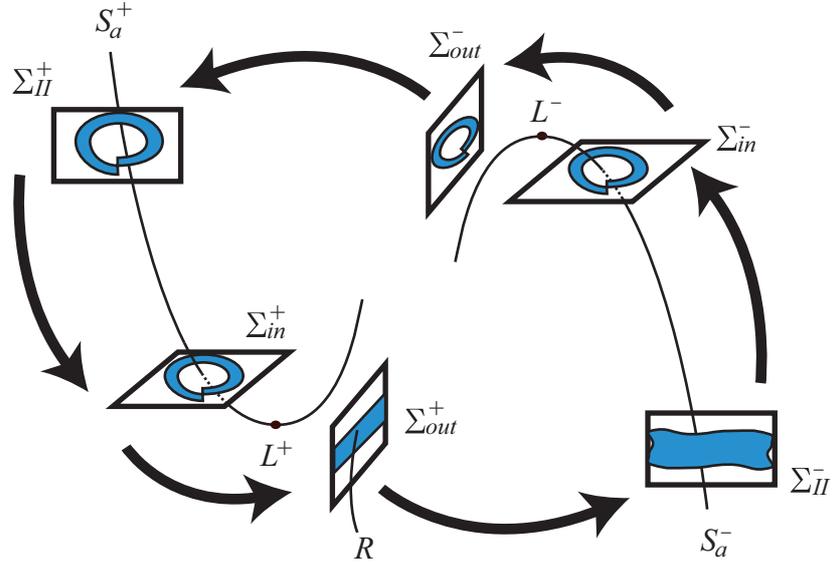}
\end{center}
\caption{Poincar\'{e} sections and a schematic view of the images of the rectangle $R$ under a succession of the 
transition maps. }
\label{fig4}
\end{figure}

To ascertain the reason why the periodic orbit or the chaotic attractor occur, we take Poincar\'{e} sections
$\Sigma ^{+}_{out},\, \Sigma ^{-}_{II},\, \Sigma ^{-}_{in},\, \Sigma ^{-}_{out},\, \Sigma ^{+}_{II}$
and $\Sigma ^{+}_{in}$ as in Fig. \ref{fig4}.

The section $\Sigma ^{+}_{out}$ is parallel to the $xz$ plane and located at the right of $L^+$.
Take a rectangle $R$ on $\Sigma ^{+}_{out}$ and consider how it behaves when it runs along solutions of Eq.(\ref{2-5}).
Since the unperturbed system of Eq.(\ref{2-5}) has the heteroclinic orbit $\alpha ^+$ connecting $L^+$ and $S^-_a$,
the rectangle $R$ also approaches to $S^-_a$ along $\alpha ^+$ and intersects the section $\Sigma ^{-}_{II}$,
as is shown in Fig. \ref{fig4}.
Since the velocity $\varepsilon \sin (5y/2)$ in the direction $z$ is positive in the vicinity of $S^-_a$ and since
$S^-_a$ consists of stable focus fixed points,
the intersection area on $\Sigma ^{-}_{II}$ moves upward, rotating around $S^-_a$.
As a result, the flow of $R$ intersects the section $\Sigma ^{-}_{in}$, which is parallel to the $xy$ plane,
to form a ring-shaped area as is shown in Figure \ref{fig4}.
Further, we can show that the ring-shaped area on $\Sigma ^{-}_{in}$ moves to $\Sigma ^{-}_{out}$ along solutions 
of Eq.(\ref{2-5}) due to Theorem 1.
The area on $\Sigma ^{-}_{out}$ goes back to the section $\Sigma ^{+}_{out}$ in a similar manner because Eq.(\ref{2-5})
has the symmetry $(x,y,z) \mapsto (-x, -y, -z)$.
Thus the Poincar\'{e} return map $\Pi$ from $\Sigma ^{+}_{out}$ into itself is well-defined and it turns out that
$\Pi (R)$ is ring-shaped.

\begin{figure}[h]
\begin{center}
\includegraphics[scale=1.0]{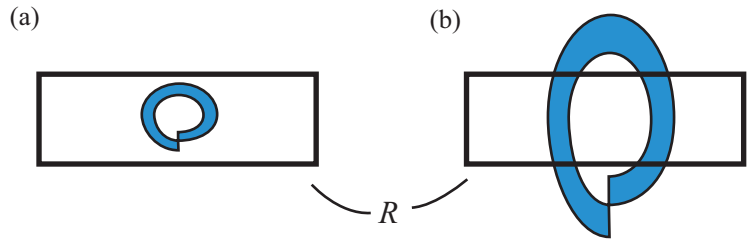}
\end{center}
\caption{Positional relationship of the rectangle $R$ with the returned ring-shaped area. }
\label{fig5}
\end{figure}

There are two possibilities of locations of the returned ring-shaped area.
If the strength of the stability of stable fixed points on $S^{\pm}_a$, say $\delta $ as in the assumption (C5),
is sufficiently large, then the radius of the ring-shaped area gets sufficiently small when passing around $S^{\pm}_a$.
As a result, the returned ring-shaped area is included in the rectangle $R$ as in Fig.\ref{fig5} (a).
It means that the Poincar\'{e} map $\Pi$ is contractive and it has a stable fixed point, which corresponds to a 
stable periodic orbit of Eq.(\ref{2-5}).
On the other hand, if the strength $\delta $ is not so large, the radius of the ring-shaped area is not so small
and it intersects with the rectangle as in Fig.\ref{fig5} (b).
In this case, the Poincar\'{e} map $\Pi$ has a horseshoe.

\section{Local analysis around the fold points}

In this section, we give a local analysis around the fold points $L^{\pm}$ 
by using the blow-up method, and calculate a transition 
map to observe how orbits of Eq.(\ref{1-5}) behave near the fold points.
To prove the existence of chaos, we will give a detailed analysis of the transition map,
which does not need for the standard proof of the existence of a periodic orbit.
The main theorem in this section (Thm.3.2) will be made in the end of Sec.3.1.
We will calculate only for $L^+$ because discussion for $L^-$ is done in the same way.


\subsection{Normal form coordinates}

At first, we transform Eq.(\ref{1-5}) into the normal form in the vicinity of $L^+ (\delta )$.
In what follows, if a (formal) power series $h$ centered at the origin begins with $n$-th degree terms (\textit{i.e.}
$\partial^i h(0)/\partial \bm{x}^i = 0 \, (i= 0,1, \cdots ,n-1)$ and $\partial^n h(0)/\partial \bm{x}^n \neq 0$),
we denote the fact as $h \sim O_p(n)$.
The notation $O(\cdot)$ is also used to the usual Landau notation. 
For example if $h(x,y,z) \sim O(x^2, y^2, z^2, xy, yz, zx)$ as $x,y,z \to 0$,
we simply denote it as $h \sim O_p(2)$.
\\[0.2cm]
\textbf{Lemma 3.1.}\, Suppose (C1), (C2) and (C4). For every $\delta \in [0, \delta _0)$,
There exists a $C^\infty$ local coordinate transformation $(x,y,z) \mapsto (X,Y,Z)$ defined near $L^+(\delta )$
such that Eq.(\ref{1-5}) is brought into the form
\begin{equation}
\left\{ \begin{array}{l}
\dot{X} = Z - Y^2 + c_1(\delta ) XY + Zh_1(X,Y,Z,\delta ) + Y^2h_2(X,Y,Z,\delta ) + \varepsilon h_3(X,Y,Z,\varepsilon,\delta  ),  \\
\dot{Y} = -X + Zh_4(X,Y,Z,\delta ) + \varepsilon h_5(X,Y,Z,\varepsilon ,\delta ),  \\
\dot{Z} = -\varepsilon + \varepsilon h_6(X,Y,Z,\varepsilon ,\delta ),
\end{array} \right.
\label{3-1}
\end{equation}
where $c_1(\delta )$ and $h_i \,\, (i = 1, \cdots ,6)$ are $C^\infty$ functions
such that $c_1 (\delta ) > 0$ for $\delta >0$ and 
\begin{eqnarray}
h_1, h_2, h_4 \sim O(X,Y,Z), \quad h_6\sim O(X,Y,Z,\varepsilon ).
\label{3-1b}
\end{eqnarray}
If we assume (C5), then $c_1 (\delta ) \sim O(\delta )$ as $\delta \to 0$.
\\[0.2cm]
In these coordinates, $L^+(\delta )$ is placed at the origin
and the branch $S^+(\delta )$ of the critical manifold is of the form $Z = Y^2 + O_p(3),\, X = O_p(2)$.
\\[0.2cm]
\textbf{Proof of the Lemma.}\, We start by calculating the normal form of the unperturbed system (\ref{2-1}).
We will use the same notation $(x,y,z)$ as the original coordinates after a succession of coordinate transformations
for simplicity. Since the Jacobian matrix of $(f_1, f_2)$ at  $L^+(\delta )$ has two zero eigenvalues due to the assumption (C1),
the normal form for the equations of $(x,y)$ is of the form (see Chow, Li and Wang~\cite{Cho})
\begin{equation}
\left\{ \begin{array}{l}
\dot{x} = a_1(\delta )z + a_2(\delta )y^2 + a_3(\delta ) xy + zh_1(x,y,z,\delta ) + y^2h_2(x,y,z,\delta ),  \\
\dot{y} = b_1(\delta )x + b_2(\delta )z + zh_4(x,y,z,\delta ),  \\
\end{array} \right.
\label{3-2}
\end{equation}
where $a_1(\delta ), a_2(\delta ), a_3(\delta ), b_1(\delta ), b_2(\delta )$ and $h_1, h_2, h_4 \sim O(x,y,z)$ are $C^\infty$ functions.
Note that $a_2(\delta ) \neq 0, b_1(\delta ) \neq 0$ for $\delta \in [0, \delta _0)$ because of the assumption (C2).
Since we can assume that $S^+(\delta )$ is locally expressed as $z\sim  y^2, x\sim 0$ without loss of generality,
by a suitable coordinate transformation,
we obtain $a_2(\delta ) = -a_1(\delta )$ and $b_2(\delta ) =0$.
Since fixed points on $S_a^+(\delta )$ are attracting and since fixed points on $S^+_r(\delta )$ are saddles
for $\delta >0$,
we obtain $a_1(\delta )b_1(\delta ) < 0$ and $a_3(\delta ) >0$ for $\delta >0$.
If we assume (C5), then $a_3 (\delta ) \sim O(\delta )$.
We can assume that $a_1(\delta ) > 0$ because we are allowed to change the coordinates as $x \mapsto -x, y\mapsto -y$ if necessary.
Thus, the normal form of Eq.(\ref{2-1}) is written as
\begin{equation}
\left\{ \begin{array}{l}
\dot{x} = a_1(\delta )(z - y^2) + a_3(\delta ) xy + zh_1(x,y,z,\delta ) + y^2h_2(x,y,z,\delta ),  \\
\dot{y} = b_1(\delta )x + zh_4(x,y,z,\delta ),  \\
\dot{z} = 0,
\end{array} \right.
\label{3-3}
\end{equation}
with $a_1(\delta ) > 0,\, b_1(\delta ) < 0$. The coordinate transformation which brings Eq.(\ref{2-1}) into Eq.(\ref{3-3}) transforms 
Eq.(\ref{1-5}) into the system of the form
\begin{equation}
\left\{ \begin{array}{l}
\dot{x} = a_1(\delta )(z - y^2) + a_3(\delta ) xy 
               + zh_1(x,y,z,\delta ) + y^2h_2(x,y,z,\delta ) + \varepsilon h_3(x,y,z,\varepsilon ,\delta ),  \\
\dot{y} = b_1(\delta )x + zh_4(x,y,z,\delta ) + \varepsilon h_5(x,y,z,\varepsilon ,\delta ),  \\
\dot{z} = \varepsilon (g_1(\delta ) + h_6(x,y,z,\varepsilon ,\delta )),
\end{array} \right.
\label{3-4}
\end{equation}
where $h_3, h_5, h_6$ are $C^\infty$ functions such that $h_6 \sim O(x,y,z,\varepsilon )$, 
and where $g_1(\delta ) := g(L^+, 0, \delta )$ is a negative constant on account of the assumption (C4).
Finally, changing coordinates and time scales as
\begin{eqnarray}
& & x= -X\frac{a_1(\delta )}{g_1(\delta )}\left( -\frac{g_1(\delta )^2}{a_1(\delta )b_1(\delta )}\right)^{4/5},\,\,
y= Y\left( -\frac{g_1(\delta )^2}{a_1(\delta )b_1(\delta )}\right)^{1/5},\,\,
z= Z\left( -\frac{g_1(\delta )^2}{a_1(\delta )b_1(\delta )}\right)^{2/5},\,\, \nonumber \\
& & t\mapsto -\frac{t}{g_1(\delta )}\left( -\frac{g_1(\delta )^2}{a_1(\delta )b_1(\delta )}\right)^{2/5},
\label{3-5}
\end{eqnarray}
and modifying the definitions of $h_i's\, (i = 1, \cdots ,6)$ appropriately, we obtain Eq.(\ref{3-1}).
Note that since $g_1(\delta ), a_1 (\delta ), b_1 (\delta ) \neq 0$ for $\delta \in [0, \delta _0)$,
this transformation is a local diffeomorphism for every $\delta \in [0, \delta _0)$.
\hfill $\blacksquare$
\\[-0.2cm]

\begin{figure}[h]
\begin{center}
\includegraphics[scale=1.0]{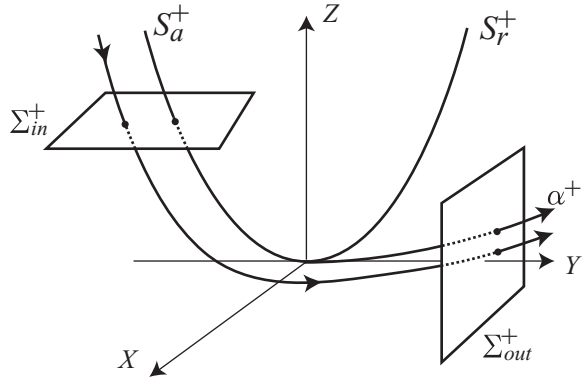}
\end{center}
\caption{Transition map $\Pi^+_{loc}$ and the heteroclinic orbit $\alpha ^+$.}
\label{fig6}
\end{figure}

Let $\rho_1$ be a small positive number and let
\begin{equation}
\Sigma^+_{in} = \{ (X,Y,\rho_1^4) \, | \, (X,Y) \in \mathbf{R}^2 \},\,\,
\Sigma^+_{out} = \{ (X,\rho_1^2, Z) \, | \,  (X,Z) \in \mathbf{R}^2 \}
\label{3-6}
\end{equation}
be Poincar\'{e} sections in the $(X,Y,Z)$ space defined near the origin (see Fig.\ref{fig6}).
The purpose of this section is to construct a transition map from $\Sigma^+_{in}$ to $\Sigma^+_{out}$.
Recall that there exists an orbit $\alpha ^+(\delta )$ emerging from $L^+(\delta )$,
where $L^+(\delta )$ corresponds to the origin in the $(X,Y,Z)$ space.
\\[0.2cm]
\textbf{Theorem 3.2.} \,\, Suppose (C1), (C2) and (C4) to (C6). 
If $\rho_1 > 0$ is sufficiently small, there exists $\varepsilon _0 > 0$ such that
the followings hold for $0 < \varepsilon < \varepsilon _0$ and $0 < \delta < \delta _0$:
\\[0.2cm]
(I) \, There exists an open set $U_\varepsilon  \subset \Sigma^+_{in}$ near the point $\Sigma^+_{in} \cap S^+_a(\delta )$ such that
the transition map $\Pi^+_{loc} : U_\varepsilon  \to \Sigma^+_{out}$ along the flow of Eq.(\ref{3-1}) is well-defined,
$C^\infty$ with respect to $X$ and $Y$, and expressed as
\begin{equation}
\Pi^+_{loc} \left(
\begin{array}{@{\,}c@{\,}}
X \\
Y \\
\rho_1^4
\end{array}
\right) = \left(
\begin{array}{@{\,}c@{\,}}
G_1(\rho_1, \delta ) \\
\rho_1^2 \\
0
\end{array}
\right) + \left(
\begin{array}{@{\,}c@{\,}}
G_2(\mathcal{X}, \mathcal{Y}, \rho_1, \delta )\varepsilon ^{4/5} + O(\varepsilon \log \varepsilon ) \\
0 \\
(\Omega + H(\mathcal{X}, \mathcal{Y}))\varepsilon ^{4/5} + O(\varepsilon \log \varepsilon )
\end{array}
\right),
\label{3-7}
\end{equation}
where $\Omega \sim -3.416$ is a negative constant, and $G_1, G_2, H$ are $C^\infty$ functions with respect to 
$\mathcal{X}, \mathcal{Y}, \delta $.
The arguments $\mathcal{X}, \mathcal{Y}$ are defined by
\begin{equation}
\left\{ \begin{array}{l}
\displaystyle \mathcal{X} = \hat{D}_1(X,Y,\rho_1, \varepsilon , \delta ) 
            \varepsilon ^{-3/5} \exp \Bigl[ -\hat{d}(\rho_1, \delta )\frac{\delta }{\varepsilon }\Bigr],  \\[0.2cm]
\displaystyle \mathcal{Y} = \hat{D}_2(X,Y,\rho_1, \varepsilon , \delta ) 
                     \varepsilon ^{-2/5} \exp \Bigl[ -\hat{d}(\rho_1, \delta )\frac{\delta }{\varepsilon }\Bigr], \\
\end{array} \right.
\label{3-8}
\end{equation}
where $\hat{D}_1, \hat{D}_2$ and $\hat{d}$ are $C^\infty$ functions with respect to $X, Y, \delta $ such that
$\hat{d} > 0$ for $\delta \geq 0$.
Functions $\hat{D}_1$ and $\hat{D}_2$ are not smooth in $\varepsilon $, 
however, they are bounded and nonzero as $\varepsilon \to 0$ and $\delta \to 0$.
\\
(II) \, The point $(G_1(\rho_1, \delta ), \rho_1^2, 0)$ is the intersection of $\alpha ^+(\delta )$ and $\Sigma^+_{out}$.
\\
(III)\, The function $H$ satisfies
\begin{equation}
H(0,0) = 0,\,\, \frac{\partial H}{\partial \mathcal{X}}(\mathcal{X}, \mathcal{Y}) \neq 0.
\label{3-9}
\end{equation}
\\
(IV)\, If $U_\varepsilon $ is sufficiently small, for each $\varepsilon \in (0, \varepsilon _0)$ and $\delta \in (0, \delta _0)$, we can suppose that
\begin{equation}
\frac{\partial \hat{D}_1}{\partial X}(X,Y,\rho_1, \varepsilon , \delta ) \neq 0,
\label{3-9b}
\end{equation}
by changing the value of $\rho_1$ if necessary.

This theorem means that an orbit of Eq.(\ref{1-5}) or Eq.(\ref{3-1}) running around $S^+_a(\delta )$
jumps near $L^+(\delta )$, goes to the right of $L^+(\delta )$ and 
the distance of the orbit and the orbit $\alpha ^+(\delta )$ is of $O(\varepsilon ^{4/5})$ (see Fig.\ref{fig6}).
In particular, it converges to $\alpha ^+(\delta )$ as $\varepsilon \to 0$.
We use the blow-up method to prove this theorem.
In Sec.3.2, we introduce the blow-up coordinates and outline the strategy of the proof of Thm.3.2.
Analysis of our system in the blow-up coordinates is done after Sec.3.3 and the proof is completed in Sec.3.6.
The constant $\Omega $ is a pole of the first Painlev\'{e} equation, as is shown in Sec.3.3.
The function $H$, which is actually an analytic function, also arises from the first Painlev\'{e} equation. 
To prove Theorems 1 and 2, it is sufficient to show that $\mathcal{X}$ and $\mathcal{Y}$ are exponentially small
 as $\varepsilon \to 0$. However, we need more precise decay rate for proving Theorem 3.
For this purpose, the factors $\varepsilon ^{-3/5}$ and $\varepsilon ^{-2/5}$ will be derived by means of the WKB theory. 
Eq.(\ref{3-9}) and (\ref{3-9b}) are also used to prove Theorem 3.
Thus our analysis involves a harder calculation than a usual treatment of fold points in fast-slow systems.
The assumption (C6) is used to assure that the domain $U_\varepsilon $ of the transition map is independent of $\delta \in (0, \delta _0)$.
The assumption (C5) is used to show that the argument of $\exp [\cdots ]$ in Eq.(\ref{3-8}) is of order $O(\delta)$.
For other parts of the theorem, we need only (C1), (C2) and (C4).


\subsection{Blow-up coordinates}

In this subsection, we introduce the blow-up coordinates to ``desingularize" the fixed point $L^+(\delta )$ having a nilpotent
linear part.
Regarding $\varepsilon $ as a dependent variable on $t$, we rewrite Eq.(\ref{3-1}) as
\begin{equation}
\left\{ \begin{array}{l}
\dot{X} = Z - Y^2 + c_1(\delta ) XY + Zh_1(X,Y,Z,\delta ) + Y^2h_2(X,Y,Z,\delta ) + \varepsilon h_3(X,Y,Z,\varepsilon ,\delta ),  \\
\dot{Y} = -X + Zh_4(X,Y,Z,\delta ) + \varepsilon h_5(X,Y,Z,\varepsilon ,\delta ),  \\
\dot{Z} = -\varepsilon + \varepsilon h_6(X,Y,Z,\varepsilon ,\delta ),\\
\dot{\varepsilon } = 0,
\end{array} \right.
\label{3-10}
\end{equation}
with the estimate (\ref{3-1b}).
For this system, we define the blow-up transformations $K_1, K_2$ and $K_3$ to be 
\begin{eqnarray}
(X,Y,Z,\varepsilon ) &=& (r_1^3x_1,\, r_1^2y_1,\, r_1^4,\, r_1^5\varepsilon _1), \label{3-11} \\[0.1cm]
(X,Y,Z,\varepsilon ) &=& (r_2^3x_2,\, r_2^2y_2,\, r_2^4z_2,\, r_2^5), \label{3-12} \\[0.1cm]
(X,Y,Z,\varepsilon ) &=& (r_3^3x_3,\, r_3^2,\, r_3^4z_3,\, r_3^5 \varepsilon _3), \label{3-13}
\end{eqnarray}
respectively, where $K_1, K_2$ and $K_3$ are defined on half spaces 
$\{ Z \geq 0\},\, \{ \varepsilon \geq 0\}$ and $\{ Y \geq 0\}$, respectively.
In what follows, we refer to the coordinates 
$(x_1, y_1, r_1, \varepsilon _1), (x_2, y_2, z_2, r_2), (x_3, r_3, z_3, \varepsilon _3)$
as $K_1, K_2, K_3$ coordinates, respectively.
Transformations $\kappa_{ij}$ from the $K_i$ coordinates to the $K_j$ coordinates are given by
\begin{equation}
\left\{ \begin{array}{l}
\kappa_{12} : (x_2, y_2, z_2, r_2)  
= (x_1 \varepsilon _1^{-3/5},\, y_1 \varepsilon _1^{-2/5},\, \varepsilon _1^{-4/5},\, r_1 \varepsilon _1^{1/5}),\\[0.1cm]
\kappa_{21} : (x_1, y_1, r_1, \varepsilon _1)  = (x_2z_2^{-3/4},\, y_2z_2^{-1/2},\, r_2z_2^{1/4},\, z_2^{-5/4}), \\[0.1cm]
\kappa_{32} : (x_2, y_2, z_2, r_2)  
= (x_3 \varepsilon _3^{-3/5},\, \varepsilon _3^{-2/5},\, z_3 \varepsilon _3^{-4/5},\, r_3 \varepsilon _3^{1/5}), \\[0.1cm]
\kappa_{23} : (x_3, r_3, z_3, \varepsilon _3)  = (x_2y_2^{-3/2},\, r_2y_2^{1/2},\, z_2y_2^{-2},\, y_2^{-5/2}),
\end{array} \right.
\label{3-14}
\end{equation}
respectively. Our next task is to write out Eq.(\ref{3-10}) in the $K_i$ coordinate.
Eqs.(\ref{3-11}) and (\ref{3-10}) are put together to provide
\begin{equation}
\left\{ \begin{array}{l}
\displaystyle \dot{x}_1 = r_1 \bigl(1-y_1^2 + c_1(\delta )r_1x_1y_1 + h_8(x_1, y_1, r_1, \delta )
    + y_1^2 h_9(x_1, y_1, r_1, \delta ) \\[0.1cm]
\displaystyle  \qquad \qquad \qquad  + r_1\varepsilon _1 h_{10}(x_1, y_1, r_1, \varepsilon _1, \delta )
                 + \frac{3}{4}x_1\varepsilon _1 (1 - h_{7}(x_1, y_1, r_1, \varepsilon _1,\delta ))\bigr), \\[0.1cm]
\displaystyle \dot{y}_1 = r_1\bigl( -x_1 + r_1 h_{11}(x_1, y_1, r_1, \delta ) + 
                  r_1^2 \varepsilon _1 h_{12}(x_1, y_1, r_1, \varepsilon _1,\delta )
             +  \frac{1}{2}y_1 \varepsilon _1 (1- h_{7}(x_1, y_1, r_1, \varepsilon _1,\delta )) \bigr), \\[0.1cm]
\displaystyle \dot{r}_1 = -\frac{1}{4}r_1^2 \varepsilon _1(1- h_{7}(x_1, y_1, r_1, \varepsilon _1,\delta )), \\[0.1cm]
\displaystyle \dot{\varepsilon }_1 = \frac{5}{4}r_1 \varepsilon _1^2(1- h_{7}(x_1, y_1, r_1, \varepsilon _1,\delta )),
\end{array} \right.
\label{3-15}
\end{equation}
where $h_i\,\, (i=7, \cdots  ,12)$ are $C^\infty$ functions such that
\begin{eqnarray}
h_7(x_1, y_1, r_1, \varepsilon _1,\delta ) = h_6(r_1^3x_1, r_1^2y_1, r_1^4, r_1^5\varepsilon _1, \delta ),
\label{3-15b}
\end{eqnarray}
and $h_8, \cdots , h_{12}$ are defined in a similar manner through $h_1 , \cdots , h_5$, respectively.
Thus in these functions, $x_1, y_1, \varepsilon _1$ are always with the factors $r_1^3, r_1^2, r_1^5$, respectively.
This fact will be used in later calculations. Note that $h_i \sim O(r_1^2)$ for $i=7,8,9,11$ because of (\ref{3-1b}). 
By changing the time scale appropriately, we can factor out $r_1$ in the right hand side of the above equations:
\begin{equation}
\mathrm{(K_1)} \left\{ \begin{array}{l}
\displaystyle \dot{x}_1 = 1-y_1^2 + c_1(\delta )r_1x_1y_1 + h_8(x_1, y_1, r_1, \delta )
    + y_1^2 h_9(x_1, y_1, r_1, \delta ) \\[0.1cm]
\displaystyle  \qquad \qquad \qquad  + r_1\varepsilon _1 h_{10}(x_1, y_1, r_1, \varepsilon _1, \delta )
                 + \frac{3}{4}x_1\varepsilon _1 (1 - h_{7}(x_1, y_1, r_1, \varepsilon _1,\delta )), \\[0.1cm]
\displaystyle \dot{y}_1 = -x_1 + r_1 h_{11}(x_1, y_1, r_1, \delta ) + 
                  r_1^2 \varepsilon _1 h_{12}(x_1, y_1, r_1, \varepsilon _1,\delta )
             +  \frac{1}{2}y_1 \varepsilon _1 (1- h_{7}(x_1, y_1, r_1, \varepsilon _1,\delta )), \\[0.1cm]
\displaystyle \dot{r}_1 = -\frac{1}{4}r_1 \varepsilon _1(1- h_{7}(x_1, y_1, r_1, \varepsilon _1,\delta )), \\[0.1cm]
\displaystyle \dot{\varepsilon }_1 = \frac{5}{4}\varepsilon _1^2(1- h_{7}(x_1, y_1, r_1, \varepsilon _1,\delta )).
\end{array} \right.
\label{3-16}
\end{equation}
Since the time scale transformation does not change the phase portrait of Eq.(\ref{3-15}),
we can use Eq.(\ref{3-16}) to calculate the transition map.

In a similar manner (\textit{i.e.} changing the coordinates and dividing by the common factors), 
we obtain the systems of equations written in the $K_2, K_3$ coordinates as
\begin{equation}
\mathrm{(K_2)} \left\{ \begin{array}{l}
\dot{x}_2 = z_2 - y_2^2 + r_2 h_{13}(x_2, y_2, z_2, r_2, \delta ), \\[0.1cm]
\dot{y}_2 = -x_2 + r_2^2h_{14}(x_2, y_2, z_2, r_2,\delta ), \\[0.1cm]
\dot{z}_2 = -1 + r_2^2 h_{15}(x_2, y_2, z_2, r_2,\delta ), \\[0.1cm]
\dot{r}_2 = 0,
\end{array} \right.
\label{3-17}
\end{equation}
and
\begin{equation}
\mathrm{(K_3)} \left\{ \begin{array}{l}
\displaystyle \dot{x}_3 = \! -1 \!+ z_3 \!+ c_1(\delta )r_3x_3 
+ \frac{3}{2}x_3h_{16}(x_3, r_3, z_3, \varepsilon _3,\delta ) + r_3^2h_{17}(x_3, r_3, z_3, \varepsilon _3,\delta ), \\
\displaystyle \dot{r}_3 = -\frac{1}{2} r_3h_{16}(x_3, r_3, z_3, \varepsilon _3,\delta ), \\[0.1cm]
\dot{z}_3 = -\varepsilon _3 + 2z_3h_{16}(x_3, r_3, z_3, \varepsilon _3,\delta )
 + r_3^2\varepsilon _3h_{18}(x_3, r_3, z_3, \varepsilon _3,\delta ), \\[0.1cm]
\displaystyle \dot{\varepsilon }_3 = \frac{5}{2} \varepsilon _3h_{16}(x_3, r_3, z_3, \varepsilon _3,\delta ),
\end{array} \right.
\label{3-18}
\end{equation}
respectively, 
where $h_{16}(x_3, r_3, z_3, \varepsilon _3,\delta ) := x_3 +  r_3^2 h_{19}(x_3, r_3, z_3, \varepsilon _3, \delta )$ and 
$h_i \,\, (i = 13, \cdots ,19)$ are $C^\infty$ functions satisfying
\begin{eqnarray*}
h_{17}, h_{18}, h_{19} \sim O(x_3, r_3, z_3, \varepsilon _3).
\end{eqnarray*}

\begin{figure}[h]
\begin{center}
\includegraphics[scale=1.0]{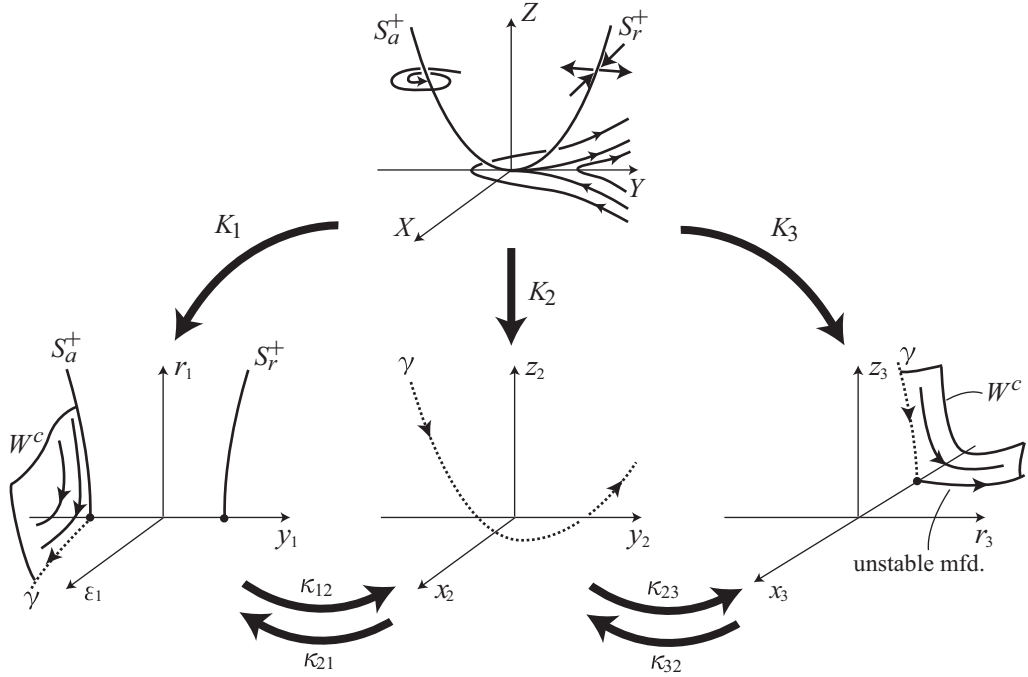}
\end{center}
\caption{The flow in the $(X,Y,Z)$ coordinates and the blow-up coordinates.
The dotted line denotes the orbit $\gamma $ of the first Painlev\'{e} equation.}
\label{fig7}
\end{figure}

Our strategy for understanding the flow of Eq.(\ref{3-1}) near the fold point $L^+(\delta )$ is as follows:
In Sec.3.3, we analyze Eq.(\ref{3-17}) in the $K_2$ coordinates.
We will find it to be a perturbed first Painlev\'{e} equation.
Since asymptotic behavior of the first Painlev\'{e} equation is well studied,
we can construct a transition map along the flow of it approximately.
In Sec.3.4, we analyze Eq.(\ref{3-16}) in the $K_1$ coordinates.
We will see that in the $K_1$ coordinates, 
$S^+_a(\delta )$ has a $2$-dimensional attracting center manifold $W^c(\delta )$ for $\delta >0$ (see Fig.\ref{fig7}).
Since it is attracting, orbits passing nearby $S^+_a(\delta )$ approaches $W^c(\delta )$.
Thus if we construct the invariant manifold $W^c(\delta )$ globally, we can well understand asymptotic behavior
of orbits passing through nearby $S^+_a(\delta )$.
Although usual center manifold theory provides the center manifold $W^c(\delta )$ only locally,
we will show that there exists an orbit $\gamma $, called the Boutroux's tritronqu\'{e}e solution,
of the first Painlev\'{e} equation in the $K_2$ coordinates such that
if it is transformed into the $K_1$ coordinates, it is attached on the edge of $W^c(\delta )$ (see Fig.\ref{fig7}).
This means that the orbit $\gamma $ of the first Painlev\'{e} equation guides the manifold $W^c(\delta )$ and provides a global structure
of it. In Sec.3.5, we analyze Eq.(\ref{3-18}) in the $K_3$ coordinates.
We will see that there exists a fixed point whose unstable manifold is $1$-dimensional.
Since the orbit $\gamma $ of the first Painlev\'{e} equation written in the $K_3$ coordinates approaches the fixed point,
the manifold $W^c(\delta )$ put on the $\gamma $ is also attached on the unstable manifold (see Fig.\ref{fig7}).
The unstable manifold corresponds to the heteroclinic orbit $\alpha ^+(\delta )$ in the $(X,Y,Z)$ coordinates if it is blown down.
This means that orbits of Eq.(\ref{3-1}) coming from a region above $L^+(\delta )$ go to the right of $L^+(\delta )$ (see Fig.\ref{fig6})
and pass near the heteroclinic orbit $\alpha ^+(\delta )$. Thus the transition map $\Pi^+_{loc}$ is well defined.
The fixed point in the $K_3$ coordinates corresponds to a pole of the solution $\gamma $ in the $K_2$ coordinates.
In this way, the value $\Omega $ of the pole appears in the transition map (\ref{3-7}).

Combining transition maps constructed on each $K_i$ coordinates and blowing it down to the $(X,Y,Z)$ coordinates,
we can prove Thm.3.2.

\subsection{Analysis in the $K_2$ coordinates}

We consider Eq.(\ref{3-17}).
Since $r_2 = \varepsilon ^{1/5}$ is a small constant, we are allowed to take the system
\begin{equation}
\left\{ \begin{array}{l}
\dot{x}_2 = z_2 - y_2^2,  \\
\dot{y}_2 = -x_2,  \\
\dot{z}_2 = -1,
\end{array} \right.
\label{3-19}
\end{equation}
as the unperturbed system of Eq.(\ref{3-17}).
This is equivalent to the first Painlev\'{e} equation :
\begin{equation}
\left\{ \begin{array}{l}
\displaystyle \frac{dx_2}{dz_2} = -z_2 + y_2^2,  \\
\displaystyle \frac{dy_2}{dz_2} = x_2, 
\end{array} \right. \quad \mathrm{or} \quad 
\frac{d^2y_2}{dz_2^2} = -z_2 + y_2^2.
\label{3-20}
\end{equation}
It is known that there exists a two parameter family of solutions of the first Painlev\'{e} equation
whose asymptotic expansions are given by
\begin{equation}
\left(
\begin{array}{@{\,}c@{\,}}
x_2(z_2) \\
y_2(z_2)
\end{array}
\right) = \left(
\begin{array}{@{\,}c@{\,}}
\displaystyle -\frac{1}{2}z_2^{-1/2}- \left( \frac{C_1}{8}z_2^{-9/8} - \sqrt{2}C_2z_2^{1/8} \right) \cos \phi
 - \left( \frac{C_2}{8}z_2^{-9/8} + \sqrt{2}C_1z_2^{1/8} \right) \sin \phi + O(z_2^{-3}) \\[0.1cm]
\displaystyle -z_2^{1/2} + C_1 z_2^{-1/8} \cos \phi + C_2 z_2^{-1/8} \sin \phi + O(z_2^{-2})
\end{array}
\right),
\label{3-21}
\end{equation}
as $z_2 \to \infty$ and
\begin{equation}
\left(
\begin{array}{@{\,}c@{\,}}
x_2(z_2) \\
y_2(z_2)
\end{array}
\right) = \left(
\begin{array}{@{\,}c@{\,}}
\displaystyle \frac{-12}{(z_2 - z_0 )^3} + \frac{z_0}{5}(z_2 - z_0) + \frac{1}{2}(z_2 - z_0)^2
 + 4C_3 (z_2 - z_0)^3 + O((z_2 - z_0)^4)\\[0.1cm]
\displaystyle \frac{6}{(z_2 - z_0)^2} + \frac{z_0}{10} (z_2 - z_0)^2  + \frac{1}{6}(z_2 - z_0 )^3 + C_3(z_2 - z_0)^4
 + O((z_2 - z_0)^5) 
\end{array}
\right),
\label{3-22}
\end{equation}
as $z_2 \to z_0 + 0$, where $\phi \sim \displaystyle \frac{4\sqrt{2}}{5}z_2^{5/4}\,\, (z_2 \to \infty) $, 
and where $C_1, C_2, C_3$ and $z_0$ are constants which depend on an initial value.
The value $z_0$ is a movable pole of the first Painlev\'{e} equation (see Ince~\cite{Inc}, Noonburg~\cite{Noo}, Conte~\cite{Con}).
In particular, there exists a unique solution $\gamma $, which corresponds to the case $C_1 = C_2= 0$,
whose asymptotic expansions
as $z_2 \to \infty$ and as $z_2 \to \Omega +0$ are of the form
\begin{equation}
\gamma  : \left(
\begin{array}{@{\,}c@{\,}}
x_2 \\
y_2
\end{array}
\right) = \left(
\begin{array}{@{\,}c@{\,}}
x_2(z_2) \\
y_2(z_2)
\end{array}
\right) = \left(
\begin{array}{@{\,}c@{\,}}
\displaystyle -\frac{1}{2}z_2^{-1/2} + O(z_2^{-3}) \\[0.1cm]
\displaystyle -z_2^{1/2} + O(z_2^{-2})
\end{array}
\right),
\label{3-23}
\end{equation}
and
\begin{equation}
\gamma  : \left(
\begin{array}{@{\,}c@{\,}}
x_2 \\
y_2
\end{array}
\right) = \left(
\begin{array}{@{\,}c@{\,}}
x_2(z_2) \\
y_2(z_2)
\end{array}
\right) = \left(
\begin{array}{@{\,}c@{\,}}
\displaystyle \frac{-12}{(z_2 - \Omega )^3} + \frac{\Omega }{5}(z_2 - \Omega ) + O((z_2 - \Omega )^2)\\[0.1cm]
\displaystyle \frac{6}{(z_2 - \Omega )^2} + \frac{\Omega }{10} (z_2 - \Omega )^2  + O((z_2 - \Omega )^3)
\end{array}
\right),
\label{3-24}
\end{equation}
respectively, where $\Omega \sim -3.416$.
The $\gamma $ is called the Boutroux's tritronqu\'{e}e solution~\cite{Bou, Jos}.

\begin{figure}[h]
\begin{center}
\includegraphics[scale=1.0]{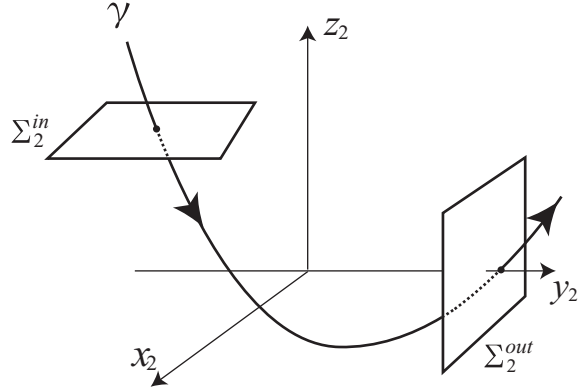}
\end{center}
\caption{The solution $\gamma $ of the first Painlev\'{e} equation and the Poincar\'{e} sections.}
\label{fig8}
\end{figure}

Let $\rho_2$ and $\rho_3$ be small positive numbers and define Poincar\'{e} sections to be
\begin{equation}
\Sigma ^{in}_2  = \{ z_2 = \rho_2^{-4/5}\} ,\,\, \Sigma ^{out}_2  = \{ y_2 = \rho_3^{-2/5}\},
\label{3-25}
\end{equation}
(see Fig.\ref{fig8}). By Eqs.(\ref{3-23}, \ref{3-24}), the intersections $P_2 = \gamma \cap \Sigma ^{out}_2,\, Q_2 = \gamma \cap \Sigma ^{in}_2$
of $\gamma $ and the sections are given by
\begin{eqnarray}
& & P_2 = (p_x, p_y, p_z) = \left( -(2/3)^{1/2}\rho_3^{-3/5} + O(\rho_3^{1/5})
,\, \rho_3^{-2/5},\, \Omega + \sqrt{6}\rho_3^{1/5} + O(\rho_3)  \right), \label{3-26} \label{gamma} \\
& & Q_2 = (q_x, q_y, q_z) = \left( -\rho_2^{2/5}/2 + O(\rho_2^{12/5})
,\, -\rho_2^{-2/5} + O(\rho_2^{8/5}),\, \rho_2^{-4/5}\right) \label{3-27},
\end{eqnarray}
respectively.

\vspace*{0.2cm}
\noindent \textbf{Proposition 3.3.} \, If $\rho_2$ and $\rho_3$ are sufficiently small positive numbers, 
there exists an open set $U_2 \subset \Sigma ^{in}_2$ such that the transition map
$\Pi_2^{loc} : U_2 \to \Sigma ^{out}_2$ along the flow of Eq.(\ref{3-17}) is well-defined and expressed as
\begin{equation}
\Pi_2^{loc} \left(
\begin{array}{@{\,}c@{\,}}
x_2 \\
y_2 \\
\rho_2^{-4/5} \\
r_2
\end{array}
\right) = \left(
\begin{array}{@{\,}c@{\,}}
p_x \\
p_y \\
p_z \\
0
\end{array}
\right) + \left(
\begin{array}{@{\,}c@{\,}}
H_1(x_2 - q_x , y_2 - q_y , \rho_2, r_2, \rho_3, \delta ) \\
0 \\
H_2 (x_2 - q_x , y_2 - q_y , \rho_2, r_2, \rho_3, \delta) \\
r_2
\end{array}
\right),
\label{3-28}
\end{equation}
where $H_1(x,y, \rho_2,r, \rho_3, \delta)$ and $H_2(x,y, \rho_2,r, \rho_3, \delta)$ are $C^\infty$ functions
with respect to $x, y, r$ and $\delta $ satisfying the equalities 
$H_1(0,0, \rho_2,0,\rho_3, \delta) =H_2(0,0, \rho_2,0,\rho_3, \delta) = 0$
for any small $\rho_2, \rho_3>0$ and $\delta \in [0, \delta _0)$.
\\[0.2cm]
\textbf{Proof.} \, This is an immediate consequence of the differentiability of solutions 
with respect to initial values $x_2, y_2$ and parameters $r_2, \delta$.
Note that at this time, we did not prove differentiability at $\rho_3 = 0$,
which will be proved in the next Lemma.
 \hfill $\blacksquare$
\\[-0.2cm]

Since $H_1$ and $H_2$ are $C^\infty$ with respect to $r$ and $\delta $, we put them in the form 
\begin{equation}
H_i(x,y, \rho_2, r, \rho_3, \delta) = \tilde{H}_i(x,y, \rho_2, \rho_3) + O(r),\,\, i=1,2,
\label{3-29}
\end{equation}
where we use the fact that when $r_2=0$, the system (\ref{3-17}) is independent of $\delta $.
Then, the value $\lim_{\rho_3 \to 0} \left( p_z + \tilde{H}_2(x-q_x,y-q_y, \rho_2,\rho_3) \right)$
gives a pole $z_0$ of a solution of Eq.(\ref{3-20}) through an initial point $(x, y, \rho_2^{-4/5})$; that is,
$x_2(z_2), y_2(z_2) \to \infty$ as $z_2\to z_0$.
Prop.3.3 implies that $\tilde{H}_i$ are $C^\infty$ in $x$ and $y$ when $\rho_3 > 0$.
Now we show that $\tilde{H}_i$ can be expanded in $\rho_3^{1/5}$ and they are $C^\infty$ even if $\rho_3 = 0$.
This means that a position of a pole is also smooth with respect to initial values.
In the proof, the Painlev\'{e} property will play a crucial role.
Part (ii) of the next Lemma is used to prove Thm.3.2 (III).
\\[0.2cm]
\textbf{Lemma 3.4.} \, (i) The functions $\tilde{H}_1$ and $\tilde{H}_2$ are analytic with respect to
$(x,y) \in U_2$, $\rho_2^{1/5} > 0$ and $\rho_3^{1/5} \geq 0$, though they are singular at $\rho_2^{1/5} = 0$.
\\
(ii) \, $ \displaystyle \tilde{H}_2(0,0, \rho_2,0)= 0, \quad 
\frac{\partial }{\partial x}\tilde{H}_2(x,y, \rho_2,0) \neq 0$.
\\[0.2cm]
\textbf{Proof.}\, Let $x_2 = x_2(z_2; \rho_2, x_0, y_0)$ and $y_2 = y_2(z_2; \rho_2, x_0, y_0)$
be a solution of the system (\ref{3-20}) with the initial condition
\begin{eqnarray*}
x_2(\rho_2^{-4/5}; \rho_2, x_0, y_0) = x_0, \quad y_2(\rho_2^{-4/5}; \rho_2, x_0, y_0) = y_0.
\end{eqnarray*}
Suppose that $y_2(z) = \rho_3^{-2/5}$ for some $z = z(x_0, y_0, \rho_2, \rho_3)$.
When $\rho_3 > 0$, the statement (i) immediately follows from the fundamental theorem of ODEs:
Since the right hand side of the system (\ref{3-20}) is analytic, any solution is analytic
in time $z_2$, initial time $\rho_2^{-4/5}$ and initial values $(x_0, y_0)$.
Applying the implicit function theorem to the equality
\begin{equation}
y_2(z(x_0, y_0, \rho_2, \rho_3); \rho_2, x_0, y_0) = \rho_3^{-2/5},
\end{equation}
one can verify that 
\begin{equation}
z(x_0, y_0, \rho_2, \rho_3) = p_z + \tilde{H}_2(x_0 - q_x, y_0 - q_y, \rho_2, \rho_3)
\label{3-30}
\end{equation}
is analytic in $x_0, y_0, \rho_2^{1/5}>0$ and $\rho_3^{1/5}>0$.
Thus
\begin{equation}
x_2(z(x_0, y_0, \rho_2, \rho_3); \rho_2, x_0, y_0) = p_x + \tilde{H}_1(x_0 - q_x, y_0 - q_y, \rho_2, \rho_3)
\label{3-31}
\end{equation}
is also analytic in the same region.
Since $z\to \infty$ as $\rho_2 \to 0$, $\tilde{H}_1$ and $\tilde{H}_2$ are singular at $\rho_2^{1/5} = 0$.

When $\rho_3 = 0$, $z(x_0, y_0, \rho_2, 0)$ gives a pole and $x_2 = y_2 = \infty$ at $z_2 = z(x_0, y_0, \rho_2, 0)$.
Thus we should change the coordinates so that a pole becomes a regular point.
For (\ref{3-20}), change the dependent variables $(x_2, y_2)$ and the independent variable $z_2$ to 
$(\xi, \eta)$ and $\tau$ by the relation
\begin{equation}
\left\{ \begin{array}{l}
\displaystyle x_2 = \frac{2\kappa^2}{\eta^3} + \frac{\kappa^2 \tau}{2}\eta 
   + \frac{\kappa^2}{2}\eta^2 - \kappa^2 \eta^3 \xi,  \\[0.2cm]
\displaystyle y_2 = - \frac{\kappa^3}{\eta^2},  \\
\end{array} \right.
\label{3-32}
\end{equation}
and $z_2 = \kappa \tau$, respectively, where $\kappa := (-6)^{1/5}  < 0$.
Then, (\ref{3-20}) is brought into the analytic system
\begin{equation}
\left\{ \begin{array}{l}
\displaystyle \frac{d\eta}{d\tau} = 1 + \frac{\tau}{4} \eta^4 + \frac{1}{4}\eta^5 - \frac{1}{2} \eta^6 \xi,  \\[0.2cm]
\displaystyle \frac{d\xi}{d\tau} = \frac{1}{8}\tau^2\eta + \frac{3}{8}\tau \eta^2
 - \left( \tau \xi - \frac{1}{4}\right) \eta^3 - \frac{5}{4}\eta^4 \xi + \frac{3}{2} \eta^5 \xi^2.  \\
\end{array} \right.
\label{3-33}
\end{equation}
Since any pole of $y_2 (z_2)$ is second order \cite{Inc}, a pole of $y_2$ is transformed into a zero of $\eta (\tau)$
of first order. 
Let $\eta = \eta (\tau; s, \eta_0, \xi_0)$ and $\xi = \xi (\tau; s, \eta_0, \xi_0)$
be a solution of the system satisfying the initial condition
\begin{eqnarray*}
\eta (s; s, \eta_0, \xi_0) = \eta_0, \quad \xi (s; s, \eta_0, \xi_0) = \xi_0,
\end{eqnarray*}
where $(\eta_0, \xi_0)$ and the initial time $s$ correspond to $(x_0, y_0)$ and $\rho_2^{-4/5}$, respectively,
by the transformation (\ref{3-32}).
Suppose that
\begin{eqnarray*}
\eta (\hat{\tau}(s, \eta_0, \xi_0, \rho_3); s, \eta_0, \xi_0) = (-\kappa^3)^{1/2} \rho_3^{1/5}
\end{eqnarray*}
for some $\tau = \hat{\tau}(s, \eta_0, \xi_0, \rho_3)$, which corresponds to a value of
$z(x_0, y_0, \rho_2, \rho_3)$ by the relation $z = \kappa \tau$ so that $y_2(z) = \rho_3^{-2/5}$
(note that when $y_2 = \rho_3^{-2/5}$, then $\eta = (-\kappa^3)^{1/2} \rho_3^{1/5}$).
Since 
\begin{eqnarray*}
\frac{\partial \eta}{\partial \tau}\Bigl|_{\eta = (-\kappa^3)^{1/2} \rho_3^{1/5}} = 1 + O(\rho_3^{4/5}),
\end{eqnarray*}
the implicit function theorem proves that $\hat{\tau}$ is analytic in $s, \eta_0, \xi_0$ and small $\rho_3^{1/5} \geq 0$.
Since the transformation $(\eta_0, \xi_0) \mapsto (x_0, y_0)$ defined through (\ref{3-32}) is analytic when $y_0 \neq 0$,
it turns out that $z(x_0, y_0, \rho_2, \rho_3)$ is analytic in $(x_0, y_0) \in U_2, \rho_2^{1/5}>0$ and $\rho_3^{1/5} \geq 0$.
Now Eqs.(\ref{3-30}, \ref{3-31})  prove the part (i) of Lemma.

To prove (ii), let us calculate the asymptotic expansion of $\hat{\tau}(s, \eta_0, \xi_0,0)$, at which $\eta = 0$.
We rewrite (\ref{3-33}) as
\begin{equation}
\left\{ \begin{array}{l}
\displaystyle \frac{d\tau}{d\eta} = \frac{1}{1 + \frac{\tau}{4} \eta^4 + \frac{1}{4}\eta^5 - \frac{1}{2} \eta^6 \xi} ,  \\[0.4cm]
\displaystyle \frac{d\xi}{d\eta} = \frac{\frac{1}{8}\tau^2\eta + \frac{3}{8}\tau \eta^2
 - \left( \tau \xi - \frac{1}{4}\right) \eta^3 - \frac{5}{4}\eta^4 \xi + \frac{3}{2} \eta^5 \xi^2}
{1 + \frac{\tau}{4} \eta^4 + \frac{1}{4}\eta^5 - \frac{1}{2} \eta^6 \xi}.  \\
\end{array} \right.
\end{equation}
A general solution of this system is obtained in a power series of $\eta$ as
\begin{equation}
\left\{ \begin{array}{l}
\displaystyle \tau = \tau_1 + \eta - \frac{\tau_1}{20} \eta^5 - \frac{1}{12} \eta^6 + \frac{\xi_1}{14}\eta^7 + O(\eta^8),  \\[0.2cm]
\displaystyle \xi = \xi_1 + \frac{\tau_1^2}{16}\eta^2 + \frac{5\tau_1}{24}\eta^3 + O(\eta^4).  \\
\end{array} \right.
\end{equation}
where $\tau_1$ and $\xi_1$ are constants to be determined from an initial condition.
By using the initial condition $(\tau, \eta, \xi) = (s, \eta_0, \xi_0)$, $\tau_1$ is determined as
\begin{equation}
\tau_1 = s - \eta_0 + \frac{s}{20}\eta_0^5 + \frac{1}{30}\eta_0^6 - \frac{\xi_0}{14} \eta_0^7 + O(\eta_0^8).
\label{3-33b}
\end{equation}
When $\eta = 0$, $\tau = \tau_1$. 
This means that the above $\tau_1$ gives the expansion of $\hat{\tau}(s, \eta_0, \xi_0, 0)$.
Then we obtain
\begin{eqnarray*}
\frac{\partial \tilde{H}_2}{\partial x_0} (x_0 - q_x, y_0 - q_y, \rho_2, 0) 
&=& \frac{\partial z}{\partial x_0}(x_0, y_0, \rho_2, 0) \\
&=& \kappa \frac{\partial \hat{\tau}}{\partial x_0}(s, \eta_0, \xi_0, 0) \\
&=& \kappa \frac{\partial \hat{\tau}}{\partial \eta_0} \frac{\partial \eta_0}{\partial x_0}
    + \kappa \frac{\partial \hat{\tau}}{\partial \xi_0} \frac{\partial \xi_0}{\partial x_0} \\
&=& \kappa \left( -\frac{1}{14}\eta_0^7 + O(\eta^8_0) \right) \cdot - \frac{1}{\kappa^2 \eta_0^3},
\end{eqnarray*}
which is not zero for small $\eta_0$ (thus for large $y_0$). 
The equality $\tilde{H}_2(0,0, \rho_2,0)= 0$ is obvious from the definition.
\hfill $\blacksquare$
\\[0.2cm]
\textbf{Remark.} Since $\tilde{H}_i$ is analytic in $\rho_3^{1/5} \geq 0$, 
it is expanded as
\begin{equation}
\tilde{H}_i(x,y,\rho_2, \rho_3) = \hat{H}_i (x,y,\rho_2) + O(\rho_3^{1/5}),
\label{3-34}
\end{equation}
for $i=1,2$.
Indeed, one can verify that
\begin{eqnarray*}
\tilde{H}_i(x,y,\rho_2, \rho_3) = \tilde{H}_i(x,y,\rho_2, 0)
 + \sqrt{6} \rho_3^{1/5} + \frac{3\sqrt{6}}{10}(\tilde{H}_i(x,y,\rho_2, 0) + p_z) \rho_3 + 3\rho_3^{6/5} + O(\rho_3^{7/5})
\end{eqnarray*}
by using the expansion (\ref{3-22}).
Further, $\tilde{H}_i$ are expanded in a Laurent series of $\rho_2^{1/5}$.
In particular, Eq.(\ref{3-33b}) show that the expansions are of the form
\begin{equation}
\hat{H}_i (x,y,\rho_2) = \hat{\hat{H}}_i (x,y) + \rho_2^{-4/5} F_i(x,y,\rho_2^{-4/5}),
\label{3-34b}
\end{equation}
because $s = \rho_2^{-4/5}/\kappa$, where $F_1, F_2$ are analytic functions.
The proof of the above lemma is based on the fact that a pole of (\ref{3-20}) can be transformed into a zero 
of the analytic system by the analytic transformation.
This property is common to Painlev\'{e} equations,
and the transformation (\ref{3-32}) is used to prove that (\ref{3-20}) has the Painlev\'{e} property \cite{Con, Inc}.


\subsection{Analysis in the $K_1$ coordinates}

\begin{figure}[h]
\begin{center}
\includegraphics[scale=1.0]{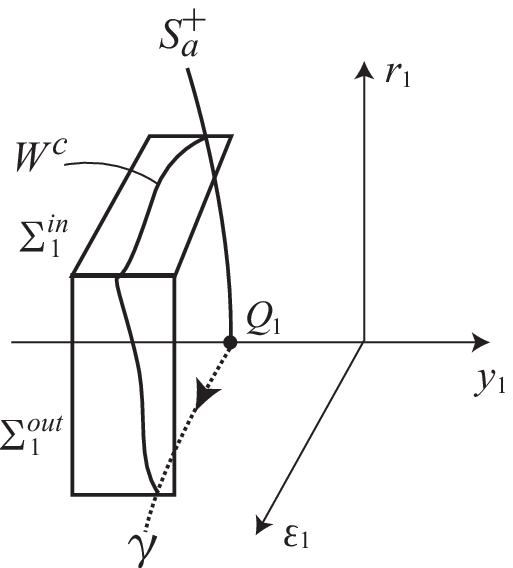}
\end{center}
\caption{Poincar\'{e} sections to define the transition map $\Pi^{loc}_1$.}
\label{fig9}
\end{figure}

We turn to Eq.(\ref{3-16}).
It is easy to verify that Eq.(\ref{3-16}) has fixed points $(x_1, y_1, r_1, \varepsilon _1) = (0, \pm 1, 0,0)$.
By virtue of the implicit function theorem, we can show that there exist two sets of fixed points which form
two curves emerging from $(0, \pm 1, 0,0)$, and
they correspond to $S_a^+(\delta )$ and $S_r^+(\delta )$, respectively (see Fig.\ref{fig7}).
On the fixed points, the Jacobian matrix of the right hand side of Eq.(\ref{3-16}) has eigenvalues given by
\begin{equation}
0,\, 0, \,\, \frac{1}{2}\left( c_1(\delta )r_1y_1 + O(r_1^3) \pm \sqrt{8y_1 - 4c_1(\delta )r_1x_1 + O(r_1^2)} \right).
\label{3-40}
\end{equation}
In particular, the eigenvalues become $0,0,\pm\sqrt{2}i$ at the fixed point $Q_1 = (0, -1,0,0)$, but at fixed points
in $S_a^+(\delta )\backslash Q_1$, they have two eigenvalues whose real parts are negative if $r_1$ is small
and $\delta >0$.
Eigenvectors associated with the two zero eigenvalues at points on $S^+_a(\delta )\backslash Q_1$ converge to those at $Q_1$,
which are given by $(0,0,1,0)$ and $(-1,0,0,2)$, as $r_1 \to 0$.
The vector $(0,0,1,0)$ is tangent to $S^+_a (\delta )$.
Thus $(-1,0,0,2)$ is a nontrivial center direction.
\\[0.2cm]
\textbf{Lemma 3.5.} \, If $\delta >0$, there exists an attracting $2$-dimensional center manifold $W^c(\delta )$ which includes
$S^+_a(\delta )$ and the orbit $\gamma $ of the first Painlev\'{e} equation written in the $K_1$ coordinates (see Fig.\ref{fig9}).
\\[0.2cm]
\textbf{Proof.}\, Let $B(a)$ be the open ball of radius $a$ centered at $Q_1$.
Since at points in $S_a^+(\delta ) \backslash B(a)$ the Jacobian matrix has two zero eigenvalues 
and the other two eigenvalues with negative real parts,
there exists an attracting $2$-dimensional center manifold $W^c(\delta ,a)$ 
emerging from $S_a^+(\delta ) \backslash B(a)$ for any small $a>0$.
Let $\gamma $ be the solution of the first Painlev\'{e} equation described in the previous subsection.
Its asymptotic expansion (\ref{3-23}) is written in the $K_1$ coordinates as
\begin{equation}
\gamma  :  \left(
\begin{array}{@{\,}c@{\,}}
x_1 \\
y_1 \\
r_1 \\
\varepsilon _1
\end{array}
\right) = \left(
\begin{array}{@{\,}c@{\,}}
\displaystyle -\frac{1}{2} z_2^{-5/4} + O(z_2^{-15/4}) \\[0.1cm]
\displaystyle -1 + O(z_2^{-5/2}) \\[0.1cm]
0 \\
z_2^{-5/4}
\end{array}
\right) \quad ( \mathrm{as} \,\, z_2 \to \infty),
\label{3-41}
\end{equation}
by the coordinate change $\kappa_{21}$ (\ref{3-14}).
The curve (\ref{3-41}) approaches the point $Q_1$ as $z_2 \to \infty$ and its tangent vector converges to the eigenvector
$(-1,0,0,2)$ at $Q_1$ as $z_2 \to \infty$.
Thus $\displaystyle W^c(\delta ) := \lim_{a\to 0} W^c(\delta ,a) \cup \gamma $ forms an invariant manifold.\hfill $\blacksquare$
\\[-0.2cm]

Note that $\gamma $ is included in the subspace $\{r_1 =0\}$.
This lemma means that the orbit $\gamma $ guides global behavior of the center manifold $W^c(\delta )$.

Let $\rho_1, \rho_2 > 0$ be the small constants referred to in Thm.3.2 and Prop.3.3, respectively.
Take two Poincar\'{e} sections $\Sigma ^{in}_1$ and $\Sigma ^{out}_1$ defined to be
\begin{eqnarray}
\Sigma ^{in}_1 = \{ (x_1, y_1, r_1, \varepsilon _1) 
\, | \, r_1 = \rho_1, |x_1| \leq \rho_1,\, |y_1 + 1| \leq \rho_1,\, 0 < \varepsilon _1 \leq \rho_2\},  \nonumber \\
\Sigma ^{out}_1 = \{ (x_1, y_1, r_1, \varepsilon _1) 
\, | \, 0\leq r_1 \leq \rho_1, |x_1| \leq \rho_1,\, |y_1 + 1| \leq \rho_1,\, \varepsilon _1 = \rho_2\},
\label{3-42}
\end{eqnarray}
respectively.
Note that $\Sigma ^{in}_1$ is included in the section $\Sigma ^{+}_{in}$ (see Eq.(\ref{3-6})) if written in the $(X,Y,Z)$
coordinates and $\Sigma ^{out}_1$ in the section $\Sigma ^{in}_2$ (see Eq.(\ref{3-25})) if written in the $K_2$ coordinates.
\\

\noindent \textbf{Proposition 3.6.} \, Suppose (C1), (C2) and (C4) to (C6).
\\
(I) \, If $\rho_1$ and $\rho_2$ are sufficiently small, the transition map $\Pi ^{loc}_1 :
\Sigma ^{in}_1 \to \Sigma ^{out}_1$ along the flow of Eq.(\ref{3-16}) is 
well-defined for every $\delta \in (0, \delta _0)$ and expressed as
\begin{equation}
\Pi^{loc}_1 \left(
\begin{array}{@{\,}c@{\,}}
x_1 \\
y_1 \\
\rho_1 \\
\varepsilon _1
\end{array}
\right) = \left(
\begin{array}{@{\,}c@{\,}}
\varphi _1(\rho_1 \varepsilon _1^{1/5} \rho_2^{-1/5}, \rho_2, \delta ) \\
\varphi _2(\rho_1 \varepsilon _1^{1/5} \rho_2^{-1/5}, \rho_2, \delta )  \\
\rho_1 \varepsilon_1 ^{1/5} \rho_{2}^{-1/5} \\
\rho_2
\end{array}
\right) + \left(
\begin{array}{@{\,}c@{\,}}
X_1 \\
Y_1 \\
0 \\
0
\end{array}
\right),
\label{3-43}
\end{equation} 
where $\varphi _1$ and $\varphi _2$ are $C^\infty$ functions such that the graph of $x_1 = \varphi _1(r_1, \varepsilon _1, \delta )$
and $y_1 = \varphi _2 (r_1, \varepsilon _1, \delta )$ gives the center manifold $W^c(\delta )$.
The second term denotes the deviation from $W^c(\delta )$, and $X_1$ and $Y_1$ are defined to be
\begin{equation}
\left\{ \begin{array}{l}
\displaystyle X_1 = D_1(x_1, y_1, \rho_1, \varepsilon _1, \rho_2, \delta ) \left( \frac{\rho_2}{\varepsilon _1} \right)^{3/5}
 \exp \Bigl[ -d(\rho_1, \varepsilon _1, \rho_2, \delta )\frac{\delta }{\varepsilon _1}\Bigr] , \\
\displaystyle Y_1
  = D_2(x_1, y_1, \rho_1, \varepsilon _1, \rho_2, \delta ) \left( \frac{\rho_2}{\varepsilon _1} \right)^{2/5}
 \exp \Bigl[ -d(\rho_1, \varepsilon _1, \rho_2, \delta )\frac{\delta }{\varepsilon _1}\Bigr],
\end{array} \right.
\label{3-44}
\end{equation}
where $D_1, D_2$ and $d$ are $C^\infty$ functions with respect to $x_1, y_1, \rho_1$ and $\delta $.
Although $D_1, D_2$ and $d$ are not $C^\infty$ in $\varepsilon _1$ and $\rho_2$, they are bounded and nonzero as $\varepsilon _1\to 0$
and $\delta \to 0$. Further, they admit the expansions of the form
\begin{eqnarray}
& & D_i(x_1, y_1, \rho_1, \varepsilon _1, \rho_2, \delta ) 
    = \hat{D}_i(x_1, y_1, \rho_1, \varepsilon _1, \delta ) + O((\varepsilon _1/ \rho_2)^{1/5}), \label{3-44b} \\
& & d(\rho_1, \varepsilon _1, \rho_2, \delta ) = \hat{d}(\rho_1, \delta )+ O((\varepsilon _1/ \rho_2)^{1/5}), \label{3-44c}
\end{eqnarray}
for $i = 1,2$.
\\
(II) \, The first term in the right hand side of Eq.(\ref{3-43}) is on the intersection of 
$\Sigma ^{out}_1$ and the center manifold $W^c(\delta )$.
In particular, as $\varepsilon _1 \to 0$, $\Pi^{loc}_1 (x_1, y_1, \rho_1, \varepsilon _1)$ converges to the intersection point
of $\Sigma ^{out}_1$ and $\gamma $.
\\
(III) \, If the initial point $(x_1, y_1, \rho_1, \varepsilon _1)$ is sufficiently close to $W^c(\delta )$,
\begin{equation}
\frac{\partial \hat{D}_1}{\partial x_1} (x_1, y_1, \rho_1, \varepsilon _1, \delta ) \neq 0
\label{3-44d}
\end{equation}
except for a countable set of values of $\varepsilon _1$.
\\[0.2cm]
\textbf{Remark.}\, To prove the existence of a periodic orbit, it is sufficient to show that $X_1$ and $Y_1$
are exponentially small as $\varepsilon _1 \to 0$.
However, to prove the existence of chaos, we need more precise estimate as the factors $(\rho_2/\varepsilon _1)^{3/5}$
and $(\rho_2/\varepsilon _1)^{2/5}$.
Eq.(\ref{3-44d}) is used to prove Eq.(\ref{3-9b}).
\\[0.2cm]
\textbf{Proof.}\, At first, we divide the right hand side of Eq.(\ref{3-16}) by $1-h_{7}$
and change the time scale accordingly.
Note that this does not change the phase portrait. Then we obtain
\begin{equation}
\left\{ \begin{array}{l}
\displaystyle \dot{x}_1 = 1-y_1^2 + c_1(\delta )r_1x_1y_1 + \frac{3}{4}x_1 \varepsilon _1 + h_8 + y_1^2 h_9 + r_1\varepsilon _1 h_{10} \\[0.1cm]
\displaystyle \qquad \qquad \qquad + (1-y_1^2 + c_1(\delta )r_1x_1y_1 + h_8 + y_1^2 h_9 + r_1\varepsilon _1 h_{10}) h_{21}, \\[0.1cm]
\displaystyle \dot{y}_1 = -x_1 + \frac{1}{2}y_1 \varepsilon _1 +  r_1 h_{11} + r_1^2 \varepsilon _1 h_{12} 
         + (-x_1 +  r_1 h_{11} + r_1^2 \varepsilon _1 h_{12})h_{21}, \\[0.1cm]
\displaystyle \dot{r}_1 = -\frac{1}{4}r_1 \varepsilon _1, \\[0.1cm]
\displaystyle \dot{\varepsilon }_1 = \frac{5}{4}\varepsilon _1^2,
\end{array} \right.
\label{3-45}
\end{equation}
where $h_{21} = \sum^\infty_{k=1}h_7^k$, and arguments of functions are omitted.
Equations for $r_1$ and $\varepsilon _1$ are solved as
\begin{equation}
r_1(t) = r_1(0) \left( \frac{4-5\varepsilon _1(0)t}{4}\right)^{1/5} , \quad
\varepsilon _1(t) = \frac{4\varepsilon _1(0)}{4-5\varepsilon _1(0)t},
\label{3-46}
\end{equation}
respectively.
Let $T$ be a transition time from $\Sigma ^{in}_1$ to $\Sigma ^{out}_1$.
Since $\varepsilon _1(T) = \rho_2$, $T$ is given by
\begin{equation}
T = \frac{4}{5\varepsilon _1 (0)}\left( 1 - \frac{\varepsilon _1(0)}{\rho_2}\right).
\label{3-47}
\end{equation}
To estimate $x_1(T)$ and $y_1(T)$, let us introduce the new time variable $\tau$ by
\begin{equation}
\tau = \left( \frac{4-5\varepsilon _1(0)t}{4}\right)^{1/5}.
\label{3-52}
\end{equation}
Then, $r_1(t) = r_1(0)\tau ,\, \varepsilon _1(t) = \varepsilon _1(0) \tau^{-5}$.
Note that when $t = 0$, $\tau = 1$ and when $t= T$, one has $\tau = (\varepsilon _1(0)/ \rho_2)^{1/5}$.
\\[0.2cm]
\textbf{Claim 1.} Any solutions $(x_1, y_1)$ of (\ref{3-45}) are of the form $x_1 = \tau^{-3}u_1(\tau),\, y_1 = \tau^{-2}u_2 (\tau)$,
where $u_1$ and $u_2$ are $C^\infty$ with respect to $\tau$.
\\[0.2cm]
\textbf{Proof.} Changing the time $t$ to $\tau$, the system (\ref{3-45}) is rewritten as
\begin{eqnarray*}
\left\{ \begin{array}{l}
\displaystyle 
-\frac{1}{4}\varepsilon _1(0) \tau^{-4} \frac{dx_1}{d\tau}
  = 1-y_1^2 + c_1(\delta )r_1(0)\tau x_1y_1 + \frac{3}{4}x_1 \varepsilon _1(0) \tau^{-5} + h_8 + y_1^2 h_9 
        + r_1(0)\varepsilon _1(0) \tau^{-4} h_{10} \\[0.1cm]
\displaystyle \qquad \qquad \qquad \qquad + (1-y_1^2 + c_1(\delta )r_1(0)\tau x_1y_1 + h_8 + y_1^2 h_9 
           + r_1(0)\varepsilon _1(0) \tau^{-4} h_{10}) h_{21}, \\[0.1cm]    
\displaystyle 
-\frac{1}{4}\varepsilon _1(0) \tau^{-4} \frac{dy_1}{d\tau}
 = -x_1 + \frac{1}{2} \varepsilon _1(0)\tau^{-5}y_1 +  r_1(0)\tau h_{11} + r_1(0)^2 \varepsilon _1(0)\tau^{-3} h_{12} \\[0.1cm]
\displaystyle \qquad \qquad \qquad \qquad
         + (-x_1 +  r_1(0)\tau h_{11} + r_1(0)^2 \varepsilon _1(0)\tau^{-3} h_{12})h_{21}. \\[0.1cm]
\end{array} \right.
\end{eqnarray*}
Putting $x_1 = \tau ^{-3}u_1, \, y_1 = \tau^{-2} u_2$ yields
\begin{eqnarray}
\left\{ \begin{array}{l}
\displaystyle 
-\frac{1}{4}\varepsilon _1(0) \frac{du_1}{d\tau}
  = \Bigl(\tau ^7 -\tau^3 u_2^2 + c_1(\delta )r_1(0)\tau^3 u_1u_2 + \tau ^7 h_8 + \tau ^3 u_2^2 h_9 
        + r_1(0)\varepsilon _1(0) \tau^{3} h_{10}\Bigr) (1 + h_{21}), \\[0.3cm]
\displaystyle 
-\frac{1}{4}\varepsilon _1(0) \frac{du_2}{d\tau}
 = \Bigl(-\tau^3 u_1 + r_1(0)\tau^7 h_{11} + r_1(0)^2 \varepsilon _1(0)\tau^{3} h_{12}\Bigr)(1 + h_{21}).
\end{array} \right.
\label{3-52b}
\end{eqnarray}
Recall that $h_7$ is defined through (\ref{3-15b}), and thus
\begin{eqnarray}
h_7(x_1, y_1, r_1, \varepsilon _1,\delta ) = h_6(r_1(0)^3u_1, r_1(0)^2u_2, r_1(0)^4\tau^4, r_1(0)^5\varepsilon _1(0), \delta ),
\end{eqnarray}
which implies that $h_7$ is $C^\infty$ with respect to $u_1, u_2, r_1(0), \varepsilon _1(0) ,\delta $ and $\tau$.
Functions $h_8 , \cdots , h_{12}$ and $h_{21}$ have the same property.
Hence the right hand side of Eq.(\ref{3-52b}) is $C^\infty$ with respect to $u_1, u_2, r_1(0), \delta$ and $\tau$,
which proves that solutions $u_1(\tau)$ and $u_2(\tau)$ are $C^\infty$ with respect to $r_1(0), \delta $ and $\tau$.
\hfill $\blacksquare$

Next thing to do is to derive the center manifold and how $x_1(t)$ and $y_1(t)$ approach to it.
The local center manifold $W^c(\delta )$ is given as a graph of $C^\infty$ functions $x_1 = \varphi _1(r_1, \varepsilon _1, \delta ),\, 
y_1 = \varphi _2 (r_1, \varepsilon _1, \delta )$. By using the standard center manifold theory, we can calculate $\varphi _1$
and $\varphi _2$ as
\begin{equation}
\varphi _1 (r_1, \varepsilon _1, \delta ) = -\frac{1}{2}\varepsilon _1 + O(r_1^2, r_1 \varepsilon _1, \varepsilon _1^2),\quad
\varphi _2(r_1, \varepsilon _1, \delta ) = -1 + O(r_1^2, r_1 \varepsilon _1, \varepsilon _1^2).
\label{3-48}
\end{equation}
To see the behavior of solutions $x_1$ and $y_1$ near the center manifold $W^c(\delta )$,
we put $x_1$ and $y_1$ in the form
\begin{equation}
x_1(\tau) =  \varphi _1(r_1(\tau) , \varepsilon _1(\tau), \delta ) + \tau ^{-3}v_1(\tau), \quad
y_1(\tau ) = \varphi _2(r_1(\tau) , \varepsilon _1(\tau), \delta ) + \tau ^{-2}v_2(\tau).
\label{3-49}
\end{equation}
Since $\tau^{3}x_1(\tau)$ and $\tau^{2}y_1(\tau)$ are $C^\infty$ in $\tau$ for every solutions $x_1$ and $y_1$,
so are solutions $\tau^{3}\varphi _1(r_1(\tau) , \varepsilon _1(\tau), \delta )$ and 
$\tau^{2}\varphi _2(r_1(\tau) , \varepsilon _1(\tau), \delta )$
on the center manifold multiplied by $\tau^{3}$ and $\tau^{2}$, respectively.
This implies that $v_1(\tau)$ and $v_2(\tau)$ are also $C^\infty$ in $\tau$.
Substituting Eq.(\ref{3-49}) into (\ref{3-45}) and expanding it in $v_1, v_2$ and $\varepsilon _1(0)$,
we obtain the system of the form
\begin{eqnarray}
\left\{ \begin{array}{l}
\displaystyle \varepsilon _1 \frac{dv_1}{d\tau}
  = -8\tau^5v_2 + 4c_1 r_1 \tau^5 v_1 + r_1^3 \tau^7 h_{22}(r_1, \tau, \delta )v_1 + r_1^2 \tau^7 h_{23}(r_1, \tau, \delta )v_2
 + g_1 (v_1, v_2, r_1, \varepsilon _1, \delta ,\tau),  \\[0.3cm]
\displaystyle \varepsilon _1 \frac{dv_2}{d\tau}
  = 4\tau^3v_1 + r_1^4 \tau^7 h_{24}(r_1, \tau, \delta )v_1 + r_1^3 \tau^7 h_{25}(r_1, \tau, \delta )v_2
 + g_2 (v_1, v_2, r_1, \varepsilon _1, \delta ,\tau),
\end{array} \right.
\label{3-55}
\end{eqnarray}
where $g_1, g_2 \sim O(v_1^2, v_1v_2, v_2^2,\varepsilon _1)$ denote higher order terms, $h_{22}, \cdots  , h_{25}$ are $C^\infty$ functions,
and where $r_1(0), \varepsilon _1(0)$ and $c_1(\delta )$ are denoted by $r_1, \varepsilon _1$ and $c_1$, respectively.
This is a singular perturbed problem with respect to $\varepsilon _1$.
\\[0.2cm]
\textbf{Claim 2.} Any nonzero solutions of this system are expressed as
\begin{equation}
v_1 = D^*_1(\tau, r_1, \varepsilon _1, \delta ; v_{10}, v_{20} ) \exp \Bigl[ -\frac{d^*(\tau, r_1, \delta )}{\varepsilon _1}\Bigr],
\quad v_2 = D^*_2(\tau, r_1, \varepsilon _1, \delta ; v_{10}, v_{20}) \exp \Bigl[ -\frac{d^*(\tau, r_1, \delta )}{\varepsilon _1}\Bigr],
\label{claim2}
\end{equation}
where $v_{10} = v_1(1)$ and $v_{20} = v_2(1)$ are initial values, and
where $D^*_1, D^*_2$ and $d^*$ are $C^\infty$ in $\tau, r_1 , v_{10}, v_{20}$ and $\delta $.
Although $D^*_1$ and $D^*_2$ are not $C^\infty$ in $\varepsilon _1$,
they are bounded and nonzero as $\varepsilon _1 \to 0, \delta \to 0$.
If $v_{10}, v_{20}, r_1$ and $\tau$ are sufficiently small,
\begin{equation}
\frac{\partial D_1^*}{\partial v_{10}} (\tau, r_1, \varepsilon _1, \delta ; v_{10}, v_{20}) \neq 0,
\label{claim2b}
\end{equation}
except for a countable set of values of $\varepsilon _1$.
\\[0.2cm]
\textbf{Proof.} At first, we consider the linearized system of (\ref{3-55}) as
\begin{eqnarray}
\left\{ \begin{array}{l}
\displaystyle \varepsilon _1 \frac{dv_1}{d\tau}
  = -8\tau^5v_2 + 4c_1 r_1 \tau^5 v_1 + r_1^3 \tau^7 h_{22}(r_1, \tau, \delta )v_1 + r_1^2 \tau^7 h_{23}(r_1, \tau, \delta )v_2
 + O(\varepsilon _1),  \\[0.3cm]
\displaystyle \varepsilon _1 \frac{dv_2}{d\tau}
  = 4\tau^3v_1 + r_1^4 \tau^7 h_{24}(r_1, \tau, \delta )v_1 + r_1^3 \tau^7 h_{25}(r_1, \tau, \delta )v_2
 + O(\varepsilon _1),
\end{array} \right.
\label{3-55b}
\end{eqnarray}
which yields the equation of $v_1$ as
\begin{equation}
\varepsilon _1^2 \frac{d^2v_1}{d\tau^2} - 
\varepsilon _1(4c_1 r_1 \tau ^5 + 4r_1^3\tau^7h_{26} + O(\varepsilon _1))\frac{dv_1}{d\tau}
 + (32 \tau^{8} - 4r_1^2\tau^7h_{27} + O(\varepsilon _1)) v_1 = 0,
\label{3-56}
\end{equation}
where $h_{26}(r_1, \tau, \delta )$ and $h_{27}(r_1, \tau, \delta )$ are $C^\infty$ functions.
According to the WKB theory, we construct a solution of this equation in the form
\begin{eqnarray*}
v_1(\tau) = \exp \Bigl[ \frac{1}{\varepsilon _1} \sum^\infty_{n=0} \varepsilon _1^n S_n(\tau) \Bigr].
\end{eqnarray*}
Substituting this into Eq.(\ref{3-56}), we obtain the equation of $S_0(\tau)$
\begin{eqnarray*}
\left( \frac{dS_0}{d\tau}\right)^2 - (4c_1 r_1 \tau^5 + 4r_1^3 \tau ^7 h_{26}) \frac{dS_0}{d\tau} + 
32 \tau^8  - 4 r_1^2\tau^7 h_{27} = 0,
\end{eqnarray*}
which is solved as $S_0 = S_0^{\pm}(\tau) = V(\tau) \pm i W(\tau)$, where
\begin{eqnarray*}
V(\tau) = \int^\tau_{1} \! (2c_1 r_1 s^5 + 2r_1^3 s^7 h_{26}) ds,
\,\, W(\tau) = \int^\tau_{1} \! (2c_1 r_1 s^5 + 2r_1^3 s^7 h_{26}) 
   \sqrt{\frac{8s^8 - r_1^2s^7h_{27} }{(c_1r_1s^5 + r_1^3s^7h_{26})^2}-1}\, ds,
\end{eqnarray*}
are real-valued functions for small $r_1$.
If $r_1>0$ is sufficiently small and if $c_1(\delta) >0,\, 0<\tau < 1$, then $V(\tau)  < 0$. 
For these $S_0^{+}(\tau)$ and $S^-_0(\tau)$, $S_1^{\pm}(\tau), S_2^{\pm}(\tau), \cdots $ are uniquely determined by induction, respectively.
Thus a general solution $v_1(\tau)$ is of the form
\begin{eqnarray*}
v_1(\tau) &=& k^+ \exp [V(\tau)/\varepsilon _1] 
\exp [iW(\tau)/\varepsilon _1] \exp \bigl[ S^+_1 + \varepsilon _1 S^+_2 + \cdots \bigr] \\
& & \quad + \, k^- \exp [V(\tau)/\varepsilon _1] 
\exp [-iW(\tau)/\varepsilon _1] \exp \bigl[ S^-_1 + \varepsilon _1 S^-_2 + \cdots \bigr],
\end{eqnarray*}
where $k^+, k^- \in \mathbf{C}$ are arbitrary constants.
Put
\begin{eqnarray*}
D^*_{1+} = \exp [iW(\tau)/\varepsilon _1] \exp \bigl[ S^+_1 + \varepsilon _1 S^+_2 + \cdots \bigr]
,\quad D^*_{1-} = \exp [-iW(\tau)/\varepsilon _1] \exp \bigl[ S^-_1 + \varepsilon _1 S^-_2 + \cdots \bigr].
\end{eqnarray*}
Then, $v_1$ is rewritten as
\begin{eqnarray*}
v_1 (\tau) = k^+ \exp [V(\tau)/\varepsilon _1]D^*_{1+} + k^- \exp [V(\tau)/\varepsilon _1]D^*_{1-},
\end{eqnarray*}
where $D^*_{1+}$ and $D^*_{1-}$ are $C^\infty$ in $ v_{10}, v_{20}, \tau, r_1$ and $\delta $.
They are not $C^\infty$ in $\varepsilon _1$ because of the factor $1/\varepsilon _1$, however,
they are bounded and nonzero as $\varepsilon _1 \to 0$.
In a similar manner, it turns out that $v_2$ is expressed as
\begin{eqnarray*}
v_2 (\tau) = k^+ \exp [V(\tau)/\varepsilon _1]D^*_{2+} + k^- \exp [V(\tau)/\varepsilon _1]D^*_{2-},
\end{eqnarray*}
where $D^*_{2+}$ and $D^*_{2-}$ are $C^\infty$ in $v_{10}, v_{20}, \tau, r_1, \delta $, and
are bounded and nonzero as $\varepsilon _1 \to 0$.
Therefore, the fundamental matrix of the linear system (\ref{3-55b}) is given as
\begin{equation}
F(\tau) = \left(
\begin{array}{@{\,}cc@{\,}}
D^*_{1+} & D^*_{1-} \\
D^*_{2+} & D^*_{2-} 
\end{array}
\right) \exp [V(\tau)/\varepsilon _1].
\end{equation}
Now we come back to the nonlinear system (\ref{3-55}). We rewrite it in the abstract form as
\begin{eqnarray*}
\varepsilon _1 \frac{d\bm{v}}{d\tau} = A(\tau)\bm{v} + g(\bm{v}, \tau),
\end{eqnarray*}
where $\bm{v} = (v_1, v_2)$, $g = (g_1, g_2)$, and $A(\tau)$ is a matrix defining the linear part of the system.
To estimate the nonlinear terms, the variation-of-constants formula is applied.
Put $\bm{v} = F(\tau)\bm{c}(\tau)$ with $\bm{c}(\tau) = (c_1(\tau ), c_2(\tau)) \in \mathbf{C}^2$.
Then, $\bm{c}(\tau)$ satisfies the equation
\begin{eqnarray}
\frac{d\bm{c}}{d\tau} = \frac{1}{\varepsilon _1} F(\tau)^{-1}g(F(\tau)\bm{c}, \tau ).
\label{3-56b}
\end{eqnarray}
Let $\bm{c} = \bm{c}(\tau, \varepsilon _1)$ be a solution of this equation.
Since $F(\tau) \sim O(e^{V(\tau)/\varepsilon _1})$ tends to zero exponentially as $\varepsilon _1 \to 0$
and since $g$ is nonlinear, the time-dependent vector field defined by the right hand side of (\ref{3-56b})
tends to zero as $\varepsilon _1 \to 0$. 
Since solutions $\bm{c}(\tau, \varepsilon _1)$ are continuous with respect to the parameter $\varepsilon _1$,
it turns out that $\bm{c}(\tau, \varepsilon _1)$ tends to a constant as $\varepsilon _1 \to 0$, which is not zero
except for the trivial solution $\bm{c}(\tau, \varepsilon _1) \equiv 0$.
This proves Eq.(\ref{claim2}) with the desired properties by putting $d^* = -V(\tau)$ and
$D^*_i = D^*_{i+} c_1 + D^*_{i-} c_2\,\, (i = 1,2)$.
Note that since the right hand side of (\ref{3-55}) is not zero at $\delta =0$,
$D^*_1 \nequiv 0, D^*_2 \nequiv 0$ as $\delta \to 0$.

When $r_1 = v_{10} = v_{20} = 0$, the derivatives $\partial v_i/ \partial v_{10},\, (i=1,2)$ with respect to the initial value $v_{10}$
satisfy the initial value problem
\begin{equation}
\left\{ \begin{array}{ll}
\displaystyle \varepsilon _1 \frac{d}{d\tau} \frac{\partial v_1}{\partial v_{10}}(\tau,0, \varepsilon _1, \delta ;0;0)
 = -8 \tau^5 \frac{\partial v_2}{\partial v_{10}}(\tau,0, \varepsilon _1, \delta ;0;0),
  & \displaystyle \frac{\partial v_1}{\partial v_{10}}(1,0, \varepsilon _1, \delta ;0;0) =1, \\
\displaystyle \varepsilon _1 \frac{d}{d\tau} \frac{\partial v_2}{\partial v_{10}}(\tau,0, \varepsilon _1, \delta ;0;0)
 = 4 \tau^3 \frac{\partial v_1}{\partial v_{10}}(\tau,0, \varepsilon _1, \delta ;0;0),
  & \displaystyle \frac{\partial v_2}{\partial v_{10}}(1,0, \varepsilon _1, \delta ;0;0) =0.  \\
\end{array} \right.
\end{equation}
This is exactly solved as
\begin{equation}
\frac{\partial v_1}{\partial v_{10}}(\tau,0, \varepsilon _1, \delta ;0;0)
 = \cos \left( \frac{4 \sqrt{2}}{5\varepsilon _1} (\tau^5 - 1)\right) .
\end{equation}
In particular, 
\begin{equation}
\frac{\partial v_1}{\partial v_{10}}(0,0, \varepsilon _1, \delta ;0;0)
 = \cos \left( \frac{4 \sqrt{2}}{5\varepsilon _1}\right)
\end{equation}
is not zero except for a countable set of values of $\varepsilon _1$.
This and the continuity of solutions of ODE prove Eq.(\ref{claim2b}).
\hfill $\blacksquare$

Let us proceed the proof of Prop.3.6.
For $v_1(\tau)$ and $v_2(\tau)$ in (\ref{claim2}), $x_1(\tau)$ and $y_1(\tau)$ are given as (\ref{3-49}).
Since $\tau = (\varepsilon _1(0)/\rho_2)^{1/5}$ when $t = T$, we obtain
\begin{eqnarray*}
x_1(T) &=& \varphi _1(r_1(0)(\varepsilon _1(0)/\rho_2)^{1/5}, \rho_2, \delta ) \\
& &    + \left( \frac{\rho_2}{\varepsilon _1(0)}\right)^{3/5} 
       D^*_1((\varepsilon _1(0)/\rho_2)^{1/5}, r_1(0), \varepsilon _1(0), \delta ; v_{10}, v_{20}) 
            \exp \Bigl[ -\frac{d^*((\varepsilon _1(0)/\rho_2)^{1/5}, r_1(0), \delta )}{\varepsilon _1(0)}\Bigr], \\
y_1(T) &=& \varphi_2(r_1(0)(\varepsilon _1(0)/\rho_2)^{1/5}, \rho_2, \delta ) \\
& &    + \left( \frac{\rho_2}{\varepsilon _1(0)}\right)^{2/5}
       D^*_2((\varepsilon _1(0)/\rho_2)^{1/5}, r_1(0), \varepsilon _1(0), \delta ; v_{10}, v_{20}) 
            \exp \Bigl[ -\frac{d^*((\varepsilon _1(0)/\rho_2)^{1/5}, r_1(0), \delta )}{\varepsilon _1(0)}\Bigr].
\end{eqnarray*} 
Put
\begin{eqnarray*}
D^*_i((\varepsilon _1(0)/\rho_2)^{1/5}, r_1(0), \varepsilon _1(0), \delta ; v_{10}, v_{20} )
 = D_i(x_1(0), y_1(0), r_1(0), \varepsilon _1(0), \rho_2, \delta )
\end{eqnarray*}
for $i = 1,2$. Since $D^*_i$ is $C^\infty$ in $v_{10}, v_{20}, r_1(0)$ and $\delta $,
$D_i$ is also $C^\infty$ in $x_1(0), y_1(0), r_1(0)$ and $\delta $.
Since $D^*_i$ is $C^\infty$ in $\tau$, $D_i$ is bounded and nonzero as $\varepsilon _1(0) \to 0$.
Finally, let us calculate
\begin{eqnarray*}
d^*((\varepsilon _1(0)/\rho_2)^{1/5}, r_1(0), \delta ) = 
\int^{1}_{(\varepsilon _1(0)/\rho_2)^{1/5}} \! (2c_1(\delta ) r_1(0) \tau^5 + 2r_1(0)^3 \tau^7 h_{26}(r_1(0),\tau, \delta )) d\tau.
\end{eqnarray*}
Due to the mean value theorem, there exists a number $\tau^*>0$ such that
\begin{eqnarray*}
d^*((\varepsilon _1(0)/\rho_2)^{1/5}, r_1(0), \delta ) = 
\frac{1}{3}c_1(\delta ) r_1(0) \left( 1 - \left( \frac{\varepsilon _1(0)}{\rho_2} \right)^{6/5}\right)
 + h_{26} (r_1(0), \tau^*, \delta ) \frac{r_1(0)^3}{4} \left( 1 - \left( \frac{\varepsilon _1(0)}{\rho_2} \right)^{8/5}\right).
\end{eqnarray*}
By the assumption (C5), an orbit of (\ref{3-45}) near the center manifold $W^c (\delta )$ approaches to $W^c (\delta )$
with the rate $O(e^{-\delta \mu^+ t})$.
By the assumption (C6), such an attraction region (basin) of $W^c (\delta )$ exists uniformly in $\delta >0$
at least near the branch $S_a^+ (\delta )$.
Thus $ h_{26} (r_1(0), \tau^*, \delta )$ is of order $O(\delta )$ as well as $c_1(\delta )$
if $\rho_2 > 0$ is sufficiently small.
Therefore, there exists a function $d$, which is $C^\infty$ with respect to $r_1(0)$ and $\delta $,
such that
\begin{eqnarray*}
d^*((\varepsilon _1(0)/\rho_2)^{1/5}, r_1(0), \delta ) = 
d(r_1(0), \varepsilon _1(0), \rho_2, \delta ) \cdot \delta .
\end{eqnarray*}
Since $\mu^+ (z,0) \neq 0$, $d(r_1(0), \varepsilon _1(0), \rho_2, 0 ) \neq 0$. 
Since $D^*_i$ and $d^*$ are $C^\infty$ in $\tau = (\varepsilon _1/\rho_2)^{1/5}$, they admit the expansions (\ref{3-44b}, \ref{3-44c}).
This proves (I) of Prop.3.6.
Proposition 3.6 (II) is clear from the definition of $\varphi _1, \varphi _2$,
and (III) follows from Eq.(\ref{claim2b}).
\hfill $\blacksquare$


\subsection{Analysis in the $K_3$ coordinates}

\begin{figure}[h]
\begin{center}
\includegraphics[scale=1.0]{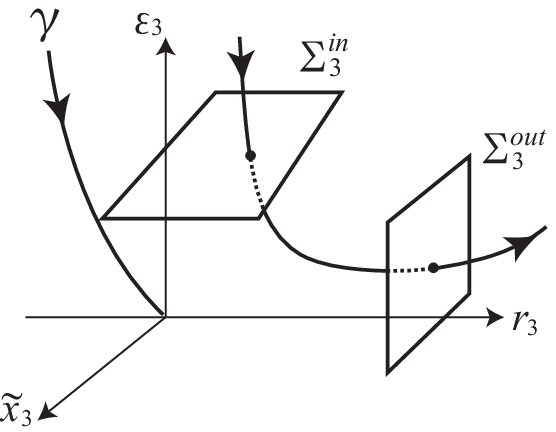}
\end{center}
\caption{Poincar\'{e} sections to define the transition map $\Pi^{loc}_3$.}
\label{fig10}
\end{figure}

We come to the system (\ref{3-18}).
This system has the fixed point $(x_3, r_3, z_3, \varepsilon _3) = (-\sqrt{2/3},0,0,0)$ (see Fig.\ref{fig7}).
To analyze the system, we divide the right hand side of Eq.(\ref{3-18}) by $-h_{16}(x_3, r_3, z_3, \varepsilon _3,\delta )$
and change the time scale accordingly.
Note that this does not change the phase portrait.
At first, note that the equality
\begin{equation}
\frac{1}{h_{16}(x_3, r_3, z_3, \varepsilon _3, \delta )}
 = -\sqrt{\frac{3}{2}}\left( 1+ \sqrt{\frac{3}{2}}\left( x_3 + \sqrt{\frac{2}{3}}\right)
 + h_{31}(x_3 + \sqrt{2/3}, r_3, z_3, \varepsilon _3,\delta ) \right)
\label{3-57}
\end{equation}
holds, where $h_{31} \sim O_p(2)$ is a $C^\infty$ function.
Using Eq.(\ref{3-57}) and introducing the new coordinate by $x_3 + \sqrt{2/3} = \tilde{x}_3$, we eventually obtain
\begin{equation}
\left\{ \begin{array}{l}
\displaystyle \dot{\tilde{x}}_3 
   = -3\tilde{x}_3 + \sqrt{\frac{3}{2}}z_3 - c_1(\delta)r_3 + h_{32}(\tilde{x}_3, r_3, z_3, \varepsilon _3,\delta ), \\[0.1cm]
\displaystyle \dot{r}_3 = \frac{1}{2}r_3, \\[0.1cm]
\displaystyle \dot{z}_3 = -2z_3 - \sqrt{\frac{3}{2}}\varepsilon _3 + \varepsilon _3h_{33}(\tilde{x}_3, z_3, \varepsilon _3,\delta )
 + \varepsilon _3 r_3h_{34}(\tilde{x}_3, r_3, z_3, \varepsilon _3,\delta ), \\[0.1cm]
\displaystyle \dot{\varepsilon }_3 = -\frac{5}{2}\varepsilon _3,
\end{array} \right.
\label{3-58}
\end{equation}
where $h_{32} \sim O_p(2)$ and $h_{33}, h_{34} \sim O_p(1)$ are $C^\infty$ functions.
Note that $h_{33}$ is independent of $r_3$.
This system has a fixed point at the origin,
and eigenvalues of the Jacobian matrix at the origin of the right hand side
of Eq.(\ref{3-58}) are given by $-3, 1/2, -2, -5/2$.
In particular, the eigenvector associated with the positive eigenvalue $1/2$ is given by $(-2c_1(\delta )/7, 1, 0,0)$
and the origin has a $1$-dimensional unstable manifold which is tangent to the eigenvector.
The asymptotic expansion (\ref{3-24}) of the solution $\gamma $ of the first Painlev\'{e} equation is 
rewritten in the present coordinates as
\begin{equation}
(\tilde{x}_3, r_3, z_3, \varepsilon _3) = (O((z_2 - \Omega )^4), \, 0,\, O((z_2 - \Omega )^4),\, O((z_2 - \Omega )^5)),
\label{3-59}
\end{equation}
which converges to the origin as $z_2 \to \Omega $ (see Fig.\ref{fig10}).

Let $\rho_1$ and $\rho_3$ be the small constants introduced in Sec.3.1 and Sec.3.3, respectively.
Define Poincar\'{e} sections $\Sigma ^{in}_3$ and $\Sigma ^{out}_3$ to be
\begin{eqnarray}
& & \Sigma ^{in}_3 = \{(\tilde{x}_3, r_3, z_3, \varepsilon _3)\, | \, |\tilde{x}_3| < \rho_1
,\, 0 < r_3 \leq \rho_1,\, |z_3| \leq \rho_1,\, \varepsilon _3 = \rho_3 \},\label{3-60} \\
& & \Sigma ^{out}_3 = \{(\tilde{x}_3, r_3, z_3, \varepsilon _3)\, | \, |\tilde{x}_3| < \rho_1
,\, r_3 = \rho_1,\, |z_3| \leq \rho_1,\, 0 < \varepsilon _3 \leq \rho_3 \}, \label{3-61}
\end{eqnarray}
respectively (see Fig.\ref{fig10}). Note that $\Sigma ^{in}_3$ is included in the section $\Sigma ^{out}_2$ (see Eq.(\ref{3-25})) 
if written in the $K_2$ coordinates and $\Sigma ^{out}_3$ in the section $\Sigma ^+_{out}$ (see Eq.(\ref{3-6}))
if written in the $(X,Y,Z)$ coordinates.

\noindent \textbf{Proposition 3.7.}
\, (I) \, If $\rho_1$ and $\rho_3$ are sufficiently small, the transition map $\Pi^{loc}_3 : \Sigma ^{in}_3 \to \Sigma ^{out}_3$
along the flow of Eq.(\ref{3-58}) is well-defined and expressed as
\begin{equation}
\Pi^{loc}_3\left(
\begin{array}{@{\,}c@{\,}}
\tilde{x}_3 \\
r_3 \\
z_3 \\
\rho_3
\end{array}
\right) = \left(
\begin{array}{@{\,}c@{\,}}
\displaystyle \beta_1(\rho_1, \delta ) + r_3^4 \beta_2(\tilde{x}_3 ,r_3, z_3, \rho_3, \rho_1, \delta  ) \\
\rho_1 \\
\displaystyle \left( z_3 - \sqrt{6} \rho_3 + \rho_3 \beta_3(\tilde{x}_3, z_3, \rho_3, \delta ) \right) \left( \frac{r_3}{\rho_1} \right) ^4
 + r_3^5 \cdot \log r_3 \cdot \beta_4(\tilde{x}_3, r_3, z_3, \rho_3, \rho_1, \delta ) \\
\displaystyle \rho_3 \left( \frac{r_3}{\rho_1} \right)^5
\end{array}
\right) ,
\label{3-62}
\end{equation}
where $\beta_1$ and $\beta_3$ are $C^\infty$ in their arguments,
$\beta_2$ and $\beta_4$ are $C^\infty$ with respect to $\tilde{x}_3,\, z_3, \rho_3$ and $\delta $
with the property that
$\beta_2$ and $\beta_4$ are bounded as $r_3 \to 0$.
\\
(II) \, As $r_3 \to 0$, $\Pi^{loc}_3(\tilde{x}_3, r_3, z_3, \rho_3)$
converges to the intersection point $(\beta(\rho_1, \delta ), \rho_1, 0 ,0)$
of $\Sigma ^{out}_3$ and the unstable manifold of the origin.
\\[0.2cm]
Before proving Prop.3.7, we need to derive the normal form of Eq.(\ref{3-58}).
\\[0.2cm]
\textbf{Lemma 3.8.}\, In the vicinity of the origin, there exists a $C^\infty$ coordinate transformation
\begin{equation}
\left(
\begin{array}{@{\,}c@{\,}}
\tilde{x}_3 \\
r_3 \\
z_3 \\
\varepsilon _3
\end{array}
\right) = \Phi (X_3, r_3, Z_3, \varepsilon _3,\delta ) := \left(
\begin{array}{@{\,}c@{\,}}
X_3 + \psi_1(X_3, Z_3, \varepsilon _3,\delta ) \\
r_3 \\
Z_3 + \varepsilon _3 \psi_2 (X_3, Z_3, \varepsilon _3,\delta ) \\
\varepsilon _3
\end{array}
\right)
\label{3-63}
\end{equation}
such that Eq.(\ref{3-58}) is transformed into
\begin{equation}
\left\{ \begin{array}{l}
\displaystyle \dot{X}_3 = -3X_3 + \sqrt{\frac{3}{2}}Z_3 - c_1(\delta ) r_3 + r_3h_{35}(X_3, r_3 ,Z_3, \varepsilon _3,\delta ), \\[0.1cm]
\displaystyle \dot{r}_3 = \frac{1}{2}r_3, \\[0.1cm]
\displaystyle \dot{Z}_3 
= -2 Z_3 - \sqrt{\frac{3}{2}}\varepsilon _3 + \varepsilon _3 r_3 h_{36}(X_3, r_3, Z_3, \varepsilon _3,\delta ), \\[0.1cm]
\displaystyle \dot{\varepsilon }_3 = -\frac{5}{2}\varepsilon _3,
\end{array} \right.
\label{3-64}
\end{equation}
where $\psi_2, h_{35}, h_{36} \sim O_p(1)$ and $\psi_{1} \sim O_p(2)$ are $C^\infty$ functions.
\\[0.2cm]
\textbf{Proof of Lemma 3.8.}\, When $r_3 = 0$, Eq.(\ref{3-58}) is written as
\begin{equation}
\left\{ \begin{array}{l}
\displaystyle \dot{\tilde{x}}_3 = -3\tilde{x}_3 + \sqrt{\frac{3}{2}}z_3 
   + h_{32}(\tilde{x}_3, 0, z_3, \varepsilon _3,\delta ),  \\[0.1cm]
\displaystyle \dot{z}_3 
= -2z_3 - \sqrt{\frac{3}{2}}\varepsilon _3 + \varepsilon _3h_{33}(\tilde{x}_3, z_3, \varepsilon _3,\delta ), \\[0.1cm]
\displaystyle \dot{\varepsilon }_3 = -\frac{5}{2}\varepsilon _3.
\end{array} \right.
\label{3-65}
\end{equation}
Since eigenvalues of the Jacobian matrix at the origin of the right hand side of the above are $-3, -2, -5/2$ and 
satisfy the non-resonance condition, there exists a $C^\infty$ transformation of the form $(\tilde{x}_3, z_3, \varepsilon _3)
\mapsto (X_3 + \psi_1 (X_3, Z_3, \varepsilon _3,\delta ), Z_3 +\tilde{\psi}_2 (X_3, Z_3, \varepsilon _3,\delta ), \varepsilon _3)$
such that Eq.(\ref{3-65}) is linearized (see Chow, Li and Wang~\cite{Cho}).
The $\tilde{\psi}_2$ is of the form $\tilde{\psi}_2 = \varepsilon _3 \psi_2$,
where $\psi_2$ is a $C^\infty$ function, because if $\varepsilon _3 = 0$, 
Eq.(\ref{3-65}) gives $\dot{z}_3 = -2z_3$ and it follows that
$Z_3 = z_3$ when $\varepsilon _3 =0$. This transformation brings Eq.(\ref{3-58}) into Eq.(\ref{3-64}).  \hfill $\blacksquare$
\\[0.2cm]
\textbf{Proof of Prop.3.7.}\, 
Note that even in the new coordinates $(X_3, r_3, Z_3, \varepsilon _3)$, 
the sections $\Sigma ^{in}_3 $ and $\Sigma ^{out}_3$ are included in the hyperplanes $\{\varepsilon _3 = \rho_3 \}$ and 
$\{ r_3 = \rho_1\}$, respectively.

Let us calculate the transition time $T$ from $\Sigma ^{in}_3$ to $\Sigma ^{out}_{3}$.
Since $r_3(t) = r_3(0)e^{t/2}$ and $\varepsilon _3(t) = \varepsilon _3(0)e^{-5t/2}$ from Eq.(\ref{3-64}), $T$ is given by
\begin{equation}
T = \log \left( \frac{\rho_1}{r_3(0)}\right)^2.
\label{3-66}
\end{equation}
By integrating the third equation of Eq.(\ref{3-64}), $Z_3(t)$ is calculated as 
\begin{equation}
Z_3(t) = Z_3(0)e^{-2t} +\sqrt{6}\rho_3 (e^{-5t/2} - e^{-2t})
 + e^{-2t}\int^t_{0} \! \rho_3 r_3(0) h_{36} (\bm{X}_3(s),\delta )ds,
\label{3-67}
\end{equation}
where $\bm{X}_3(s) = (X_3(s), r_3(s), Z_3(s), \varepsilon _3(s))$.
Owing to the mean value theorem, there exists $0 \leq \tau = \tau(t) \leq t$ such that
Eq.(\ref{3-67}) is rewritten as
\begin{equation}
Z_3(t) = (Z_3(0) - \sqrt{6}\rho_3) e^{-2t} + \sqrt{6} \rho_3 e^{-5t/2}
 + \rho_3 r_3(0) e^{-2t}h_{36} (\bm{X}_3(\tau ),\delta )t.
\label{3-68}
\end{equation}
This and Eq.(\ref{3-66}) are put together to yield
\begin{equation}
Z_3(T) = (Z_3(0) - \sqrt{6}\rho_3) \left( \frac{r_3(0)}{\rho_1} \right)^4
 + \sqrt{6} \rho_3 \left( \frac{r_3(0)}{\rho_1} \right)^5
 + \rho_3 \frac{r_3(0)^5}{\rho_1^4}h_{36} (\bm{X}_3(\tau (T)),\delta )T.
\label{3-69}
\end{equation}

Next, let us estimate $X_3(T)$. Since $(X_3, r_3)$-plane is invariant,
the unstable manifold of the origin is included in this plane and given as a graph of the $C^\infty$ function
\begin{equation}
X_3 = \phi (r_3, \delta ) = -\frac{2}{7} c_1(\delta ) r_3 + O(r_3^2).
\end{equation}
To measure the distance between $X_3(t)$ and the unstable manifold, put $X_3 = \phi (r_3, \delta ) + u$.
Then, the first equation of (\ref{3-64}) is rewritten as
\begin{eqnarray*}
\dot{u} = (-3 + h_{37}(u, r_3, Z_3, \varepsilon _3, \delta ))u
 + Z_3h_{38}(u, r_3, Z_3, \varepsilon _3, \delta ) + \varepsilon _3 h_{39}(u, r_3, Z_3, \varepsilon _3, \delta ),
\end{eqnarray*}
where $h_{37} \sim O_p(1)$ and $h_{38}, h_{39}$ are $C^\infty$ functions.
This is integrated as
\begin{eqnarray}
\!\! u(t) = e^{-3t}E(t) \left( u(0)\! +\!\! \int^t_{0} \!\! e^{3s}E(s)^{-1}
 (Z_3(s)h_{38}(\bm{u}(s), \delta ) + \varepsilon _3(s) h_{39}(\bm{u}(s), \delta ) ) ds \right),
\label{3-69b}
\end{eqnarray}
where $\bm{u}(s) = (u(s), r_3(s), Z_3(s), \varepsilon _3(s))$ and $E(t) = \exp [\int^t_{0} \! h_{37}(\bm{u}(s), \delta )ds ]$.
Substituting Eq.(\ref{3-68}) and $\varepsilon _3(t) = \rho_3 e^{-5t/2}$ and estimating with the aid of 
the mean value theorem, one can verify that $u(T)$ is of the form
\begin{equation}
u(T) = r_3(0)^4 h_{40} (X_3(0), r_3(0), Z_3(0), \rho_3, \rho_1, \delta ),
\end{equation}
where $h_{40}$ is bounded as $r_3(0) \to 0$
(the factor $Z_3(s)$ in Eq.(\ref{3-69b}) yields the factor $r_3(0)^4$, and other terms are of $O(r_3(0)^5 \log r_3(0)$).
Since the transition time $T$ is not $C^\infty$ in $\rho_1$ and $r_3(0)$,
$h_{40}$ is $C^\infty$ only in $X_3(0), Z_3(0), \rho_3 $ and $\delta$.
Thus the transition map $\tilde{\Pi}_3^{loc}$ from $\Sigma ^{in}_3$ to $\Sigma ^{out}_3$ along the flow of Eq.(\ref{3-64})
is given by
\begin{equation}
\tilde{\Pi}_3^{loc}\left(
\begin{array}{@{\,}c@{\,}}
X_3 \\
r_3 \\
Z_3 \\
\rho_3
\end{array}
\right) = \left(
\begin{array}{@{\,}c@{\,}}
\displaystyle \phi (\rho_1, \delta ) + r_3^4 h_{40} (X_3, r_3, Z_3, \rho_3, \rho_1,\delta ) \\
\displaystyle \rho_1 \\
\displaystyle (Z_3 - \sqrt{6}\rho_3) \left( \frac{r_3}{\rho_1} \right)^4
 + \sqrt{6} \rho_3 \left( \frac{r_3}{\rho_1} \right)^5
 -2 \rho_3 \frac{r_3^5}{\rho_1^4} \log \left( \frac{r_3}{\rho_1} \right)h_{41}(X_3, r_3, Z_3, \rho_3, \rho_1,\delta ) \\
\displaystyle \rho_3 \left( \frac{r_3}{\rho_1} \right)^5
\end{array}
\right) ,
\label{3-71}
\end{equation}
where $h_{41}(X_3, r_3, Z_3, \rho_3, \rho_1,\delta ) = h_{36} (\bm{X}_3(\tau (T)),\delta )$ is  bounded as $r_3 \to 0$
because $\bm{X}_3(\tau (T))$ is bounded.
Since the transition time $T$ is not $C^\infty$ in $\rho_1$ and $r_3(0)$,
$h_{41}$ is $C^\infty$ in $X_3(0), Z_3(0), \rho_3 $ and $\delta$.
Now Eq.(\ref{3-62}) is verified by calculating $\Phi \circ \tilde{\Pi}_3^{loc} \circ \Phi^{-1}$. 
Note that $\beta_3$ in Eq.(\ref{3-62}) is independent of $r_3$ and $\rho_1$ because it comes from the inverse of the 
transformation (\ref{3-63}), which is of the form
\begin{eqnarray*}
\Phi^{-1} (\tilde{x}_3, r_3, z_3, \rho_3)
 = \left(
\begin{array}{@{\,}c@{\,}}
\tilde{x}_3 + \beta_5 (\tilde{x}_3, z_3, \rho_3, \delta ) \\
r_3 \\
z_3 + \rho_3 \beta_3 (\tilde{x}_3, z_3, \rho_3, \delta )\\
\rho_3
\end{array}
\right)
\end{eqnarray*}
with $C^\infty$ functions $\beta_3$ and $\beta_5$.
The unstable manifold $\beta_1(\rho_1, \delta )$ in $(\tilde{x}_3, r_3, z_3, \varepsilon _3)$ coordinate
is obtained from that in $(X_3, r_3, Z_3, \varepsilon _3)$ coordinate as $\beta_1 (\rho_1, \delta )
 = \phi (\rho_1, \delta ) + \psi_1(\phi (\rho_1, \delta ), 0, 0, \delta )$.
 This proves Prop.3.7 (I).
To prove (II) of Prop.3.7, note that the hyperplane $\{ r_3 =0\}$ 
is invariant and included in the stable manifold of the origin.
Since a point $(\tilde{x}_3, r_3, z_3, \rho_3)$ converges to the stable manifold as $r_3 \to 0$,
$\Pi^{loc}_3(\tilde{x}_3, r_3, z_3, \rho_3)$ converges to the unstable manifold as $r_3 \to 0$ 
on account of the $\lambda $-lemma. This proves Prop.3.7 (II).
 \hfill $\blacksquare$


\subsection{Proof of Theorem 3.2}

We are now in a position to prove Theorem 3.2.
Let $\tau_x : (x,r,z, \varepsilon ) \mapsto (x - \sqrt{2/3}, r ,z, \varepsilon )$ be the translation 
in the $x$ direction introduced in Sec.3.5. Eq.(\ref{3-7}) is obtained by writing out the map 
$\tilde{\Pi}^+_{loc} := \tau_x \circ \Pi^{loc}_3 \circ \tau_x^{-1} \circ \kappa_{23} 
\circ \Pi^{loc}_2 \circ \kappa_{12} \circ \Pi^{loc}_1$ and blowing it down to the $(X,Y,Z)$ coordinates.
At first, $\Pi^{loc}_2 \circ \kappa_{12} \circ \Pi^{loc}_1$ is calculated as
\begin{eqnarray}
& & \left(
\begin{array}{@{\,}c@{\,}}
x_1 \\
y_1 \\
\rho_1 \\
\varepsilon _1
\end{array}
\right) \stackrel{\Pi^{loc}_1}{\longmapsto} 
\left(
\begin{array}{@{\,}c@{\,}}
\varphi _1 + X_1 \\
\varphi _2 + Y_1 \\
\rho_1 \varepsilon _1^{1/5}\rho_2^{-1/5} \\
\rho_2
\end{array}
\right) \stackrel{\kappa_{12}}{\longmapsto} 
\left(
\begin{array}{@{\,}c@{\,}}
\rho_2^{-3/5} \varphi _1 + \rho_2^{-3/5} X_1\\
\rho_2^{-2/5} \varphi _2 + \rho_2^{-2/5} Y_1 \\
\rho_2^{-4/5} \\
\rho_1 \varepsilon _1^{1/5}
\end{array}
\right) \nonumber \\
&\stackrel{\Pi^{loc}_2}{\longmapsto} & \left(
\begin{array}{@{\,}c@{\,}}
p_x+ H_1(\rho_2^{-3/5} \varphi _1 + \rho_2^{-3/5} X_1 - q_x,
\rho_2^{-2/5} \varphi _2 + \rho_2^{-2/5} Y_1-q_y , \rho_2,  \rho_1 \varepsilon _1^{1/5}, \rho_3,\delta ) \\
\rho_3^{-2/5} \\
p_z+ H_2(\rho_2^{-3/5} \varphi _1 + \rho_2^{-3/5} X_1 - q_x,
\rho_2^{-2/5} \varphi _2 + \rho_2^{-2/5} Y_1- q_y  , \rho_2, \rho_1 \varepsilon _1^{1/5}, \rho_3, \delta ) \\
\rho_1 \varepsilon _1^{1/5}
\end{array}
\right) , \nonumber \\
\label{3-72}
\end{eqnarray}
where $\varphi _1 = \varphi _1(\rho_1\varepsilon _1^{1/5}\rho_2^{-1/5}, \rho_2, \delta ),\, 
\varphi _2 = \varphi _2(\rho_1\varepsilon _1^{1/5}\rho_2^{-1/5}, \rho_2, \delta ) $, and $X_1, Y_1$
are defined by Eq.(\ref{3-44}).
In what follows, we omit the arguments of $H_1$ and $H_2$.
The last term in the above is further mapped to 
\begin{eqnarray}
&\stackrel{\kappa_{23}}{\longmapsto} & \left(
\begin{array}{@{\,}c@{\,}}
\rho_3^{3/5}p_x+ \rho_3^{3/5}H_1 \\
\rho_1 \rho_3^{-1/5}\varepsilon _1^{1/5} \\
\rho_3^{4/5}p_z+ \rho_3^{4/5}H_2 \\
\rho_3
\end{array}
\right) \stackrel{\tau_x^{-1}}{\longmapsto} 
\left(
\begin{array}{@{\,}c@{\,}}
\sqrt{2/3} + \rho_3^{3/5} p_x + \rho_3^{3/5}H_1 \\
\rho_1 \rho_3^{-1/5}\varepsilon _1^{1/5} \\
\rho_3^{4/5}p_z+ \rho_3^{4/5}H_2 \\
\rho_3
\end{array}
\right) 
:= \left(
\begin{array}{@{\,}c@{\,}}
\tilde{x}_3 \\
r_3 \\
z_3 \\
\rho_3
\end{array}
\right).
\label{3-73}
\end{eqnarray}
Let us denote the resultant as $(\tilde{x}_3, r_3, z_3, \rho_3)$ as above.
Then, $\tilde{\Pi}^+_{loc}$ proves to be given by
\begin{eqnarray}
\tilde{\Pi}^+_{loc} \left(
\begin{array}{@{\,}c@{\,}}
x_1 \\
y_1 \\
\rho_1 \\
\varepsilon _1 
\end{array}
\right) &=& \left(
\begin{array}{@{\,}c@{\,}}
\displaystyle -\sqrt{2/3} + \beta_1(\rho_1, \delta ) + r_3^4 \beta_2 (\tilde{x}_3 , r_3, z_3, \rho_3, \rho_1,\delta ) \\
\rho_1 \\
\displaystyle (z_3 - \sqrt{6} \rho_3 + \rho_3 \beta_3 (\tilde{x}_3, z_3, \rho_3,\delta ))\left( \frac{r_3}{\rho_1}\right)^4
 + r_3^5 \cdot \log r_3 \cdot \beta_4 (\tilde{x}_3,r_3, z_3, \rho_3, \rho_1,\delta ) \\
\displaystyle \rho_3 \left( \frac{r_3}{\rho_1}\right)^5
\end{array}
\right) \nonumber \\
&=& \left(
\begin{array}{@{\,}c@{\,}}
\displaystyle -\sqrt{2/3} + \beta_1(\rho_1, \delta ) + \rho_1^4 \rho_3^{-4/5}\varepsilon _1^{4/5} 
\beta_2 (\tilde{x}_3 ,\rho_1 \rho_3^{-1/5}\varepsilon _1^{1/5}, z_3, \rho_3, \rho_1,\delta ) \\
\rho_1 \\
\displaystyle (z_3 - \sqrt{6} \rho_3 + \rho_3 \beta_3 (\tilde{x}_3, z_3, \rho_3,\delta ))
             \left( \frac{\varepsilon _1}{\rho_3}\right)^{4/5}
 + O(\varepsilon _1 \log \varepsilon _1) \\
\varepsilon _1
\end{array}
\right) .
\label{3-74}
\end{eqnarray}
By using the definition of $p_z$ in (\ref{gamma}), the third component of the above is calculated as
\begin{eqnarray}
& & (z_3 - \sqrt{6} \rho_3 + \rho_3 \beta_3 (\tilde{x}_3, z_3, \rho_3,\delta ))\left( \frac{\varepsilon _1}{\rho_3}\right)^{4/5}
 + O(\varepsilon _1 \log \varepsilon _1) \nonumber \\
&=& \left( \Omega + O(\rho_3) + H_2(\hat{X},\hat{Y} , \rho_2,\rho_1 \varepsilon _1^{1/5}, \rho_3,\delta )
 + \rho_3^{1/5} \beta_3(\tilde{x}_3, z_3, \rho_3,\delta ) \right) \varepsilon _1^{4/5}
 + O(\varepsilon _1 \log \varepsilon _1), \quad
\label{3-75}
\end{eqnarray}
where 
\[ \hat{X} = \rho_2^{-3/5} \varphi _1 + \rho_2^{-3/5} X_1 - q_x, \quad
\hat{Y} = \rho_2^{-2/5} \varphi _2 + \rho_2^{-2/5} Y_1 - q_y. \]
From Eqs.(\ref{3-29}) and (\ref{3-34}), Eq.(\ref{3-75}) is rewritten as
\begin{equation}
\left( \Omega + O(\rho_3) + \hat{H}_2(\hat{X}, \hat{Y} , \rho_2 ) + O(\rho_3^{1/5})
 + \rho_3^{1/5} \beta_3(\tilde{x}_3, z_3, \rho_3, \delta ) \right) \varepsilon _1^{4/5}
 + O(\varepsilon _1 \log \varepsilon _1).
\label{3-76}
\end{equation}
Since $\tilde{\Pi}^{+}_{loc} (x_1, y_1, \rho_1, \varepsilon _1)$ 
is independent of $\rho_3$, which is introduced to define the intermediate sections $\Sigma ^{out}_2$
and $\Sigma ^{in}_3$, all terms including $\rho_3$ are canceled out and Eq.(\ref{3-76}) has to be of the form
\begin{equation}
\left( \Omega + \hat{H}_2(\hat{X}, \hat{Y}, \rho_2) \right) \varepsilon _1^{4/5}
 + O(\varepsilon _1 \log \varepsilon _1).
\label{3-77}
\end{equation}
Now we look into $\hat{X}$ and $\hat{Y}$.
Since $\varphi _1(r_1, \varepsilon _1,\delta )$ and $\varphi _2(r_1, \varepsilon _1,\delta )$ give the graph of 
the center manifold $W^c(\delta )$ and since the orbit $\gamma $ of the first Painlev\'{e} equation is 
attached on the edge of $W^c(\delta )$, $x_1 = \varphi _1(0, \varepsilon _1,\delta )$ and $y_1 = \varphi _2(0, \varepsilon _1, \delta )$
coincide with $\gamma $ written in the $K_1$ coordinates.
Thus we obtain
\begin{eqnarray}
\hat{X} &=& \rho_2^{-3/5} \varphi _1(\rho_1 \varepsilon _1^{1/5} \rho_2^{-1/5}, \rho_2,\delta ) 
            + \rho_2^{-3/5} X_1 - q_x \nonumber \\
&=& \left( \rho_2^{-3/5} \varphi _1(0, \rho_2,\delta ) - q_x \right) 
      + \rho_2^{-3/5}X_1 + \rho_2^{-3/5} O((\varepsilon _1/ \rho_2)^{1/5}) \nonumber \\
&=& \rho_2^{-3/5}X_1 + \rho_2^{-3/5} O((\varepsilon _1/ \rho_2)^{1/5}) \nonumber \\
&=& \rho_2^{-3/5}D_1 \left( \frac{\rho_2}{\varepsilon _1}\right)^{3/5} 
    \exp \Bigl[ -\frac{d\delta }{\varepsilon _1 } \Bigr]
    + \rho_2^{-3/5}O((\varepsilon _1/ \rho_2)^{1/5}) \nonumber \\
&=& D_1 \varepsilon _1^{-3/5} \exp \Bigl[ -\frac{d\delta }{\varepsilon _1 } \Bigr]
     + \rho_2^{-3/5}O((\varepsilon _1/ \rho_2)^{1/5}).
\label{3-78}
\end{eqnarray}
The $\hat{Y}$ is calculated in the same manner.
Functions $D_1, D_2$ and $d$ are expanded as Eqs.(\ref{3-44b}, \ref{3-44c}),
and $\hat{H}_2$ is expanded as (\ref{3-34b}).
Since Eq.(\ref{3-77}) should be independent of $\rho_2$, which is introduced to define the intermediate sections
$\Sigma ^{out}_1$ and $\Sigma ^{in}_2$, Eq.(\ref{3-77}) is rewritten as
\begin{equation}
\left( \Omega + \hat{\hat{H}}_2( \mathcal{X} , \mathcal{Y}) \right) \varepsilon _1^{4/5}
 + O(\varepsilon _1 \log \varepsilon _1).
\label{3-79}
\end{equation}
where
\begin{equation}
\mathcal{X} = \hat{D}_1(x_1, y_1, \rho_1, \varepsilon _1,\delta )\varepsilon _1^{-3/5} 
\exp \Bigl[ -\frac{\hat{d}(\rho_1, \delta )\delta }{\varepsilon _1 }\Bigr],
\,\, \mathcal{Y} = \hat{D}_2(x_1, y_1, \rho_1, \varepsilon _1,\delta ) \varepsilon _1^{-2/5} 
\exp \Bigl[ -\frac{\hat{d}(\rho_1, \delta )\delta }{\varepsilon _1 } \Bigr] .
\end{equation}
Similarly, since the first component of Eq.(\ref{3-74}) is independent of $\rho_2$ and $\rho_3$,
we find that it is expressed as
\begin{eqnarray}
-\sqrt{2/3} + \beta_1 (\rho_1, \delta )
  + \hat{G} (\mathcal{X}, \mathcal{Y}, \rho_1, \delta ) \varepsilon _1^{4/5} +  O(\varepsilon _1 \log \varepsilon _1)
\label{3-80}
\end{eqnarray}
with some $C^\infty$ function $\hat{G}$.

Our final task is to blow down Eq.(\ref{3-74}) with Eqs.(\ref{3-79}, \ref{3-80}) to the $(X,Y,Z)$ coordinates to obtain Eq.(\ref{3-7}).
By the transformation (\ref{3-11}), a point $(X,Y,\rho_1^4, \varepsilon )$ in $(X,Y,Z, \varepsilon )$-space
is mapped to the point $(X\rho_1^{-3}, Y\rho_1^{-2}, \rho_1, \varepsilon \rho_1^{-5})$ in $K_1$-space.
Further, it is mapped by the transition map $\tilde{\Pi}^+_{loc}$ to
\begin{eqnarray*}
\left(
\begin{array}{@{\,}c@{\,}}
-\sqrt{2/3} + \beta_1 (\rho_1, \delta )
  + \hat{G} (\mathcal{X}, \mathcal{Y}, \rho_1, \delta ) \varepsilon ^{4/5}\rho_1^{-4} +  O(\varepsilon  \log \varepsilon ) \\
\rho_1 \\
\left( \Omega + \hat{\hat{H}}_2( \mathcal{X} , \mathcal{Y}) \right) \varepsilon ^{4/5} \rho_1^{-4}
 + O(\varepsilon \log \varepsilon ) \\
\varepsilon \rho_1^{-5}
\end{array}
\right),
\end{eqnarray*}
in $K_3$-space, in which
\begin{eqnarray*}
& & \mathcal{X} = \hat{D}_1(X\rho_1^{-3}, Y\rho_1^{-2}, \rho_1, \varepsilon \rho_1^{-5}, \delta )
     \varepsilon ^{-3/5} \rho_1^3
         \exp \Bigl[ -\frac{\hat{d}(\rho_1, \delta )\delta }{\varepsilon  \rho_1^{-5}}\Bigr], \\
& & \mathcal{Y} = \hat{D}_2(X\rho_1^{-3}, Y\rho_1^{-2}, \rho_1, \varepsilon \rho_1^{-5}, \delta )
     \varepsilon^{-2/5} \rho_1^2
         \exp \Bigl[ -\frac{\hat{d}(\rho_1,  \delta )\delta }{\varepsilon \rho_1^{-5}}\Bigr].
\end{eqnarray*}
Finally, it is blown down by (\ref{3-13}) as
\begin{eqnarray*}
\left(
\begin{array}{@{\,}c@{\,}}
-\sqrt{2/3}\rho_1^3 + \beta_1 (\rho_1, \delta )\rho_1^3
  + \rho_1^{-1} \hat{G} (\mathcal{X}, \mathcal{Y}, \rho_1, \delta ) \varepsilon _1^{4/5} +  O(\varepsilon _1 \log \varepsilon _1) \\
\rho_1^2 \\
\left( \Omega + \hat{\hat{H}}_2( \mathcal{X} , \mathcal{Y}) \right) \varepsilon _1^{4/5}
 + O(\varepsilon _1 \log \varepsilon _1) \\
\varepsilon
\end{array}
\right).
\end{eqnarray*}
By changing the definitions of $\hat{D}_1, \hat{D}_2$ and $\hat{d}$ appropriately, we obtain Theorem 3.2 (I)
with
\begin{eqnarray*}
G_1 = -\sqrt{2/3}\rho_1^3 + \beta_1 (\rho_1, \delta )\rho_1^3,
\, G_2 = \rho_1^{-1} \hat{G},\, H = \hat{\hat{H}}_2.
\end{eqnarray*}

Theorem 3.2 (II) follows from the fact that the unstable manifold
described in Prop.3.7 (II) coincides with the heteroclinic orbit $\alpha ^+(\delta )$ if written in the $(X,Y,Z)$ coordinates.
Theorem 3.2 (III) follows from Lemma 3.4,
and (IV) follows from Eq.(\ref{3-44d}) because $\varepsilon _1$ in Eq.(\ref{3-44d}) is now replaced by $\varepsilon \rho_1^{-5}$. 
This complete the proof of Theorem 3.2 \hfill $\blacksquare$


\section{Global analysis and the proof of main theorems}

In this section, we construct a global Poincar\'{e} map by combining a succession of transition maps (see Fig.\ref{fig4})
and prove Theorems 1,2 and 3.

\subsection{Global coordinate}

\begin{figure}[h]
\begin{center}
\includegraphics[scale=1.0]{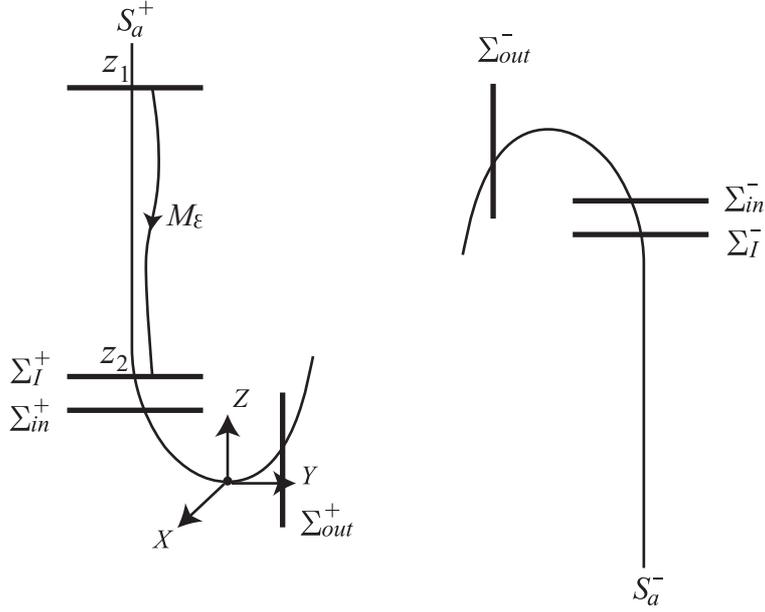}
\end{center}
\caption{Coordinate for calculating the global Poincar\'{e} map, and a slow manifold $M_\varepsilon $
corresponding to the segment $S^+_a(z_2, z_1) = \{ (X,Y,Z) \in S^+_a \, | \, z_2 \leq Z \leq z_1 \}$
of the critical manifold.}
\label{fig11}
\end{figure}

Let us introduce a global coordinate to calculate the global Poincar\'{e} map.
In what follows, we suppose without loss of generality that the branch $S^+$ and $S^-$ of the critical manifold 
are convex downward and upward, respectively, as is shown in Fig.\ref{fig1}.
Recall that $(X,Y,Z)$ coordinate is defined near the fold point $L^+$ and that the sections $\Sigma ^+_{in}$
and $\Sigma ^+_{out}$ are defined in Eq.(\ref{3-6}).
We define a global coordinate transformation $(x,y,z) \mapsto (X,Y,Z)$ satisfying following:
We suppose that in the $(X,Y,Z)$ coordinate, $L^+(\delta ) = (0,0,0),\, L^-(\delta ) = (0, y_0, z_0)$ with $y_0 >0, z_0>0$,
and that $Y$ coordinates of $S^-_a$ are larger than those of $S^+_a$ just as shown in Fig.\ref{fig11}.
Let $z_1 > z_0$ be a number and put $z_2 = \rho_1^{4} + e^{-1/\varepsilon ^2}$.
Define the new section
\begin{equation}
\Sigma ^+_{I} = \{ Z = \rho_1^{4} + e^{-1/\varepsilon ^2}\},
\end{equation}
which lies slightly above $\Sigma ^+_{in}$.
Change the coordinates so that the segment of $S^+_a$ in the region $z_2 \leq Z \leq z_1$ is expressed as
\begin{equation}
\{ X = 0,\, Y = -\eta,\, z_2 \leq Z \leq z_1 \},
\end{equation}
where $\eta$ is a sufficiently small positive constant (if $\rho_1$ is sufficiently small).
We can define such a coordinate without changing the local coordinate near $L^+$ and
the expression of $\Pi^{+}_{loc}$ given in Eq.(\ref{3-7})
by using a partition of unity. We can change the coordinates near $S^-_a \cup \{ L^-\}$ in a similar manner without
changing the coordinate expression near $S^+_a \cup \{ L^+\}$. Let
\begin{equation}
\left\{ \begin{array}{l}
\dot{X} = f_1(X,Y,Z, \varepsilon , \delta ),  \\
\dot{Y} = f_2(X,Y,Z, \varepsilon , \delta ),  \\
\dot{Z} = \varepsilon g(X,Y,Z, \varepsilon , \delta ) ,
\end{array} \right.
\label{4-0}
\end{equation}
be the system (\ref{1-5}) written in the resultant coordinate,
where the definitions of $f_1, f_2$ and $g$ are accordingly changed.

\subsection{Flow near the slow manifold}

Put $S^+_a(z_2, z_1) = \{ (X,Y,Z) \in S^+_a \, | \, z_2 \leq Z \leq z_1 \}$.
Then, $S^+_a(z_2, z_1)$ is a compact attracting normally hyperbolic invariant manifold of the unperturbed system of (\ref{4-0}),
see Fig.\ref{fig11}.
In this subsection, we construct an approximate flow around the slow manifold $M_\varepsilon $
corresponding to $S^+_a(z_2, z_1)$.
If the parameter $\delta $ is a constant, the existence of the slow manifold immediately follows from Fenichel's theorem:
\\[0.2cm]
\textbf{Theorem \, (Fenichel \cite{Fen2}).}

Let $N$ be a $C^r$ manifold ($r \geq 1$), and $\mathcal{X}^r(N)$ the set of $C^r$ vector fields
on $N$ with the $C^1$ topology. 
Let $F$ be a $C^r$ vector field on $N$ and suppose that $M \subset N$ is a compact normally hyperbolic $F$-invariant manifold.
Then, there exists a neighborhood $\mathcal{U} \subset \mathcal{X}^r(N)$ of the origin
such that if $\varepsilon $ is a small positive number so that $\varepsilon G \in \mathcal{U}$ 
for a given vector field $G \in \mathcal{X}^r(N)$, then the vector field $F + \varepsilon G$ has
a locally invariant manifold $M_\varepsilon $ within an $\varepsilon $-neighborhood of $M$.
It is diffeomorphic to $M$ and has the same stability as that of $M$.
\\[-0.2cm]

Further, Fenichel \cite{Fen3, Fen4} proved that $M_\varepsilon $ admits a fibration:
there exists a family of smooth manifolds $\{\mathcal{F}_\varepsilon (p)\}_{p\in M_\varepsilon }$ such that
\\[0.2cm]
(i) if $p \neq p'$, then $\mathcal{F}_\varepsilon (p) \cap \mathcal{F}_\varepsilon (p') = \emptyset$.
\\
(ii) $\mathcal{F}_\varepsilon (p) \cap M_\varepsilon = \{p\}$.
\\
(iii) the family $\{\mathcal{F}_\varepsilon (p) \}$ is invariant in the sense that $\phi _t (\mathcal{F}_\varepsilon (p))
\subset \mathcal{F}_\varepsilon (\phi _t(p))$, where $\phi _t$ is a flow generated by $F + \varepsilon G \in \mathcal{X}^r(N)$.
\\
(iv) there exist $C >0, \lambda >0$ such that for $q\in \mathcal{F}_\varepsilon (p)$, $|| \phi _t(p) - \phi _t(q) || < Ce^{-\lambda t}$,
where we suppose for simplicity that $M$ (and thus $M_\varepsilon $) is attracting.
\\[0.2cm]
See also Wiggins \cite{Wig2} for Fenichel theory.
These theorems are applied to fast-slow systems by Fenichel \cite{Fen} to obtain a slow manifold $M_\varepsilon $
and a flow around $M_\varepsilon $.
Roughly speaking, these theorems state that for a fast-slow system, there is a locally invariant manifold $M_\varepsilon $,
called the slow manifold, within an $\varepsilon $-neighborhood of the critical manifold $M$ if $\varepsilon >0$ is 
sufficiently small. A flow near $M_\varepsilon $ is given as the sum of the slow motion (dynamics on $M_\varepsilon $)
and the fast motion.
If $M_\varepsilon $ is attracting, the fast motion decays exponentially to zero and eventually a flow is well approximated
by the dynamics on $M_\varepsilon $.

Applying these results to our fast-slow system (\ref{4-0}), when $\delta $ is independent of $\varepsilon $, we obtain an
attracting slow manifold $M_\varepsilon $ and we can construct an approximate flow around $M_\varepsilon $.
However, if $\delta $ depends on $\varepsilon $, Fenichel theory is no longer applicable in general even if $\varepsilon <\!< \delta $.
To see this, let us recall how the existence of $M_\varepsilon $ is proved.

For simplicity of exposition, suppose that vector fields are defined on $\mathbf{R}^m \times \mathbf{R}^n$.
We denote a point on this space as $(x,z) \in \mathbf{R}^m \times \mathbf{R}^n$.
Suppose that a given unperturbed vector field $F$ has an attracting compact normally hyperbolic invariant manifold $M$
on the subspace $\{x = 0\}$. We denote a flow $\phi _t$ generated by the perturbed vector field $F + \varepsilon G$ by
\begin{eqnarray*}
\phi _t(x,z,\varepsilon ) = (\phi _t^1(x,z,\varepsilon ), \phi _t^2(x,z,\varepsilon )).
\end{eqnarray*}
From the assumption of normal hyperbolicity, we can show that there exists a positive constant $T$ such that
\begin{equation}
\Bigl|\Bigl| \frac{\partial \phi _T^1}{\partial x}(0,z,0) \Bigl|\Bigl| 
      \cdot \Bigl|\Bigl| \frac{\partial \phi ^2_T}{\partial z}(0,z,0)^{-1} \Bigl|\Bigl| < \frac{1}{4},
\quad \mathrm{for}\,\, (0,z) \in M,
\label{4-1}
\end{equation}
because $\partial \phi _t^1/\partial x$ decays faster than $\partial \phi _t^2/ \partial z$.
Since $M \subset \{x=0\}$ is $F$-invariant, we have
\begin{equation}
\phi _T^1(0,z,0) = 0, \quad \frac{\partial \phi _T^1}{\partial z}(0,z,0) = 0,
\quad \mathrm{for}\,\, (0,z) \in M.
\label{4-2}
\end{equation}
Since the flow is continuous with respect to $x,z$ and $\varepsilon $, for given small positive numbers $\eta_1$ and $\eta_2$,
there exist $\varepsilon _0 >0$ and an open set $V \supset M$ such that the inequalities
\begin{eqnarray}
& & \Bigl|\Bigl| \frac{\partial \phi _T^1}{\partial x}(x,z,\varepsilon ) \Bigl|\Bigl| 
      \cdot \Bigl|\Bigl| \frac{\partial \phi ^2_T}{\partial z}(x,z,\varepsilon )^{-1} \Bigl|\Bigl| < \frac{1}{2},\label{4-3} \\[0.1cm]
& & || \phi _T^1(x,z,\varepsilon )  || < \eta_1, \label{4-4} \\[0.1cm]
& & \Bigl|\Bigl| \frac{\partial \phi ^1_T}{\partial z}(x,z,\varepsilon ) \Bigl|\Bigl| < \eta_2, \label{4-5}
\end{eqnarray}
hold for $0<\varepsilon <\varepsilon _0$ and $(x,z) \in V$.
Let $S$ be the set of Lipschitz functions from $M$ into the $x$-space with a suitable norm.
Let $S_C$ be the subset of $S$ consisting of functions $h$ such that $(h(z),z)\in V$ and 
their Lipschitz constants are smaller than some constant $C>0$.
We now define the map $G : S_c \to S$ through
\begin{eqnarray*}
(Gh)(\phi _T^2(h(z),z,\varepsilon )) = \phi _T^1(h(z),z,\varepsilon ).
\end{eqnarray*}
By using inequalities (\ref{4-3}, \ref{4-4}, \ref{4-5}) (and several inequalities which trivially follow from compactness of $M$),
we can show that $G$ is a contraction map from $S_C$ into $S_C$. 
See Lemma 3.2.9 of Wiggins \cite{Wig2}, in which all inequalities for proving Fenichel's theorem are collected.
Thus $G$ has a fixed point $h_\varepsilon $ satisfying
$h_\varepsilon (\phi _T^2(h_\varepsilon (z),z,\varepsilon )) = \phi _T^1(h_\varepsilon (z),z,\varepsilon )$.
This proves that the graph of $x = h_\varepsilon (z)$, which defines $M_\varepsilon $, is invariant under the flow 
$\phi _t(\, \cdot\, ,\, \cdot \,, \varepsilon )$.
The existence of a fibration $\{\mathcal{F}_\varepsilon (p)\}_{p\in M_\varepsilon }$ can be proved in a similar manner.

If the unperturbed vector field $F = F_\delta $ smoothly depends on $\delta $ and if $\delta $ depends on $\varepsilon $,
the above discussion is not valid even if $\varepsilon <\!<\delta $.
The inequality (\ref{4-1}) for $F_\delta $ does not imply the inequality (\ref{4-3}) for $F_\delta + \varepsilon G$ in general.
For example, consider the linear system $\dot{\bm{x}} = A_0 \bm{x} + \delta A_1 \bm{x}$ with matrices
\begin{eqnarray*}
A_0 = \left(
\begin{array}{@{\,}cc@{\,}}
0 & 1 \\
0 & 0
\end{array}
\right), \quad A_1 = \left(
\begin{array}{@{\,}cc@{\,}}
-1 & 0 \\
0 & -1
\end{array}
\right).
\end{eqnarray*}
Suppose that $\delta = \sqrt{\varepsilon }$.
Eigenvalues of $A_0 + \delta A_1$ are given by $-\delta $ (double root), so that the derivative of the flow at the origin
is exponentially small for $t>0$.
Next, add the perturbation $\varepsilon A_2 \bm{x}$ to this system, where
\begin{eqnarray*}
A_2 = \left(
\begin{array}{@{\,}cc@{\,}}
0 & 0 \\
4 & 0
\end{array}
\right).
\end{eqnarray*}
Although $\varepsilon A_2$ is quite smaller than $A_0 + \delta A_1$ if $\varepsilon $ is sufficiently small,
the eigenvalues of $A_0 + \delta A_1 + \varepsilon A_2$ are $\delta $ and $-3 \delta $, so that 
the derivative of the flow of the perturbed system diverges as $t\to \infty$.
This shows that Eq.(\ref{4-1}) does not imply Eq.(\ref{4-3}) in general if $\delta $ depends on $\varepsilon $.
Further, the open set $V$ above also depends on $\varepsilon $ through $\delta $ and it may shrink as $\varepsilon \to 0$.
For this linear system, it is easy to see that such a stability change does not occur if $A_0$ has no Jordan block.
For our fast-slow system, the assumption (C5) allows us to prove that such a stability change does not occur.
\\[0.2cm]
\textbf{Lemma 4.1.} Let $A(\delta , z)$ and $B(\varepsilon , \delta , z)$ be $2 \times 2$ matrices
which are $C^\infty$ in their arguments.
Suppose that eigenvalues of $A(\delta , z)$ are given by $-\delta \mu (z, \delta ) \pm \sqrt{-1} \omega (z, \delta )$
with the conditions $\mu (z, \delta ) > 0$ and $\omega (z, \delta ) \neq 0$ for $\delta \geq 0$.
Further suppose that $\delta $ depends on $\varepsilon $ as $\varepsilon \sim o(\delta )$
\, (that is, $\varepsilon < \!< \delta $ as $\varepsilon \to 0$).
Then, eigenvalues of $A(\delta , z) + \varepsilon B(\varepsilon , \delta , z)$ are given by
\begin{equation}
-\delta \mu (z, \delta ) \pm \sqrt{-1} \omega (z, \delta ) + O(\varepsilon )
\label{4-6}
\end{equation} 
as $\varepsilon \to 0$.
\\[0.2cm]
\textbf{Proof.} Straightforward calculation. \hfill $\blacksquare$
\\[-0.2cm]

Now we return to our fast-slow system (\ref{4-0}).
Put $\bm{X} = (X,Y), \bm{f} = (f_1, f_2)$ and rewrite Eq.(\ref{4-0}) as
\begin{equation}
\dot{\bm{X}} = \bm{f}(\bm{X}, Z, \varepsilon , \delta ), \quad \dot{Z} = \varepsilon g(\bm{X}, Z, \varepsilon , \delta ).
\label{4-7}
\end{equation}
The flow generated by this system is denoted as
\begin{equation}
\phi _t(\bm{X},Z,\varepsilon ,\delta ) = (\phi _t^1(\bm{X},Z,\varepsilon ,\delta ) , \phi _t^2(\bm{X},Z,\varepsilon ,\delta )).
\end{equation}
Recall that $S_a^+(z_2, z_1)$ is expressed as $X = 0, Y = -\eta$;
that is, $\bm{f}(0, -\eta, Z, 0, \delta ) = 0$ for $z_2 \leq Z \leq z_1$.
When $\varepsilon =0$, $\phi _t^2 (\bm{X},Z,0 ,\delta ) = Z$, which proves that 
$|| (\partial \phi _t^2 (\bm{X},Z,0 ,\delta )/ \partial Z)^{-1} || = 1$.
Next, the derivative of $\phi _t^1$ satisfies the variational equation
\begin{eqnarray*}
\frac{d}{dt}\frac{\partial \phi _t^1}{\partial \bm{X}}(\bm{X},Z, 0,\delta )
 = \frac{\partial \bm{f}}{\partial \bm{X}}(\phi _t^1(\bm{X},Z, 0,\delta ),Z,0,\delta  ) 
      \frac{\partial \phi _t^1}{\partial \bm{X}}(\bm{X},Z, 0,\delta ).
\end{eqnarray*}
On $S_a^+(z_2, z_1)$, this is reduced to the autonomous system
\begin{eqnarray*}
\frac{d}{dt}\frac{\partial \phi _t^1}{\partial \bm{X}}(0, -\eta, Z, 0,\delta )
 = \frac{\partial \bm{f}}{\partial \bm{X}}(0, -\eta ,Z, 0,\delta  ) 
      \frac{\partial \phi _t^1}{\partial \bm{X}}(0, -\eta, Z, 0,\delta ).
\end{eqnarray*}
The assumption (C5) implies that the eigenvalues of the matrix 
$\displaystyle \frac{\partial \bm{f}}{\partial \bm{X}}(0, -\eta,Z,0,\delta )$ are given by 
$-\delta \mu^+ (z, \delta ) \pm \sqrt{-1} \omega^+ (z, \delta )$.
Thus 
\begin{eqnarray*}
\frac{\partial \phi _t^1}{\partial \bm{X}}(0, -\eta, Z, 0,\delta ) \sim O(e^{-\delta t})
\end{eqnarray*}
on $S_a^+(z_2, z_1)$.
This proves the inequality
\begin{equation}
\Bigl|\Bigl| \frac{\partial \phi _T^1}{\partial \bm{X}}(0, -\eta,Z,0, \delta ) \Bigl|\Bigl| 
      \cdot \Bigl|\Bigl| \frac{\partial \phi ^2_T}{\partial Z}(0, -\eta,Z,0, \delta )^{-1} \Bigl|\Bigl| < \frac{1}{4},
\label{4-8}
\end{equation}
for some large $T>0$.
In general, this does not imply Eq.(\ref{4-3}) as was explained. 
However, in our situation, by applying Lemma 4.1 to 
\begin{eqnarray*}
A(\delta , z) = \frac{\partial \bm{f}}{\partial \bm{X}}(0, -\eta,Z,0,\delta ),
\end{eqnarray*}
it turns out that eigenvalues of the matrix 
$\displaystyle \frac{\partial \bm{f}}{\partial \bm{X}}(0, -\eta,Z,\varepsilon ,\delta )$ are of the form (\ref{4-6})
for small $\varepsilon >0$.
Therefore, $\displaystyle \frac{\partial \phi _t^1}{\partial \bm{X}}(0, -\eta, Z, \varepsilon ,\delta )$
also decays with the rate $O(e^{-\delta t})$ on $S_a^+(z_2, z_1)$.
Further, the assumption (C6) proves that there exists a neighborhood $V^+$ of $S_a^+(z_2, z_1)$,
which is independent of $\delta $, such that real parts of eigenvalues of 
$\partial \bm{f}/\partial \bm{X}$ are also of order $O(-\delta )$ on $V^+$.
This yields the inequality (\ref{4-3}) on $V^+$.
Inequalities (\ref{4-4}) and (\ref{4-5}) are easily obtained.
In this manner, all inequalities for proving Fenichel's theorem are obtained, and the existence of the slow manifold $M_\varepsilon $
and a fibration on $M_\varepsilon $ for our system are proved in the standard way as long as $\varepsilon <\!< \delta $
(To prove Theorem 3, we will suppose that $\delta \sim O(\varepsilon (-\log \varepsilon )^{1/2}) > \!> \varepsilon $).
Note that the existence of a neighborhood $V^+$ of the critical manifold, on which eigenvalues of $\partial \bm{f}/\partial \bm{X}$
have negative real parts, are also assumed in the classical approach for singular perturbed problems to estimate the dynamics
of fast motion, see O'Malley \cite{Oma} and Smith \cite{Smi}.

We have seen that a solution of (\ref{4-7}) on $V^+$ is written as the sum of the slow motion on the slow manifold
and the fast motion which decays exponentially. 
To calculate them, it is convenient to introduce the slow time scale by $\tau = \varepsilon t$,
which provides
\begin{equation}
\varepsilon \frac{d\bm{X}}{d\tau} = \bm{f}(\bm{X}, Z, \varepsilon , \delta ), 
     \quad \frac{dZ}{d\tau} = g(\bm{X}, Z, \varepsilon , \delta ).
\label{4-9}
\end{equation}
A solution of this system is given by
\begin{equation}
\left\{ \begin{array}{l}
\bm{X}(\tau, \varepsilon , \delta ) = \bm{x}_s(\tau, \varepsilon , \delta ) + \bm{x}_f(\tau, \varepsilon , \delta ),  \\
Z(\tau, \varepsilon , \delta )  = z_s(\tau, \varepsilon ,\delta ) + z_f(\tau , \varepsilon , \delta),  \\
\end{array} \right.
\label{4-10}
\end{equation}
where $\bm{x}_s, z_s$ describe the slow motion and $\bm{x}_f, z_f$ describe the fast motion.
They are $C^\infty$ in $\varepsilon $ (see Fenichel \cite{Fen}) and their expansions with respect to $\varepsilon $
are obtained step by step according to O'Malley \cite{Oma} as follows:
We expand them as
\begin{eqnarray*}
& & \bm{x}_s(\tau, \varepsilon , \delta ) = \sum^\infty_{k=0} \varepsilon ^k \bm{x}_s^{(k)}(\tau, \delta ), \quad
\bm{x}_f(\tau, \varepsilon , \delta ) = \sum^\infty_{k=0} \varepsilon ^k \bm{x}_f^{(k)}(\tau, \delta ), \\
& & z_s(\tau, \varepsilon , \delta ) = \sum^\infty_{k=0} \varepsilon ^k z_s^{(k)}(\tau, \delta ), \quad
z_f(\tau, \varepsilon , \delta ) = \sum^\infty_{k=0} \varepsilon ^k z_f^{(k)}(\tau, \delta ),
\end{eqnarray*}
with the initial condition
\begin{eqnarray*}
\bm{X}(0, \varepsilon , \delta ) = \bm{x}_0(\delta ) + O(\varepsilon ), \quad
Z(0, \varepsilon , \delta ) = z_0(\delta ) + O(\varepsilon ),
\end{eqnarray*}
in $V^+$.
At first, $\bm{x}_s^{(0)}$ and $z_s^{(0)}$ are determined to satisfy the system (\ref{4-9}) for $\varepsilon =0$.
Thus $\bm{x}^{(0)}_s$ is given by $\bm{x}^{(0)}_s = (0, -\eta)$ and $z^{(0)}_s$ is given as the solution of the equation
\begin{equation}
\frac{dz^{(0)}_s}{d\tau} = g(0, -\eta , z^{(0)}_s, 0, \delta )
\label{4-10b}
\end{equation}
with the initial condition $z^{(0)}_s(0, \delta ) = z_0(\delta )$.
This system is called the slow system.
Next, from the system (\ref{4-7}) for $\varepsilon =0$, we obtain $z^{(0)}_f \equiv 0$,
and $\bm{x}_f^{(0)}$ is governed by the system
\begin{equation}
\frac{d\bm{x}_f^{(0)}}{dt} = \frac{d\bm{X}}{dt}(t,0,\delta ) - \frac{d\bm{x}_s^{(0)}}{dt}(\tau, \delta )
= \bm{f}((0, -\eta) + \bm{x}_f^{(0)}, z_s^{(0)}(\tau) , 0, \delta )
\label{4-10c}
\end{equation}
with the initial condition
\begin{equation}
\bm{x}_f^{(0)}(0, \delta ) = \bm{x}_0(\delta ) - \bm{x}_s(0, 0, \delta ) = \bm{x}_0(\delta ) - (0, -\eta).
\end{equation}
Fenichel's theorem (Part (iv) above) shows that if $\bm{x}_f^{(0)}(0, \delta ) \in V^+$,
then $\bm{x}_f^{(0)}$ decays exponentially as $t\to \infty$.
In the classical approach \cite{Oma}, the existence of $V^+$ is used to estimate Eq.(\ref{4-10c}) directly
to prove that $\bm{x}_f^{(0)}$ decays exponentially, see also Smith \cite{Smi}.
To investigate behavior of a solution as $\varepsilon \to 0$, we rewrite Eq.(\ref{4-10c}) as
\begin{equation}
\frac{d\bm{x}_f^{(0)}}{d\tau}
 = \frac{1}{\varepsilon }\frac{\partial \bm{f}}{\partial \bm{X}}(0, -\eta , z_s^{(0)}(\tau ), 0, \delta )\bm{x}_f^{(0)}
 + \frac{1}{\varepsilon }Q_1(\bm{x}_f^{(0)},\delta ),
\label{4-10d}
\end{equation}
where $Q_1 \sim O((\bm{x}_f^{(0)})^2)$ is a $C^\infty$ function.
\\[0.2cm]
\textbf{Lemma 4.2.} A solution of the system (\ref{4-10d}) is given by
\begin{eqnarray}
& & \left(
\begin{array}{@{\,}cc@{\,}}
K_1(\tau, \varepsilon ) \cos \Bigl[ \frac{1}{\varepsilon } W(\tau) \Bigr]
 + K_2(\tau, \varepsilon ) \sin \Bigl[ \frac{1}{\varepsilon } W(\tau) \Bigr] 
 & K_3(\tau, \varepsilon ) \cos \Bigl[ \frac{1}{\varepsilon } W(\tau) \Bigr]
 + K_4(\tau, \varepsilon ) \sin \Bigl[ \frac{1}{\varepsilon } W(\tau) \Bigr]  \\[0.2cm]
K_5(\tau, \varepsilon ) \cos \Bigl[ \frac{1}{\varepsilon } W(\tau) \Bigr]
 + K_6(\tau, \varepsilon ) \sin \Bigl[ \frac{1}{\varepsilon } W(\tau) \Bigr]  
 & K_7(\tau, \varepsilon ) \cos \Bigl[ \frac{1}{\varepsilon } W(\tau) \Bigr]
 + K_8(\tau, \varepsilon ) \sin \Bigl[ \frac{1}{\varepsilon } W(\tau) \Bigr] 
\end{array}
\right) \times \nonumber \\[0.2cm]
& & \quad \exp \Bigl[ - \frac{\delta }{\varepsilon }\int^\tau_{0} \! \mu^+ (z^{(0)}_s (s), \delta ) ds \Bigr] 
\bigl( \bm{x}_f^{(0)}(0, \delta ) + \bm{u}(\tau, \varepsilon , \delta ; \bm{x}_f^{(0)}(0, \delta )) \bigr) ,
\label{4-11}
\end{eqnarray}
where $W(\tau) = \int^\tau_{0} \! \omega ^+ (z^{(0)}_s (s), \delta ) ds$, $K_i\, (i = 1, \cdots , 8)$ are $C^\infty$ functions,
and $\bm{u} \sim O(\bm{x}_f^{(0)}(0, \delta )^2)$
denotes higher order terms with respect to the initial value.
\\[0.2cm]
\textbf{Proof.} We use the WKB analysis.
Put $\bm{x}_f^{(0)} = (v_1, v_2)$ and 
\begin{equation}
\frac{\partial \bm{f}}{\partial \bm{X}}(0, -\eta , z_s^{(0)}(\tau ), 0, \delta ) =
\left(
\begin{array}{@{\,}cc@{\,}}
a(\tau) & b(\tau) \\
c(\tau) & d(\tau)
\end{array}
\right).
\label{4-11c}
\end{equation}
Let us consider the linearized system
\begin{equation}
\frac{d}{d\tau}\left(
\begin{array}{@{\,}c@{\,}}
v_1 \\
v_2
\end{array}
\right) = \frac{\partial \bm{f}}{\partial \bm{X}}(0, -\eta , z_s^{(0)}(\tau ), 0, \delta ) \left(
\begin{array}{@{\,}c@{\,}}
v_1 \\
v_2
\end{array}
\right) = \left(
\begin{array}{@{\,}cc@{\,}}
a(\tau) & b(\tau) \\
c(\tau) & d(\tau)
\end{array}
\right) \left(
\begin{array}{@{\,}c@{\,}}
v_1 \\
v_2
\end{array}
\right).
\label{4-11d}
\end{equation}
Then, $v_1(\tau)$ proves to satisfy the equation
\begin{equation}
\varepsilon ^2 v_1'' - \left( \varepsilon (a+d) + \varepsilon ^2 \frac{b'}{b} \right) v_1'
 + \left( ad - bc + \varepsilon (\frac{ab'}{b}-a') \right) v_1 = 0.
\label{4-11b}
\end{equation}
We construct a formal solution of the form
\begin{eqnarray*}
v_1(\tau) = \exp \Bigl[ \frac{1}{\varepsilon } \sum^\infty_{n=0} \varepsilon ^n S_n(\tau) \Bigr].
\end{eqnarray*}
Substituting it into Eq.(\ref{4-11b}), we obtain an equation of $S_0$
\begin{eqnarray*}
(S_0')^2 - (a + d) S_0' + (ad - bc) =0.
\end{eqnarray*}
This is solved as
\begin{eqnarray*}
S_0(\tau)  = \int^\tau_{0} \! \lambda _+(s) ds, \,\,  \int^\tau_{0} \! \lambda _-(s) ds,
\end{eqnarray*}
where
\begin{eqnarray*}
\lambda _{\pm}(\tau) = -\delta \mu (z^{(0)}_s (\tau), \delta ) \pm \sqrt{-1} \omega ^+ (z^{(0)}_s (\tau), \delta )
\end{eqnarray*}
are eigenvalues of the matrix (\ref{4-11c}).
For each $\int^\tau_{0} \! \lambda _+(s) ds$ and $\int^\tau_{0} \! \lambda _-(s) ds$,
$S_1, S_2, \cdots $ are uniquely determined. 
Thus a general solution $v_1(\tau)$ is given by
\begin{eqnarray*}
v_1(\tau ) = C_1 \exp \Bigl[ \frac{1}{\varepsilon } \int^\tau_{0} \! \lambda _+(s) ds \Bigr] K_{11}(\tau, \varepsilon )
 + C_2 \exp \Bigl[ \frac{1}{\varepsilon } \int^\tau_{0} \! \lambda _-(s) ds \Bigr] K_{12}(\tau, \varepsilon ),
\end{eqnarray*}
where $C_1, C_2\in \mathbf{C}$ and $K_{11}, K_{12}$ are $C^\infty$ functions.
In a similar manner, it turns out that $v_2$ is expressed as
\begin{eqnarray*}
v_2(\tau ) = C_1 \exp \Bigl[ \frac{1}{\varepsilon } \int^\tau_{0} \! \lambda _+(s) ds \Bigr] K_{21}(\tau, \varepsilon )
 + C_2 \exp \Bigl[ \frac{1}{\varepsilon } \int^\tau_{0} \! \lambda _-(s) ds \Bigr] K_{22}(\tau, \varepsilon ).
\end{eqnarray*}
Therefore, a general solution of the system (\ref{4-11d}) is written as
\begin{eqnarray*}
\left(
\begin{array}{@{\,}c@{\,}}
v_1 \\
v_2
\end{array}
\right) = \left(
\begin{array}{@{\,}cc@{\,}}
\exp \Bigl[ \frac{1}{\varepsilon } \int^\tau_{0} \! \lambda _+(s) ds \Bigr] K_{11}(\tau, \varepsilon ) 
 & \exp \Bigl[ \frac{1}{\varepsilon } \int^\tau_{0} \! \lambda _-(s) ds \Bigr] K_{12}(\tau, \varepsilon ) \\[0.2cm]
\exp \Bigl[ \frac{1}{\varepsilon } \int^\tau_{0} \! \lambda _+(s) ds \Bigr] K_{21}(\tau, \varepsilon )
 & \exp \Bigl[ \frac{1}{\varepsilon } \int^\tau_{0} \! \lambda _-(s) ds \Bigr] K_{22}(\tau, \varepsilon )
\end{array}
\right) \left(
\begin{array}{@{\,}c@{\,}}
C_1 \\
C_2
\end{array}
\right) .
\end{eqnarray*}
The fundamental matrix of (\ref{4-11d}) is given by
\begin{eqnarray*}
\left(
\begin{array}{@{\,}cc@{\,}}
\exp \Bigl[ \frac{1}{\varepsilon } \int^\tau_{0} \! \lambda _+(s) ds \Bigr] K_{11}(\tau, \varepsilon ) 
 & \exp \Bigl[ \frac{1}{\varepsilon } \int^\tau_{0} \! \lambda _-(s) ds \Bigr] K_{12}(\tau, \varepsilon ) \\[0.2cm]
\exp \Bigl[ \frac{1}{\varepsilon } \int^\tau_{0} \! \lambda _+(s) ds \Bigr] K_{21}(\tau, \varepsilon )
 & \exp \Bigl[ \frac{1}{\varepsilon } \int^\tau_{0} \! \lambda _-(s) ds \Bigr] K_{22}(\tau, \varepsilon )
\end{array}
\right) \left(
\begin{array}{@{\,}cc@{\,}}
K_{11}(0, \varepsilon ) & K_{12}(0, \varepsilon ) \\[0.2cm]
K_{21}(0, \varepsilon ) & K_{22}(0, \varepsilon )
\end{array}
\right) ^{-1}.
\end{eqnarray*}
This shows that each component of the fundamental matrix is a linear combination of
\begin{eqnarray*}
\exp \Bigl[ - \frac{\delta }{\varepsilon }\int^\tau_{0} \! \mu^+ (z^{(0)}_s (s), \delta ) ds \Bigr] 
\cos \Bigl[ \frac{1}{\varepsilon } W(\tau) \Bigr]
\quad \mathrm{and} \quad \exp \Bigl[ - \frac{\delta }{\varepsilon }\int^\tau_{0} \! \mu^+ (z^{(0)}_s (s), \delta ) ds \Bigr] 
\sin \Bigl[ \frac{1}{\varepsilon } W(\tau) \Bigr].
\end{eqnarray*}
Finally, the variation-of-constants formula is applied to the nonlinear system (\ref{4-10d}) to prove Lemma 4.2.
\hfill $\blacksquare$
\\[-0.2cm]

With this $\bm{x}_f^{(0)}$, the zeroth order approximate solution is constructed as
\begin{eqnarray}
\left(
\begin{array}{@{\,}c@{\,}}
\bm{X}(\tau, \varepsilon ,\delta ) \\
Z(\tau, \varepsilon ,\delta )
\end{array}
\right) = \left(
\begin{array}{@{\,}c@{\,}}
O(\varepsilon ) \\
-\eta + O(\varepsilon )\\
z_s^{(0)}(\tau, \delta ) + O(\varepsilon ),
\end{array}
\right) + \left(
\begin{array}{@{\,}c@{\,}}
\bm{x}_f^{(0)} (\tau, \delta ) + O(\varepsilon ) \\
O(\varepsilon )
\end{array}
\right),
\label{4-12}
\end{eqnarray}
as long as the orbit is in $V^+$.
The first term in the right hand side denotes the position on $M_\varepsilon $
and the second term denotes the deviation from $M_\varepsilon $.
It is known that all terms $\bm{x}_f^{(k)}, z_f^{(k)}$ in the expansions of the fast motion decay exponentially
as well as $\bm{x}_f^{(0)}$ (\cite{Fen, Oma, Smi}).

Combining this approximate solution near the slow manifold with the transition map near the fold point,
Theorem 1 is easily proved.
\\[0.2cm]
\textbf{Proof of Theorem 1.} To prove Theorem 1, $\delta $ is assumed to be fixed.
For the system (\ref{local1}), take an initial value in $V^+$.
Then, a solution is given by (\ref{4-12}) with (\ref{4-11}).
These expressions show that when $t>0$, the solution lies sufficiently close to the critical manifold $S^+_a$ 
if $\varepsilon $ is sufficiently small. Because of the assumption (A3), $z_s$ decreases (where we suppose that
$S^+$ is convex downward) with the velocity of order $\varepsilon $
(with respect to the original time scale $t$).
Thus the solution reaches the section $\Sigma^+_{in}$ after some time, which is of order $O(1/\varepsilon )$.
The intersection point is mapped into $\Sigma^+_{out}$ by the transition map $\Pi^+_{loc}$ given in Thm.3.2,
and it proves that after passing through $\Sigma^+_{out}$ the distance between the solution and the orbit $\alpha ^+$
is of order $O(\varepsilon ^{4/5})$. \hfill $\blacksquare$


\subsection{Global Poincar\'{e} map}

In Sec.3, the transition map $\Pi^+_{loc}$ around the fold point $L^+(\delta )$ had been constructed.
The transition map around the fold point $L^-(\delta )$ is obtained in the same way.
The sections $\Sigma^-_{in}$ and $\Sigma^-_{out}$ are defined in a similar way to 
$\Sigma^+_{in}$ and $\Sigma^+_{out}$ (see Fig.\ref{fig4}),
respectively, and the transition map $\Pi^-_{loc}$ from an open set in $\Sigma^-_{in}$ into $\Sigma^-_{out}$
along the flow of (\ref{4-0}) proves to take the same form as $\Pi^+_{loc}$, although functions $G_1, G_2$ and 
higher order terms denoted as $O(\varepsilon \log \varepsilon )$ may be different from one another
(note that $\Omega $ and $H$ are common for $\Pi^+_{loc}$ and $\Pi^-_{loc}$ because they arise from the first 
Painlev\'{e} equation).

\begin{figure}[h]
\begin{center}
\includegraphics[scale=1.0]{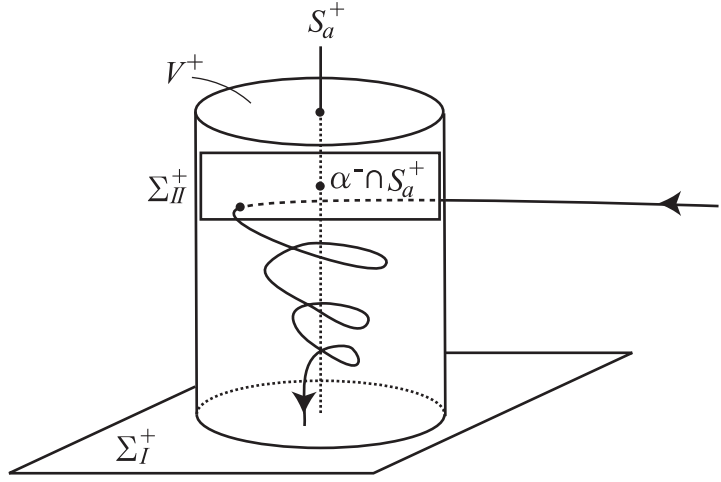}
\end{center}
\caption{The sections $\Sigma ^+_{I}, \Sigma ^+_{II}$ and an orbit of Eq.(\ref{4-0}).}
\label{fig12}
\end{figure}

Since the unperturbed system has a heteroclinic orbit $\alpha ^-$ connecting $L^- (\delta )$ with a point on $S^+_a (\delta )$
and since $S^+_a (\delta )$ has an attraction basin $V^+$ which is independent of $\delta $, there is an open set 
$U^-_{out}\subset \Sigma^-_{out} $, which is independent of $\delta $ and $\varepsilon $, such that orbits of (\ref{4-0})
starting from $U^-_{out}$ go into $V^+$ and are eventually approximated by Eq.(\ref{4-12}).
Let $z_0$ be the $Z$ coordinate of $L^-(\delta )$.
Define the section $\Sigma^+_{II}$ to be
\begin{equation}
\Sigma^+_{II} = V^+ \cap \{(X,Y,Z) \, | \, Y = -\eta ,\, |Z- z_0| \leq \rho_4 \},
\end{equation}
where $\rho_4$ is a small positive number so that a solution of (\ref{4-0}) starting from $U^-_{out}$
intersects $\Sigma^+_{II}$ only once (see Fig.\ref{fig12}).

The global Poincar\'{e} map is constructed as follows:
Let $\Pi^+_{II, out},\, \Pi^+_{I, II},\, \Pi^+_{in, I}$ be transition maps from $U^-_{out} \subset \Sigma ^-_{out}$
into $\Sigma ^+_{II}$, $\Sigma ^+_{II}$ into $\Sigma ^+_{I}$, $\Sigma ^+_{I}$ into $\Sigma ^+_{in}$, respectively.
Then, the transition map $\Pi^+$ from $U^-_{out}$ into $\Sigma ^+_{out}$ is given by
\begin{eqnarray*}
\Pi^+ = \Pi^+_{loc} \circ \Pi^+_{in, I} \circ \Pi^+_{I, II} \circ \Pi^+_{II, out}.
\end{eqnarray*}
The transition map $\Pi^-$ from an open set in $\Sigma^+_{out} $ into $\Sigma^-_{in}$ is calculated in a similar manner and 
it has the same form as $\Pi^+$.
The global Poincar\'{e} map is given  by $\Pi^+ \circ \Pi^-$. 
However, it is sufficient to investigate one of them by identifying $\Sigma^+_{out} $ and $\Sigma^-_{out} $.
If $\Pi^+ : U^-_{out} \to \Sigma^+_{out} $ is a contraction map, so is $\Pi^+ \circ \Pi^-$,
and if $\Pi^+$ has a horseshoe, so is $\Pi^+ \circ \Pi^-$ because $\Pi^+$ and $\Pi^-$ have the same properties.
To identify two sections $\Sigma ^-_{out}$ and $\Sigma ^+_{out}$, recall that
$L^- = (0, y_0, z_0)$ in the $(X,Y,Z)$-coordinate, and define $\Sigma ^-_{out}$ to be $\{ Y = y_0 - \rho_1^2\}$.
Let $U^-_{out}$ be an open set in $\Sigma ^-_{out}$ such that the transition map $\Pi^+_{II, out} : U^-_{out} \to \Sigma ^+_{II}$
is well-defined.
The set $U^-_{out}$ includes the point $\Sigma ^-_{out} \cap \alpha ^-$.
We identify $U^-_{out}$ with an open set $U^+_{out}$ in $\Sigma ^+_{out}$ by the translation
\begin{equation}
\mathcal{T} : \left(
\begin{array}{@{\,}c@{\,}}
X \\
\rho_1^2 \\
Z
\end{array}
\right) \mapsto
\left(
\begin{array}{@{\,}c@{\,}}
X \\
y_0- \rho_1^2 \\
Z + z_0
\end{array}
\right).
\end{equation}
Then, the transition map $\tilde{\Pi}^+_{II,out}$ from $U^+_{out}\subset \Sigma ^+_{out}$ into $\Sigma ^+_{II}$ 
is obtained by combining the translation and $\Pi^+_{II, out}$.
Since the velocity in the $Z$ direction is of order $\varepsilon $, it is expressed as
\begin{equation}
\tilde{\Pi}^+_{II,out}\left(
\begin{array}{@{\,}c@{\,}}
X \\
\rho_1^2 \\
Z
\end{array}
\right) = 
\Pi^+_{II, out} \circ \mathcal{T}\left(
\begin{array}{@{\,}c@{\,}}
X \\
\rho_1^2 \\
Z
\end{array}
\right) = \left(
\begin{array}{@{\,}c@{\,}}
P^+(X,Z,\varepsilon ,\delta ) \\
-\eta \\
Z + z_0 + O(\varepsilon )
\end{array}
\right),
\end{equation}
where $P^+$ is a $C^\infty$ function. 
Since $\tilde{\Pi}^+_{II,out}$ is $C^\infty$, we expand it as
\begin{equation}
\tilde{\Pi}^+_{II,out}\left(
\begin{array}{@{\,}c@{\,}}
X \\
\rho_1^2 \\
Z
\end{array}
\right) = 
\left(
\begin{array}{@{\,}c@{\,}}
p(\delta ) + O(X,Z, \varepsilon ) \\
-\eta \\
Z + z_0 + O(\varepsilon )
\end{array}
\right),
\end{equation}
To prove Theorem 3, we will use the fact that there exists a positive constant $p_0 > 0$ such that
$|p(\delta )| > p_0$ for $0 < \delta  < \delta _0$, which is proved as follows:
Since $\delta $ controls the strength of the stability of $S^+_a$, if $\delta $ is sufficiently small,
orbits which converge to $(0, -\eta, z_0)$\, (the intersection of the heteroclinic orbit $\alpha ^-$ and $S^+_a$)
rotate around this point so many times. 
In particular, they intersect with $\Sigma^+_{II}$ before reaching $(0, -\eta, z_0)$.
If $p(\delta )$ were zero, the right hand side above tends to $(0, -\eta, z_0)$ as $X,Z,\varepsilon \to 0$,
which yields a contradiction.

Next thing to do is to combine the above $\tilde{\Pi}^+_{II,out}$ with $\Pi^+_{I, II}$.
By Eq.(\ref{4-12}), the transition map $\Pi^+_{I, II}$ from $\Sigma^+_{II}$ into $\Sigma^+_{I}$ is given by
\begin{equation}
\Pi^+_{I, II}\left(
\begin{array}{@{\,}c@{\,}}
X \\
-\eta \\
Z
\end{array}
\right) = \left(
\begin{array}{@{\,}c@{\,}}
O(\varepsilon ) \\
-\eta + O(\varepsilon )\\
z_2
\end{array}
\right) + \left(
\begin{array}{@{\,}c@{\,}}
\bm{x}_f^{(0)}(\tau (X,Z, \varepsilon , \delta ), \delta ) + O(\varepsilon )\\
0
\end{array}
\right),
\end{equation}
where $\bm{x}_f^{(0)} = \bm{x}_f^{(0)}(\tau, \delta )$ is given by (\ref{4-11})
with the initial condition $\bm{x}_f^{(0)}(0, \delta ) = (X, 0)$, 
$z_2 = \rho_1^4 + e^{-1/\varepsilon ^2}$ is the $Z$ coordinate of the section $\Sigma^+_{I}$ as defined before, 
and $\tau = \tau (X,Z, \varepsilon , \delta )$ is a transition time (with respect to the slow time scale)
from a point $(X,-\eta , Z)$ to $\Sigma^+_{I}$.
This transition time $\tau$ is determined as follows:
Let $z_s^{(0)}(\tau, \delta )$ be a solution of Eq.(\ref{4-10b}) with the initial condition $z^{(0)}_s(0, \delta ) = Z$.
Then, Eq.(\ref{4-12}) implies that $\tau = \tau (X,Z, \varepsilon , \delta )$ is given as a root of the equation
\begin{eqnarray*}
z_2 = z^{(0)}_s(\tau, \delta ) + O(\varepsilon ).
\end{eqnarray*}
Let $\hat{\tau}$ be a root of the equation $z_2 = z^{(0)}_s(\tau, \delta ) $.
By virtue of the implicit function theorem, $\tau$ is written as $\tau = \hat{\tau} + O(\varepsilon )$.
Since Eq.(\ref{4-10b}) is independent of $X$ and $\varepsilon $, so is $\hat{\tau}$.
Thus we obtain
\begin{equation}
\tau (X, Z, \varepsilon , \delta ) = \hat{\tau} (Z, \delta ) + O(\varepsilon ).
\end{equation}
Further, $\hat{\tau}$ is bounded as $\delta \to 0$ because $g\neq 0$ on $S^+_a$ uniformly in $0 \leq \delta < \delta _0$.
Therefore, $\Pi^+_{I, II}$ proves to be of the form
\begin{eqnarray}
& & \Pi^+_{I, II}\left(
\begin{array}{@{\,}c@{\,}}
X \\
-\eta \\
Z
\end{array}
\right) = \left(
\begin{array}{@{\,}c@{\,}}
O(\varepsilon ) \\
-\eta + O(\varepsilon ) \\
z_2
\end{array}
\right) \nonumber \\
& & + \left(
\begin{array}{@{\,}c@{\,}}
X \left( K_1(\hat{\tau}, \varepsilon ) \cos \Bigl[ \frac{1}{\varepsilon } W(\hat{\tau}) \Bigr]
 + K_2(\hat{\tau}, \varepsilon ) \sin \Bigl[ \frac{1}{\varepsilon } W(\hat{\tau}) \Bigr]  \right)
  \exp \Bigl[ - \frac{ \delta}{\varepsilon }\int^{\hat{\tau}}_{0} \! \mu^+ (z^{(0)}_s (s), \delta ) ds \Bigr] (1 + O(\varepsilon ,X ))\\[0.2cm]
X \left( K_5(\hat{\tau}, \varepsilon ) \cos \Bigl[ \frac{1}{\varepsilon } W(\hat{\tau}) \Bigr]
 + K_6(\hat{\tau}, \varepsilon ) \sin \Bigl[ \frac{1}{\varepsilon } W(\hat{\tau}) \Bigr]  \right)
  \exp \Bigl[ - \frac{ \delta}{\varepsilon }\int^{\hat{\tau}}_{0} \!\mu^+ (z^{(0)}_s (s), \delta ) ds \Bigr] (1 + O(\varepsilon ,X))\\[0.2cm]
0
\end{array}
\right).\nonumber \\
\label{4-13b}
\end{eqnarray}
The first line denotes the intersection point $M_\varepsilon  \cap \Sigma^+_{I}$ and thus it is independent of $X$ and $Z$.
The second line denotes the deviation from the intersection.
Note that the transition map $\Pi^+_{in, I}$ from $\Sigma ^+_{I}$ into $\Sigma ^+_{in}$ is $O(e^{-1/\varepsilon ^2})$-close
to the identity map.
Thus $\Pi^+_{in, I}\circ \Pi^+_{I,II} \circ \Pi^+_{II, out} \circ \mathcal{T}$ is calculated as
\begin{eqnarray}
& & \Pi^+_{in, I}\circ \Pi^+_{I,II} \circ \Pi^+_{II, out} \circ \mathcal{T} \left(
\begin{array}{@{\,}c@{\,}}
X \\
\rho_1^2 \\
Z
\end{array}
\right) = \left(
\begin{array}{@{\,}c@{\,}}
O(\varepsilon ) \\
-\eta + O(\varepsilon )\\
\rho_1^4
\end{array}
\right)  \nonumber \\
& +&  \left(
\begin{array}{@{\,}c@{\,}}
p(\delta ) \left( K_1(\hat{\tau}, \varepsilon ) \cos \Bigl[ \frac{1}{\varepsilon } W(\hat{\tau}) \Bigr]
 + K_2(\hat{\tau}, \varepsilon ) \sin \Bigl[ \frac{1}{\varepsilon } W(\hat{\tau}) \Bigr]  \right)
     \exp \Bigl[ - \frac{ \delta}{\varepsilon }
          \int^{\hat{\tau}}_{0} \! \mu^+ (z^{(0)}_s (s), \delta ) ds \Bigr] (1 + O(\varepsilon ,X , Z )) \\[0.2cm]
p(\delta ) \left( K_5(\hat{\tau}, \varepsilon ) \cos \Bigl[ \frac{1}{\varepsilon } W(\hat{\tau}) \Bigr]
 + K_6(\hat{\tau}, \varepsilon ) \sin \Bigl[ \frac{1}{\varepsilon } W(\hat{\tau}) \Bigr] \right)
     \exp \Bigl[ - \frac{ \delta}{\varepsilon }
          \int^{\hat{\tau}}_{0} \!  \mu^+ (z^{(0)}_s (s), \delta ) ds \Bigr](1 + O(\varepsilon ,X , Z )) \\[0.2cm]
0
\end{array}
\right) , \nonumber \\
\label{4-13}
\end{eqnarray}
where $\hat{\tau} = \hat{\tau}(Z + z_0, \delta )$ and $z^{(0)}_s(\tau)$ is a solution of (\ref{4-10b})
satisfying the initial condition $z^{(0)}_s (0) = Z + z_0$.
Finally, the transition map 
\begin{eqnarray*}
\Pi^+ = \Pi^+_{loc} \circ \Pi^+_{in, I}\circ \Pi^+_{I,II} \circ \Pi^+_{II, out} \circ \mathcal{T}
\end{eqnarray*}
from $U^-_{out}$ into $\Sigma^+_{out} $ is obtained by combining the above map with $\Pi^+_{loc}$.
\\[-0.2cm]

At this stage, we can prove Theorem 2.
\\
\textbf{Proof of Theorem 2.} To prove Theorem 2, it is sufficient to show that the map $\Pi^+$
has a hyperbolically stable fixed point. Then, the global Poincar\'{e} map (without identifying 
$\Sigma ^+_{out}$ and $\Sigma ^-_{out}$) has the same property
because $\Pi^-$ takes the same form as $\Pi^+$.
Indeed, if $\varepsilon $ is sufficiently small for fixed $\delta $, Them.3.2 and Eq.(\ref{4-13}) show that
the image of the map $\Pi^+$ is exponentially small, and thus $\Pi^+$ is a contraction map.
Further, eigenvalues of the derivative of $\Pi^+$ is of order $O(e^{-1/\varepsilon })$,
which proves that $\Pi^+$ has a hyperbolically stable fixed point. \hfill $\blacksquare$


\subsection{Derivative of the transition map}

If $\delta $ is fixed, it is obvious that the transition map $\Pi^+$ is of order $O(e^{-1/\varepsilon })$
as $\varepsilon \to 0$.
However, when $\delta $ is small as well as $\varepsilon $, the action of $\Pi^+$ becomes more complex.
In what follows, we suppose that $\delta $ depends on $\varepsilon $ and $\varepsilon  \sim o(\delta )
\, (\varepsilon <\! < \delta )$ as $\varepsilon \to 0$.
A straightforward calculation shows that the derivative of $\Pi^+$ is of the form
\begin{eqnarray}
\frac{\partial \Pi^+}{\partial (X,Z)} &=&
\left(
\begin{array}{@{\,}cc@{\,}}
L_1(X,Z,\varepsilon ,\delta )\varepsilon ^{1/5} & L_2(X,Z,\varepsilon ,\delta )\varepsilon ^{-4/5} \\
L_3(X,Z,\varepsilon ,\delta )\varepsilon ^{1/5} & L_4(X,Z,\varepsilon ,\delta )\varepsilon ^{-4/5}
\end{array}
\right) \times \nonumber \\
 & & \quad \exp \Bigl[ - \hat{d}(\rho, \delta ) \frac{\delta }{\varepsilon }\Bigr] \cdot
\exp \Bigl[ -\frac{\delta }{\varepsilon } \int^{\hat{\tau}}_{0} \! \mu^+(z^{(0)}_s(s), \delta )ds \Bigr] 
(1 + L_5(X,Z,\varepsilon ,\delta )),
\label{4-14}
\end{eqnarray}
where $L_i\, (i = 1, \cdots  ,4)$ are bounded as $\varepsilon \to 0$, and $L_5$ denotes higher order terms
such that $L_5 \sim o(1)$ as $X,Z,\varepsilon \to 0$.

Eigenvalues of the derivative are given by
\begin{equation}
\lambda _1 = L_4 \varepsilon ^{-4/5} \exp \Bigl[ - \hat{d}(\rho, \delta ) \frac{\delta }{\varepsilon }\Bigr] \cdot
\exp \Bigl[ -\frac{\delta }{\varepsilon } \int^{\hat{\tau}}_{0} \! \mu^+(z^{(0)}_s(s), \delta )ds \Bigr] (1 + o(1)),
\end{equation}
and 
\begin{equation}
\lambda _2 = \frac{L_1L_4-L_2L_3}{L_4} \varepsilon ^{1/5} \exp \Bigl[ - \hat{d}(\rho, \delta ) \frac{\delta }{\varepsilon }\Bigr] \cdot
\exp \Bigl[ -\frac{\delta }{\varepsilon }\int^{\hat{\tau}}_{0} \! \mu^+(z^{(0)}_s(s), \delta )ds  \Bigr] (1 + o(1)).
\end{equation}
If $\delta $ is fixed, they are exponentially small as $\varepsilon \to 0$, although if $\delta $ is small as well as
$\varepsilon $, $|\lambda _1|$ may become large.
For example, if $\delta  = C \varepsilon (-\log \varepsilon )^{1/2}$ with a positive constant $C$,
and if $L_4(X,Z,\varepsilon ,\delta )\neq 0$, $|\lambda _1|$ is of order $\varepsilon ^{-4/5}e^{-C(-\log \varepsilon )^{1/2}}$,
which is larger than $1$ if $\varepsilon $ is sufficiently small.
On the other hand, $|\lambda _2|$ is always smaller than $1$.
The function $L_4$ is given by
\begin{eqnarray}
L_4(X,Z, \varepsilon , \delta )
&=& \frac{\partial H}{\partial X}
(\hat{D}_1 \varepsilon ^{-3/5}e^{-\hat{d}\delta /\varepsilon }, \hat{D}_2 \varepsilon ^{-2/5}e^{-\hat{d}\delta /\varepsilon }) 
\cdot \frac{\partial \hat{D}_1}{\partial X} \cdot p(\delta ) \cdot \frac{\partial }{\partial Z}W(\hat{\tau}) \times \nonumber\\
& & \left( -K_1 (\hat{\tau}, \varepsilon ) \sin \Bigl[ \frac{1}{\varepsilon } W(\hat{\tau}) \Bigr]
 + K_2(\hat{\tau}, \varepsilon ) \cos \Bigl[ \frac{1}{\varepsilon } W(\hat{\tau}) \Bigr]\right),
\end{eqnarray} 
in which arguments of $\hat{D}_i = \hat{D}_i(\,\cdot \,, \,\cdot \,, \rho_1, \varepsilon ,\delta )$ are given by the first
and second components of Eq.(\ref{4-13}).
From Thm.3.2 (III) and (IV), we obtain $\partial H/ \partial X \neq 0,\, \partial \hat{D}_1/ \partial X\neq 0$.
The value $p(\delta )$ is also not zero as was explained above.
Recall that $\hat{\tau}(Z+z_0, \delta )$ is defined as a transition time along the flow of Eq.(\ref{4-10b}).
Since $g<0$ uniformly on $S^+_a$ and $0\leq \delta < \delta _0$, $\hat{\tau}$ is monotonically increasing with respect to $Z$.
Further, $W(\hat{\tau})$ is monotonically decreasing or monotonically increasing because $\omega ^+ \neq 0$ uniformly.
This proves $\partial W(\hat{\tau})/\partial Z \neq 0$.
Therefore, $L_4 = 0$ if and only if
\begin{eqnarray*}
& & -K_1 (\hat{\tau}(Z+z_0, \delta ), \varepsilon ) \sin \Bigl[ \frac{1}{\varepsilon } W(\hat{\tau}(Z+z_0, \delta )) \Bigr]
 + K_2(\hat{\tau}(Z+z_0, \delta ), \varepsilon ) \cos \Bigl[ \frac{1}{\varepsilon } W(\hat{\tau}(Z+z_0, \delta )) \Bigr]\\
&=& -K_1 (\hat{\tau}(z_0, \delta ), \varepsilon ) \sin \Bigl[ \frac{1}{\varepsilon } W(\hat{\tau}(Z+z_0, \delta )) \Bigr]
 + K_2(\hat{\tau}(z_0, \delta ), \varepsilon ) \cos \Bigl[ \frac{1}{\varepsilon } W(\hat{\tau}(Z+z_0, \delta )) \Bigr] + O(Z)
\end{eqnarray*}
is zero.
If there exists $Z$ such that the above value is zero, then it is zero for a countable set of values of $Z$ because of the 
periodicity.
For these ``bad" $Z$, $\lambda _1$ degenerates and $|\lambda _1|$ may become smaller than $1$.
Now we have the same situation as the proof of the existence of chaos in Silnikov's systems.
In the proof of Silnikov's chaos, an eigenvalue of a transition map degenerates if and only if an expression
$k_1 \sin (\log (z/\varepsilon )) + k_2 \cos (\log (z/\varepsilon ))$ is zero, where $k_1$ and $k_2$ are some constants,
see Wiggins \cite{Wig}.


\subsection{Proof of Theorem 3}

Now we are in a position to prove Theorem 3.
The proof is done in the same way as the proof of Silnikov's chaos.
At first, we show that the transition map $\Pi^+$ has a topological horseshoe:
We show that an image of a rectangle under $\Pi^+$ becomes a ring-shaped area and it appropriately intersects with 
the rectangle.
Next, to prove that the horseshoe is hyperbolic, we investigate the derivative of $\Pi^+$.
We can avoid ``bad" $Z$, at which the derivative degenerates, because they are at most countable.
\\[0.2cm]
\textbf{Proof of Thm.3.}\,\, Suppose that $\delta  = C_1 \varepsilon (-\log \varepsilon )^{1/2}$ with some positive constant $C_1$.
Recall that there exists a slow manifold within an $\varepsilon $ neighborhood of $S^+_a$.
Since it is one dimension, the slow manifold is a solution orbit of the system (\ref{4-7}).
By virtue of Thm.3.2, this orbit intersects with $\Sigma ^+_{out}$ near $\alpha ^+$.
Let $Q\in \Sigma ^+_{out}$ be the intersection point of this orbit and $\Sigma ^+_{out}$.
Take a rectangle $R$ on $\Sigma ^+_{out}$ including the point $Q$,
whose boundaries are parallel to the $X$ axis and the $Z$ axis (see Fig.\ref{fig4}).
Let $h_R = C_2\varepsilon $ be the height of $R$, where $C_2$ is a positive constant to be determined.
The image of $R$ under the map $\tilde{\Pi}^+_{out, II} = \Pi^+_{out, II} \circ \mathcal{T}$
is a deformed rectangle whose ``height" is also of order $O(\varepsilon )$ since $dZ/dt \sim O(\varepsilon )$.

\begin{figure}[h]
\begin{center}
\includegraphics[scale=1.0]{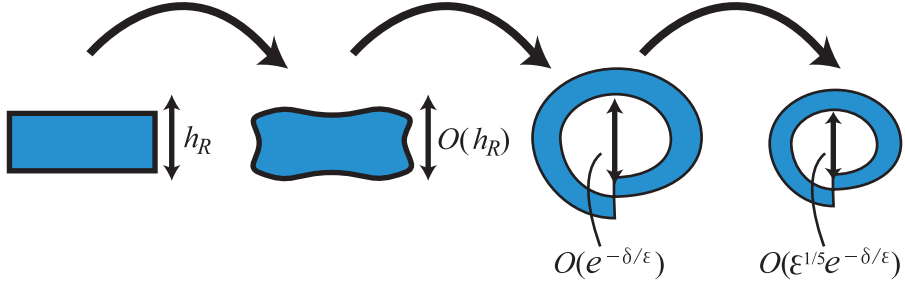}
\end{center}
\caption{Images of the rectangle $R$ under a succession of transition maps.}
\label{fig13}
\end{figure}

Next thing to consider is the shape of $\Pi^+_{II,I} \circ \tilde{\Pi}^+_{out, II}(R)$.
It is easy to show by using Eq.(\ref{4-13b}) that the image of $\tilde{\Pi}^+_{out, II}(R)$ under the map $\Pi^+_{II,I}$
becomes a ring-shaped area whose radius is of order $e^{-\delta /\varepsilon }$.
Since the ``height" of $\tilde{\Pi}^+_{out, II}(R)$ is of order $\varepsilon $, the rotation angle of the ring-shaped area is estimated as
\begin{equation}
\frac{1}{\varepsilon }W(\hat{\tau}(Z+z_0 + O(\varepsilon ))) - \frac{1}{\varepsilon }W(\hat{\tau}(Z+z_0 )) 
 = \frac{1}{\varepsilon }\int^{\hat{\tau}(Z+z_0 + O(\varepsilon ))}_{\hat{\tau}(Z+z_0 )} \! \omega ^+(z^{(0)}_s(s),\delta )ds \sim O(1).
\label{4-16}
\end{equation}
Thus we can choose $C_2$ so that the rotation angle of the ring-shaped area is sufficiently close
to $2\pi$ as is shown in Fig.\ref{fig13}.

Finally, we consider the shape of $\Pi^+ (R)$ by using Thm.3.2.
Since $\partial H/ \partial \mathcal{X}(0,0) \neq 0$,
the expansion of $H$ is estimated as
\begin{equation}
H(\mathcal{X}, \mathcal{Y}) 
   \sim \mathcal{X} \varepsilon ^{-3/5} \exp [-\hat{d}\delta /\varepsilon ] (1 + O(\varepsilon ^{1/5})). 
\label{4-17}
\end{equation}
This and Eq.(\ref{3-7}) show that
the radius of $\Pi^+(R)$ is of order $O(\varepsilon ^{1/5}e^{-\delta /\varepsilon })$.
Since we put $\delta  = C_1 \varepsilon (-\log \varepsilon )^{1/2}$, the inequality
\begin{equation}
h_R = C_2 \varepsilon  < \! < O(\varepsilon ^{1/5}e^{- \delta /\varepsilon })
\label{4-18}
\end{equation}
holds if $\varepsilon $ is sufficiently small.
Further, the ring $\Pi^+(R)$ surrounds the point $Q$ because the image of the rectangle $R$ under the flow
rotates around the slow manifold when passing between the section $\Sigma ^{+}_{II}$ and $\Sigma ^{+}_{in}$.
This means that two horizontal boundaries of $R$ intersect with the ring $\Pi^+(R)$ as is shown in Fig.\ref{fig5} (b).
It is obvious that the vertical boundaries of $R$ are mapped to the inner and outer boundaries of the ring,
and the horizontal boundaries are mapped to the other boundaries in radial direction.
This proves that the map $\Pi^+$ creates a horseshoe and thus has an invariant Cantor set.

To prove that this invariant set is hyperbolic, it is sufficient to show that there exist two disjoint rectangles $H_1$ and $H_2$ in $R$,
whose horizontal boundaries are parallel to the $X$ axis and vertical boundaries are included in those of $R$, such that
the inequalities 
\begin{eqnarray}
& & || D_x\Pi^+_1 || < 1, \\
& & || (D_z\Pi^+_2)^{-1} || < 1, \\
& & 1 - || (D_z\Pi^+_2)^{-1} || \cdot || D_x\Pi_1^+ ||
   > 2 \sqrt{|| D_z\Pi^+_1 || \cdot || D_x\Pi^+_2|| \cdot || (D_z\Pi^+_2)^{-1} ||^2}, \\
& & 1 - (|| D_x \Pi_1^+ ||+|| (D_z\Pi^+_2)^{-1} ||) + || D_x\Pi^+_1 || \cdot || (D_z\Pi^+_2)^{-1} ||
 > || D_x\Pi^+_2 || \cdot || D_z\Pi^+_1 || \cdot || (D_z\Pi^+_2)^{-1} ||,\quad \quad
\end{eqnarray}
hold on $H_1 \cup H_2$, where $\Pi^+_1$ and $\Pi^+_2$ denote the $X$ and $Z$ components of $\Pi^+$, respectively,
and $D_x$ and $D_z$ denote the derivatives with respect to $X$ and $Z$, respectively.
See Wiggins \cite{Wig} for the proof.
We can take such $H_1$ and $H_2$ so that ``bad" $Z$, at which $L_4 = 0$, are not included.
Then, inequalities above immediately follows from Eq.(\ref{4-14}):
$|| D_x\Pi^+_1 ||$ and $|| D_x\Pi^+_2 ||$ are sufficiently small, and
$|| D_z\Pi^+_1 ||$ and $|| D_z\Pi^+_2 ||$ are sufficiently large as $\varepsilon \to 0$.
This proves Theorem 3.\hfill $\blacksquare$


\section{Concluding remarks}

Our assumption of Bogdanov-Takens type fold points is not generic in the sense that the Jacobian matrix
has two zero eigenvalues.
However, this assumption is not essential for existence of periodic orbits or chaotic invariant sets.

At first, we remark that Theorems 2 and 3 hold even if we add a small perturbation to Eq.(\ref{1-5}),
since hyperbolic invariant sets remain to exist under small perturbations.

\begin{figure}[h]
\begin{center}
\includegraphics[scale=1.0]{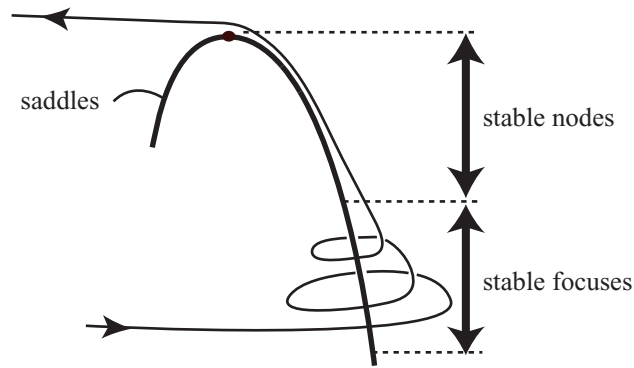}
\end{center}
\caption{Critical manifold consisting of a saddle-node type fold point, stable nodes, and stable focuses
and an orbit near it.}
\label{fig14}
\end{figure}

Second, we can consider the case that one of the connected components of critical manifolds consists of 
stable nodes, stable focuses and a saddle-node type fold point 
(\textit{i.e.} a saddle-node bifurcation point of a unperturbed system), as in Fig.\ref{fig14}.
In this case, Theorem 2 is proved in a similar way and Theorem 3 still holds if the length of the subset
of the critical manifold consisting of stable focuses is of order $O(1)$.
However, analysis of saddle-node type fold points is well performed in~\cite{Kru1, Mis, Gil} and thus we do not deal with
such a situation in this paper.

\begin{figure}[h]
\begin{center}
\includegraphics[scale=1.0]{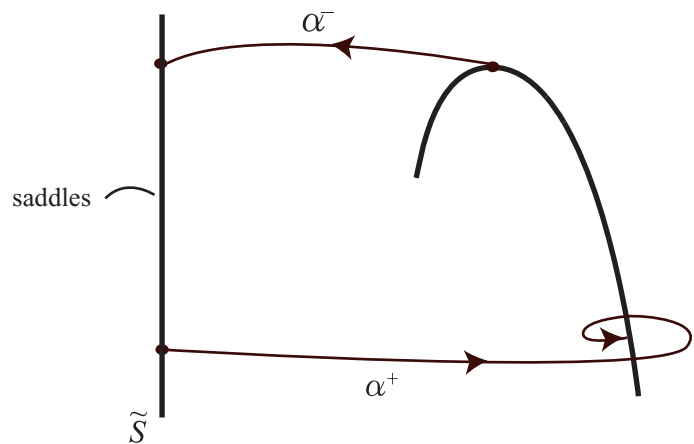}
\end{center}
\caption{Two connected components of critical manifolds.
One is similar to that of our system (\ref{1-5}) and the other consists of only saddles.}
\label{fig15}
\end{figure}

We can also consider the case that one of the connected components $\tilde{S}$ of critical manifolds has
no fold points but consists of saddles with heteroclinic orbits $\alpha ^{\pm}$, see Fig.\ref{fig15}.
In this case, analysis around the $\tilde{S}$ is done by using the exchange lemma (see Jones~\cite{Jon}) and 
we can prove theorems similar to Theorems 2 and 3.
Such a situation arises in an extended prey-predator system.
In~\cite{Kuw}, a periodic orbit and chaos in an extended prey-predator system are numerically investigated
with the aid of the theory of the present paper.

\vspace*{1.5cm}
\textbf{Acknowledgments}

The author would like to thank Professor Toshihiro Iwai and Professor Masataka Kuwamura for critical reading of the
manuscript and for useful comments.
\\


\end{document}